%% file: _main.tex
\pdfoutput=1 % for arxiv submission 
\documentclass[pdflatex,sn-mathphys-num]{sn-jnl-preprint}% Math and Physical Sciences Numbered Reference Style

\usepackage{graphicx}%
\usepackage{multirow}%
\usepackage{amsmath,amssymb,amsfonts}%
\usepackage{mathrsfs}%
\usepackage[title]{appendix}%
\usepackage[dvipsnames]{xcolor} % added dvipnames
\usepackage{textcomp}%
\usepackage{manyfoot}%
\usepackage{booktabs}%
\usepackage{listings}%
\input{./macros}

\usepackage{titlesec}
\newcommand{\mysecspace}{\hskip.5em}
\titlelabel{\thetitle.\mysecspace}

\colorlet{siaminlinkcolor}{green!50!black}
\colorlet{siamexlinkcolor}{red!50!black}
\colorlet{siamreviewcolor}{black!50}
\hypersetup{allcolors=siaminlinkcolor,urlcolor=siamexlinkcolor}

\theoremstyle{thmstyleone}%
\newtheorem{proposition}{Proposition}% to get separate numbers for theorem and proposition etc.

\theoremstyle{thmstyletwo}%
\newtheorem{example}{Example}%
\newtheorem{remark}{Remark}%

\theoremstyle{thmstylethree}%
\numberwithin{equation}{section}
\numberwithin{theorem}{section}

\ifpdf
\DeclareGraphicsExtensions{.eps,.pdf,.png,.jpg,.jpeg}
\else
\DeclareGraphicsExtensions{.eps}
\fi
\graphicspath{ {figures/} }

\ifpdf
\hypersetup{
	pdftitle={Operator learning meets inverse problems: A probabilistic perspective},
	pdfauthor={N. H. Nelsen and Y. Yang}
}
\fi

\begin{document}

\title[Operator learning meets inverse problems]{\phantom{A }Operator learning meets inverse problems:\\A probabilistic perspective}

\author*[1,2,3]{\fnm{Nicholas~H.} \sur{Nelsen}}\email{nnelsen@oden.utexas.edu}

\author*[2]{\fnm{Yunan} \sur{Yang}}\email{yunan.yang@cornell.edu}

\affil[1]{\normalsize \orgdiv{Department of Mathematics}, \orgname{Massachusetts Institute of Technology}, \orgaddress{\city{Cambridge}, \state{MA} \postcode{02139}, \country{USA}}}

\affil[2]{\normalsize \orgdiv{Department of Mathematics}, \orgname{Cornell University}, \orgaddress{\city{Ithaca}, \state{NY} \postcode{14853}, \country{USA}}}

\affil[3]{\normalsize \orgdiv{Oden Institute for Computational Engineering and Sciences \& Department of Aerospace Engineering and Engineering Mechanics}, \orgname{The University of Texas at Austin}, \orgaddress{\city{Austin}, \state{TX} \postcode{78712}, \country{USA}}}

\abstract{Operator learning offers a robust framework for approximating mappings between infinite-dimensional function spaces. It has also become a powerful tool for solving inverse problems in the computational sciences. This chapter surveys methodological and theoretical developments at the intersection of operator learning and inverse problems. It begins by summarizing the probabilistic and deterministic approaches to inverse problems, and pays special attention to emerging measure-centric formulations that treat observed data or unknown parameters as probability distributions. The discussion then turns to operator learning by covering essential components such as data generation, loss functions, and widely used architectures for representing function-to-function maps. The core of the chapter centers on the end-to-end inverse operator learning paradigm, which aims to directly map observed data to the solution of the inverse problem without requiring explicit knowledge of the forward map. It highlights the unique challenge that noise plays in this data-driven inversion setting, presents structure-aware architectures for both point predictions and posterior estimates, and surveys relevant theory for linear and nonlinear inverse problems. The chapter also discusses the estimation of priors and regularizers, where operator learning is used more selectively within classical inversion algorithms. 
}



\maketitle

\setcounter{tocdepth}{4}
\tableofcontents

\clearpage

%%%%%%%%%%%%%%%%%%%%%%%%%%%%%%%%%%%%%%%%%%%
\section{Introduction}\label{sec:introduction}
This chapter centers on operator learning for inverse problems with infinite-dimensional unknown parameters. This setting pervades numerous fields of science and engineering. The observed data can be finite-dimensional, which is realistic, or infinite-dimensional, which is a mathematical idealization useful for theoretical applications or when data is abundant. Mathematically, the relationship between the parameters and data is usually determined by (possibly nonlinear) partial- or integro-differential equations. For example, the parameter could be a spatially varying permeability coefficient function appearing in a physics model for groundwater flow. The measurements could be noisy observations of the associated pressure field at a sparse set of sensors in the domain. While all inverse problems are characterized by some notion of ill-posedness in the sense of Hadamard~\cite{hadamard1902problemes}, i.e., lacking existence, uniqueness, or continuity concerning the data, the primary classes of inverse problems that the present chapter surveys are those with inherent instability. This makes the numerical resolution of such problems challenging.

Classical methods for solving inverse problems often rely on hand-crafted regularization techniques, explicit models of the forward process, and iterative optimization routines. 
These methods are mathematically rigorous, well-established, and offer interpretability and theoretical guarantees. In practice, they may require substantial computation and careful modeling. Their performance can be sensitive to noise or to forward models that are noninvertible, ill-conditioned, or only partially specified.

In recent times, deep learning has demonstrated an exceptional ability to accurately solve inverse problems and even outperform classical reconstruction methods in linear and nonlinear imaging inverse problems \cite{lucas2018using,mccann2017convolutional,ongie2020deep}. The subsequent development and success of operator learning methods for the emulation of forward maps of partial differential equations (PDEs) has spurred significant interest in the application of such tools to inverse problems formulated in function spaces \cite{ovadia2024vito,cho2025physics,jiang2024resolution,jiao2024solving,kaltenbach2023semi,long2025invertible,gao2024adaptive,guo2021construct,guotransformer,wang2024latent,yang2021seismic,yang2023rapid,zhang2024bilo,zhang2025coefficient}.
Indeed, operator learning serves as an attractive alternative and a complement to classical inversion. Architectures such as neural operators are designed to map between infinite-dimensional spaces of functions. As a byproduct, the tunable parameters of these models are decoupled from the resolution of the discretized data used to train them. This property affords exceptional flexibility at inference time and opens up new directions for multiresolution regularization. 
This chapter focuses on applying operator learning methods through two complementary strategies: (\emph{i}) directly approximating the underlying inverse map from data alone, referred to as ``end-to-end learning,'' and (\emph{ii}) learning regularizers that respect known physical structure in the forward map, to then be incorporated into existing inverse solvers requiring regularization.

The chapter distinguishes between solving a single inverse problem with machine learning \cite{raissi2019physics,chen2021solving,jagtap2022physics,mishra2023estimates}, e.g., requiring the evolution of a Markov chain or running an optimization solver, and simultaneously solving or probing families of inverse problems with the same forward model but different observed data, which is enabled by operator learning methods.
There is an important and popular line of work that employs operator learning to build fast forward map surrogates for subsequent downstream use in various inverse problems \cite{burman2024stabilizing,cao2025derivative,herrmann2020deep,li2021fourier,li2024physics,lunz2021learned,wu2023large,zhou2024ai}. Under certain conditions, the error incurred by replacing the true forward operator with the surrogate operator controls the propagated error in the inverse problem solution \cite{helin2024intro,lie2018random,stuart2018posterior,cao2023residual,cao2024lazydino}; also see~\cite[Chp.~3]{bach2024inverse}. However, due to its emphasis on the forward map, such work is outside of the scope of the present chapter. Instead, the chapter focuses on the machine learning of quantities specific to inverse problems and their underlying mathematical structure.

Another area of research that is contemporary in data-driven inverse problems but beyond the scope of the chapter is generative modeling based on the dynamical transport of measure. This research direction is timely and ever-expanding, encompassing the rapid advances in score-based diffusion models~\cite{ho2020denoising,song2021score,sun2024provable,wu2024principled}, flow matching~\cite{lipman2022flow}, and stochastic interpolants~\cite{albergo2023building,albergo2023stochastic,chen2024probabilistic}. These frameworks enable sampling from posterior distributions, interpolating between distributions, or learning regularizers implicitly from data. They can offer scalable alternatives to classical sampling algorithms or variational inference. Instead, the present chapter restricts its attention to more established Bayesian and numerical analysis methodologies for these tasks. The main topics of the chapter are (\emph{i}) the learning of end-to-end inverse problem solvers directly from data and (\emph{ii}) the learning of priors and regularizers for use in traditional inverse problem solvers.

\begin{figure}[tb]
	\centering
 	\captionsetup{skip=10pt}
	\begin{tikzpicture}[
		node distance=0.5cm and 0.5cm,
		box/.style={rectangle, draw, rounded corners, align=center, minimum width=0.5cm, minimum height=0.5cm},
		arrow/.style={-{Stealth[]}, thick}, scale = 0.5
		]
		
		% Nodes
		\node[box] (inv) {Ill-Posed Inverse Problems\\(\Cref{sec:back_ip})};
		\node[box, below left=of inv] (prob) {Probabilistic \\Methods\\ (\Cref{sec:back_ip_prob})};
		\node[box, below right=of inv] (reg) {Regularization \\ Methods\\(\Cref{sec:back_ip_reg})};
		\node[box, below=1.5cm of inv] (op) {Operator \\ Learning\\(\Cref{sec:back_ol})};
		\node[box, below=2cm of prob] (sol) {Learning Inverse\\Problem Solvers \\(\Cref{sec:ee})};
		\node[box, below=2cm of reg] (prior) {Learning Prior \\ Distributions  \& Regularizers \\(\Cref{sec:reg})};
		
		% Connections
		\draw[arrow, line width=1.5pt] (inv) -- (prob);
		\draw[arrow, line width=1.5pt] (inv) -- (reg);
		\draw[densely dotted, arrow, line width=1.5pt] (prob) -- (sol);
		\draw[dashed, arrow, line width=1.5pt] (reg) -- (prior);
		\draw[densely dotted, arrow, line width=1.5pt] (op) -- (sol);
		\draw[dashed, arrow, line width=1.5pt] (op) -- (prior);
		\draw[densely dotted, arrow, bend left=25, line width=1.5pt] (reg) to (sol);
		\draw[dashed,arrow, bend right=25, line width=1.5pt] (prob) to (prior);
		
	\end{tikzpicture}
	\caption{Overall structure of the chapter. The dotted lines indicate the connections associated with end-to-end learning, whereas the dashed lines correspond to regularizer learning.}
	\label{fig:structure}
\end{figure}

\paragraph*{Outline}
This chapter reviews the current landscape of data-driven inverse problem solvers and regularizers based on operator learning. It places a particular focus on probabilistic formulations and theoretical insights.
Unless otherwise specified, the chapter works in the continuum. This means that all functions, operators, and measures are viewed as elements of infinite-dimensional metric spaces, and no appeal is made to a particular type of discretization. As illustrated in \cref{fig:structure}, the chapter is organized as follows. 

\Cref{sec:back_ip,sec:back_ol} provide the necessary preliminaries on inverse problems and operator learning, respectively. On the inverse problems side, \cref{sec:back_ip} discusses probabilistic and deterministic solution methods. It places particular emphasis on the expanding \emph{measure-centric} viewpoint, which casts certain inverse problems into ones where either the data or parameter are themselves probability distributions. Subsequently, \cref{sec:back_ol} reviews the fundamental components of the operator learning framework, including training data assumptions, loss functions, and popular architectures. 

The bulk of the chapter is concentrated in \cref{sec:ee}, which concerns the end-to-end operator learning approach that seeks to directly map observed data to the unknown parameters that solve the inverse problem. It reveals features unique to inverse map learning, describes several structure-exploiting operator learning architectures for both point prediction and posterior estimation, and surveys some theoretical guarantees.

\Cref{sec:reg} studies prior and regularizer learning from the perspectives of variational inference and denoising. This subject adopts a more specialized viewpoint than the previous section, focusing solely on machine learning to estimate sub-components within a larger analytical inversion procedure. The theoretical analysis discusses the stability properties of the prior-to-posterior map under prior perturbations and the convergence properties of learned Tikhonov regularizers. 

\Cref{sec:conclusion} concludes the chapter with an outlook toward challenges and future opportunities at the interface of operator learning and inverse problems.

%%%%%%%%%%%%%%%%%%%%%%%%%%%%%%%%%%%%%%%%%%%%%%%%%%%%%%%%%%%%%%%%%%%%%%%%%%%%%%%%%%%%%%%%
\section{Background on inverse problems}\label{sec:back_ip}
Inverse problems concern the recovery of unknown parameters from indirect, and often noisy, observations. Consider a data model of the form
\begin{align}\label{eqn:ip_main}
	y = \mathcal{G}(u) + \eta\,,
\end{align}
where $\mathcal{G} \colon \mathcal{U} \to \mathcal{Y}$ is a forward operator (also known as the parameter-to-data map), $u \in \mathcal{U}$ is the unknown parameter to be inferred, $y \in \mathcal{Y}$ is the observed data, and $\eta$ models measurement noise or model error. The inverse problem is to recover $u$ from $y$. In many applications arising from science and engineering, $\mathcal{G}$ is nonlinear. The parameter space $\cU$ is contained in an infinite-dimensional Banach space; this is referred to as a \emph{nonparametric} setting. The data space $\cY$ is high- or infinite-dimensional as well. The inherent nonlinearity and the high- or infinite-dimensional nature of many practical inverse problems can introduce significant mathematical and computational challenges.

More fundamentally, the inverse problem of recovering $u$ from $y$ is often \emph{ill-posed} in the sense of Hadamard: solutions may not exist, may not be unique, or may not depend continuously on the data \cite{hadamard1902problemes}. Regarding the last difficulty, the inverse problem solution's sensitivity to noise in $y$ makes na\"ive inversion schemes unstable and unreliable. To address this challenge, one may introduce additional information about the parameter $u$ in \eqref{eqn:ip_main}, such as smoothness, sparsity, or statistical structure. This can be done either through \emph{regularization} techniques or by adopting a \emph{probabilistic} framework to capture uncertainty and stabilize the solution. This section reviews approaches in both directions.

A canonical example of a severely ill-posed nonlinear inverse problem is \emph{electrical impedance tomography} (EIT) \cite{cheney1999electrical}.
EIT is a noninvasive imaging technique used to reconstruct the internal conductivity field of an object by making electrical measurements on the object's boundary.
Spatial variations in the conductivity correspond to different regions of interest in the object, such as malignant biological tissue in hospital patients or distinct Earth subsurface structures in a geophysical context.
Due to its well-understood mathematical properties, EIT has emerged as a prototypical nonlinear benchmark inverse problem for both traditional inverse solvers and deep learning methods \cite{agnelli2020classification,BERETTA2025114162,colibazzi2022learning,fan2020solving,guotransformer,guo2021construct,bui2022bridging,hamilton2018deep,hamilton2019beltrami,tanyu2023electrical,knudsen2009regularized}.
This chapter works with the idealized continuum measurement version of the problem, which is known as \emph{Calder\'on's problem} \cite{calderon2006inverse}. For convenience, we will use the Calder\'on problem and EIT terminology interchangeably. The following EIT example will appear frequently throughout the chapter.

\begin{figure}[tb]
	\centering
	\captionsetup{skip=10pt}
	\begin{subfigure}[]{0.334\textwidth}
		\centering
		\includegraphics[width=\textwidth]{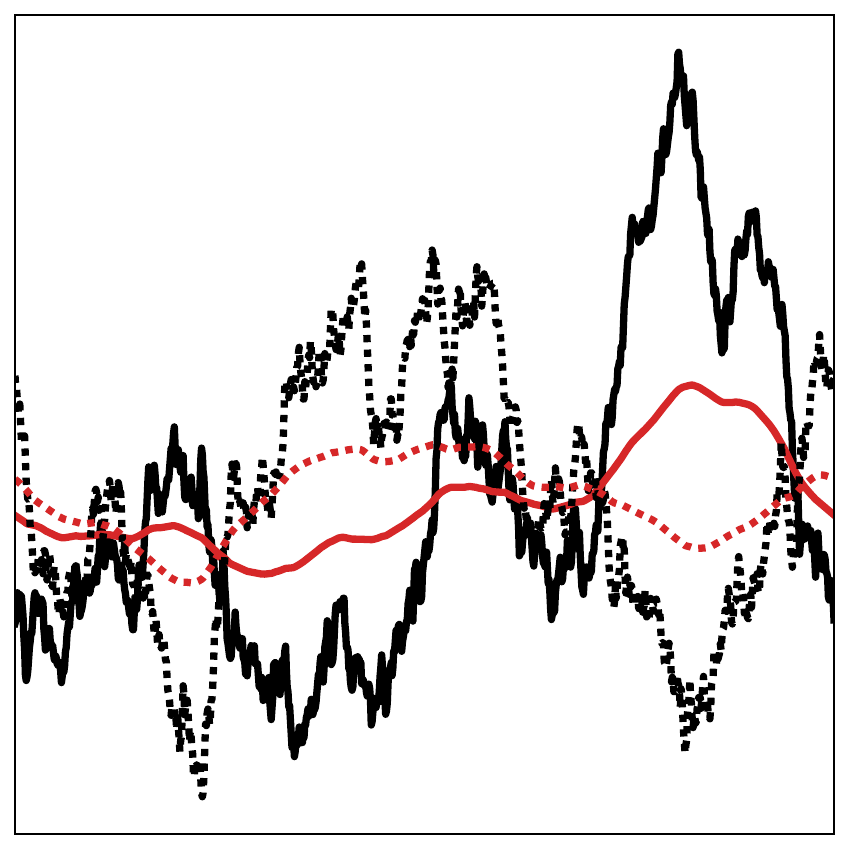}
		\caption{$\{(g_m, \Lambda_\gamma g_m)\}_{m=1}^2$}
		\label{subfig:eit_finite}
	\end{subfigure}
	\hfill%
	\begin{subfigure}[]{0.325\textwidth}
		\centering
		\includegraphics[width=\textwidth]{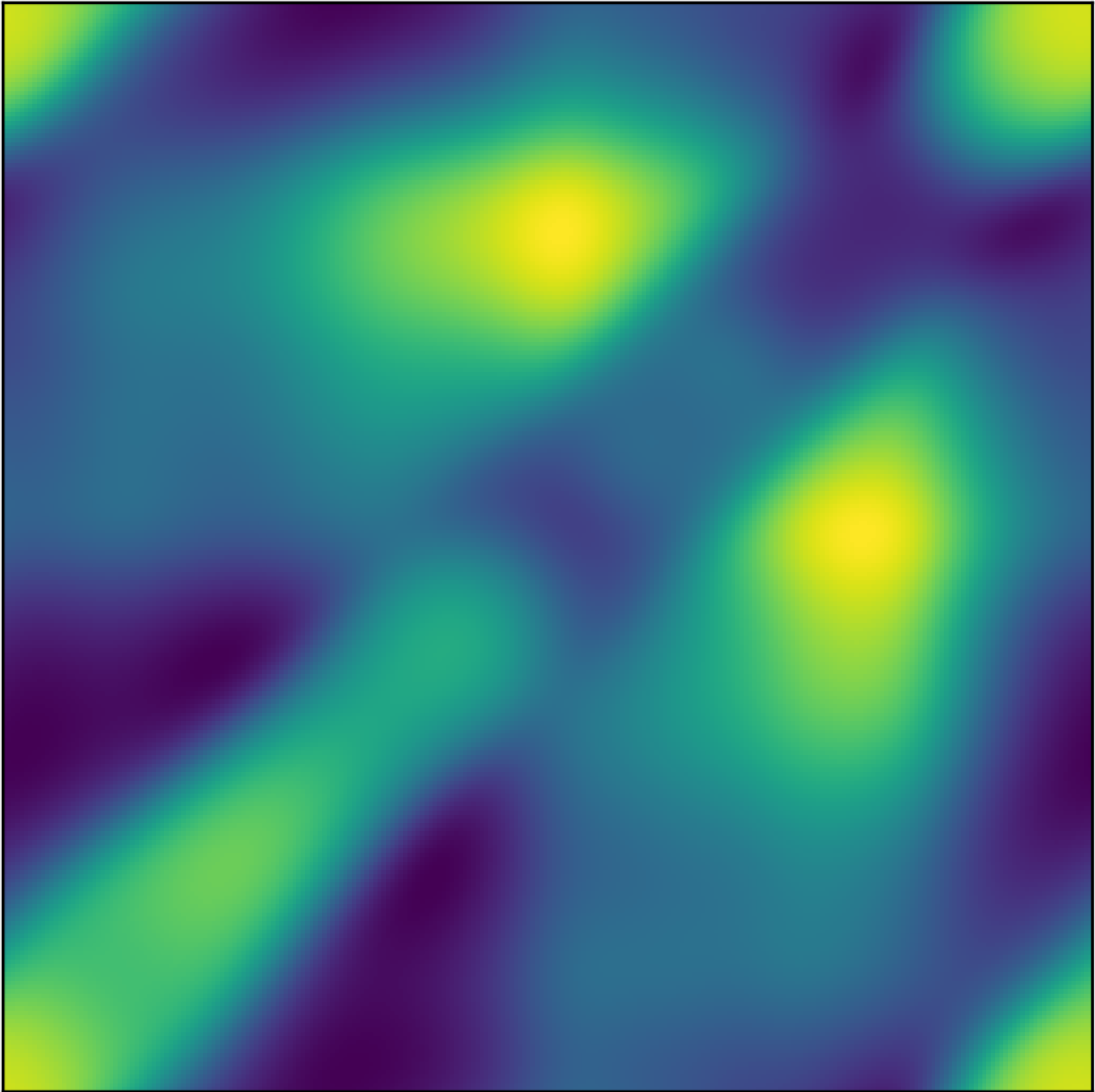}
		\caption{$\Lambda_\gamma$}
		\label{subfig:eit_in}
	\end{subfigure}
	\hfill%
	\begin{subfigure}[]{0.325\textwidth}
		\centering
		\includegraphics[width=\textwidth]{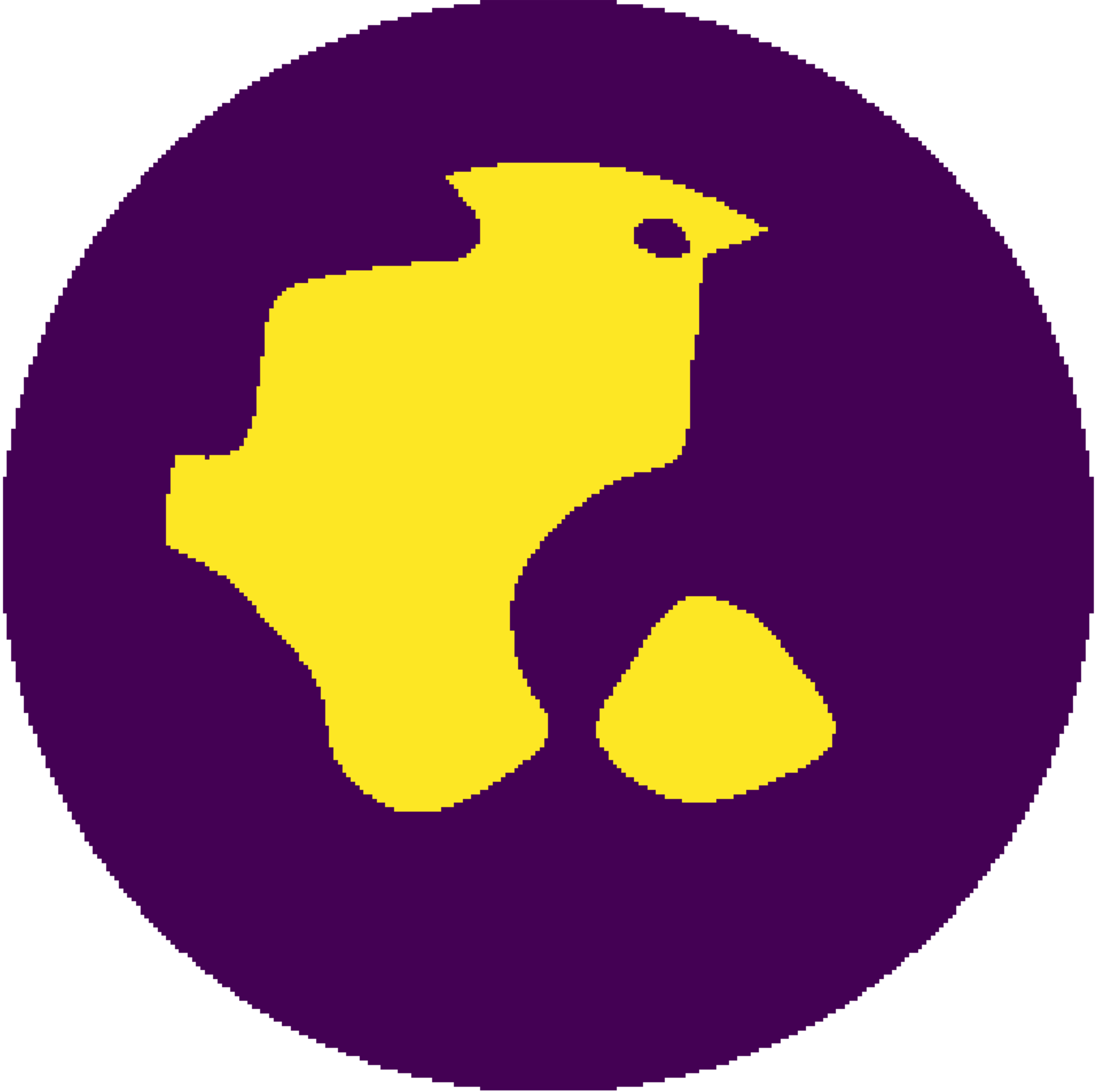}
		\caption{$\gamma$}
		\label{subfig:eit_out}
	\end{subfigure}
	\caption{Electrical impedance tomography from \cref{ex:eit}. The domain $\Omega$ is the unit disk in $\R^2$. The boundary manifold $\partial\Omega$ is identified as the one-dimensional unit torus $\T$ by transforming into polar coordinates. \Cref{subfig:eit_finite} displays two pairs of randomly sampled Neumann data ($g_m$, black) and corresponding Dirichlet data ($\Lambda_\gamma g_m$, red) on the torus. In \cref{subfig:eit_in}, the full NtD map $\Lambda_\gamma$ is represented numerically using its $L^2(\partial\Omega\times\partial\Omega)$ integral kernel function; this kernel function is further shifted, re-scaled, and transformed to $\T^2$ to aid visualization. The true realization of the conductivity $\gamma$ generating these measurements under the forward map $\gamma\mapsto\Lambda_\gamma$ is shown in \cref{subfig:eit_out}. Solving the inverse conductivity problem amounts to evaluating the inverse map $\Lambda_\gamma\mapsto\gamma$.
	}
	\label{fig:eit_viz}
\end{figure}

\begin{example}[electrical impedance tomography]\label{ex:eit}
	The mathematical model of EIT begins with the linear elliptic PDE
	\begin{align}\label{eqn:elliptic_eit}
		\begin{cases}
			\begin{alignedat}{2}
				\nabla\cdot(\gamma \nabla \phi)&=0  && \qin \Omega\,,\\
				\gamma\dfrac{\partial \phi}{\partial \mathsf{n}}&=g && \qon \partial \Omega\,.
			\end{alignedat}
		\end{cases}
	\end{align}
	In \eqref{eqn:elliptic_eit}, $\Omega\subset\R^d$ is the medium of interest, $\gamma\in L^\infty(\Omega)$ is the strictly positive almost everywhere electrical conductivity field, $\phi$ is the electric potential, $g$ is the applied boundary current, and $\partial/\partial\mathsf{n}$ is the outward normal derivative to the boundary $\pOmega$. We assume that $\pOmega$ is smooth. The PDE defines a solution operator that maps a given pair $(\gamma,g)$ to $\phi$ solving \eqref{eqn:elliptic_eit}. However, in EIT, interior quantities such as $\phi$ are not accessible. Instead, for a fixed $\gamma$, one can vary the boundary current density $g$ and measure the corresponding boundary voltages $\phi|_{\pOmega}$. This measurement process leads to a linear operator $\Lambda_\gamma$ called the \emph{Neumann-to-Dirichlet (NtD) map}, which is defined by
	\begin{align}\label{eqn:ntd_eit_defn}
		g\mapsto \Lambda_\gamma g\defeq \phi|_{\partial\Omega}
	\end{align}
	subject to \eqref{eqn:elliptic_eit}. It is known that the NtD map $\Lambda_\gamma\in\cL(H^{-1/2}(\partial\Omega);H^{1/2}(\partial\Omega))$ is a continuous linear operator from the dual Sobolev space $H^{-1/2}(\partial\Omega)=(H^{1/2}(\partial\Omega))^*$ into the trace Sobolev space $H^{1/2}(\partial\Omega)$ of $H^1(\Omega)$ boundary values \cite{mueller2012linear}. In the notation of \eqref{eqn:ip_main}, we take $\cU\defeq L^\infty(\Omega)$, $\cY\defeq\cL(H^{-1/2}(\partial\Omega);H^{1/2}(\partial\Omega))$, $u\defeq\gamma$, and $\cG(\gamma)\defeq \Lambda_\gamma$. The Calder\'on inverse problem is to reconstruct the unknown $\gamma$ from $\Lambda_\gamma$. \Cref{fig:eit_viz} illustrates the main components of EIT.
	
	The inverse problem of finding $\gamma$ from $\Lambda_\gamma$ is hard. Indeed, we wish to recover a function, often in an infinite-dimensional space, on the interior of $\Omega$ only by probing the boundary $\partial\Omega$. More quantitatively, EIT is severely ill-posed from a stability perspective. That is, for sufficiently smooth conductivities $\gamma_0$ and $\gamma_1$ bounded away from zero and infinity, it holds for some $\beta>0$ that
	\begin{align}\label{eqn:stability_eit_inside_range}
		\norm{\gamma_0-\gamma_1}_{L^\infty(\Omega)}\lesssim \log^{-\beta}\Bigl(\norm{\Lambda_{\gamma_0} - \Lambda_{\gamma_1}}_{\cL(H^{-1/2}(\partial\Omega);H^{1/2}(\partial\Omega))}^{-1}\Bigr)
	\end{align}
	whenever $\delta=\norm{\Lambda_{\gamma_0} - \Lambda_{\gamma_1}}_{\cL(H^{-1/2}(\partial\Omega);H^{1/2}(\partial\Omega))}$ is small enough~\cite{alessandrini1988stable,barcelo2007stability}. This means that for the parameter error $\norm{\gamma_0-\gamma_1}_{L^\infty(\Omega)}$ to be less than $\ep$, it suffices for the measurement error $\delta$ to be exponentially small: $\delta\lesssim \exp(-\ep^{-1/\beta})$. This estimate only holds in the range of the forward map $\cG\colon \gamma\mapsto \Lambda_{\gamma}$. See \cref{fig:ip_in} for a visualization. Results similar to \eqref{eqn:stability_eit_inside_range} remain valid outside of the range of $\cG$, as \cref{sec:ee_math_noisy,sec:ee_theory_eit} will discuss. These extensions are useful to handle realistic noisy measurements~\cite{de2025extension}.
\end{example}

\begin{figure}[tb]
	\centering
 	\captionsetup{skip=10pt}
	\begin{tikzpicture}[scale=1.0]
		
		% Boxes
		\draw[thick] (-4.5,-3) rectangle (-0.5,3);
		\draw[thick] (1.5,-3) rectangle (5.5,3);
		
		% Labels
		\node at (-2.5,3.3) {Parameter Space};
		\node at (3.5,3.3) {Data Space};
		\node at (-1,2.8) {$\cU$};
		\node at (5.05,2.8) {$\cY$};
		
		% Curved sets
		\draw[thick, domain=-2:1, smooth, variable=\y] plot ({-4+0.3*(\y+2)^2}, \y);
		\node at (-3.7,-2.4) {$\domain(\cG)$};
		
		\draw[thick, domain=-3:4]
		plot[smooth] coordinates {(1.8,2.1) (2,2) (3.2,1.4) (5,1) (2,-2)};
		\node at (2.6,-2.4) {$\cG\bigl(\domain(\mathcal{G})\bigr)$};
		
		% Points
		\filldraw (-2.1,0.5) circle (1.2pt) node[below] {$\gamma_1$};
		\filldraw (3.2,1.4) circle (1.2pt) node[above=1.5pt] {$\Lambda_{\gamma_1}$};
		\filldraw (-3.3,-0.5) circle (1.2pt) node[below] {$\gamma_0$};
		\filldraw (4.5,1.3) circle (1.2pt) node[above] {$\Lambda_{\gamma_0}$};

		% Arrows
		\draw[arrows = {-Latex[width=5pt, length=10pt]}] (-2.1,0.5) -- (3.2,1.4) node[midway, above, sloped] {$\mathcal{G}$};
		\draw[arrows = {-Latex[width=5pt, length=10pt]}] (-3.3,-0.5) -- (4.5,1.3) node[midway, below, sloped] {$\mathcal{G}$};
		
		% Spaces
		\node at (-4.6,2.6) {};
		\node at (5.4,-2.6) {};
		
	\end{tikzpicture}
	\caption{Stability within the range of $\cG\colon\gamma\mapsto\Lambda_\gamma$ from \cref{ex:eit}, \cref{eqn:stability_eit_inside_range}.}
	\label{fig:ip_in}
\end{figure}

The remainder of this section provides an overview of inverse problem solution methods. \Cref{sec:back_ip_prob} considers techniques to solve inverse problems that heavily rely on ideas from probability and statistics. In contrast, \cref{sec:back_ip_reg} discusses a complementary viewpoint based on the functional-analytic regularization of inverse problems. We refer the reader to \cite{calvetti2018inverse} for a more targeted review of the probabilistic and regularization approaches to inverse problems.

\subsection{Probabilistic methods}\label{sec:back_ip_prob}
In inverse problems, the observed data are typically indirect, incomplete, heterogeneous, and corrupted by noise. Such data challenges lead to inherent ambiguity in the recovered parameters. To address such ambiguity, probabilistic approaches formulate inverse problems in terms of probability measures, which promotes a systematic treatment of uncertainty in addition to resolving ill-posedness. The central goal is not only to recover an estimate of the unknown parameter $u$ from the model~\eqref{eqn:ip_main}, but also to quantify the uncertainty associated with that estimate. Such probabilistic approaches offer richer solution concepts than pointwise methods do. Indeed, uncertainty quantification (UQ) offers a principled approach to assessing the reliability and confidence of inferred inverse problem solutions, which is crucial for downstream decision-making, model validation, and risk assessment.

The \emph{Bayesian approach} treats the unknown parameter as a random variable (\Cref{sec:back_ip_prob_bayes}). It computes the posterior distribution given the data, combining a prior distribution with a likelihood informed by the forward model and measurement process. This allows for both point estimates, such as the maximum a posteriori (MAP) point or posterior mean, and principled UQ through posterior variances, credible sets, and predictive distributions.

Sampling methods rooted in the \emph{transportation of measure} offer flexible ways to represent and visualize uncertainty, particularly in high-dimensional, multi-modal, and non-Gaussian settings (\Cref{sec:back_ip_prob_transport}). They can be viewed as a scalable alternative to traditional Bayesian posterior sampling algorithms such as Markov chain Monte Carlo (MCMC).

In contrast, \emph{frequentist approaches} to statistical inverse problems draw on nonparametric statistics, minimax theory, and empirical risk minimization (\Cref{sec:back_ip_prob_stat}). These methods yield confidence intervals, model selection procedures, and adaptive estimators with frequentist guarantees.

The \emph{measure-centric perspective} formulates new inverse problems, or reinterprets existing ones, by modeling parameters or data as elements of the space of probability measures (\Cref{sec:back_ip_prob_measure}).
As a special case, the so-called stochastic inverse problem framework focuses on pushing forward or pulling back probability measures through the forward map $\cG$. It emphasizes consistency between distributions over the parameter and data spaces and enables uncertainty propagation even in poorly identified regimes.

Each of these methodologies contributes a complementary lens on the incorporation of uncertainty into the solution procedure for inverse problems.

\subsubsection{Bayesian inverse problems}\label{sec:back_ip_prob_bayes}
The Bayesian approach to inverse problems is perhaps the most well-developed of all the probabilistic methods mentioned earlier. In recent decades, there has been significant progress and success in integrating Bayesian inference with statistical inverse problems formulated in function spaces \cite{stuart2010inverse,nickl2023bayesian}. This has led to comprehensive well-posedness theory, approximation theory, and posterior consistency theory \cite{agapiou2013posterior,dashti2017} on the one hand, and scalable computational algorithms on the other hand \cite{cotter2013mcmc,sanz2023inverse}. By leveraging an inherently infinite-dimensional framing, the nonparametric Bayesian approach is naturally applicable to inverse problems from continuum applied mathematics and scientific computing. 

The Bayesian statistician models all quantities in the inverse problem~\eqref{eqn:ip_main} probabilistically. The noise $\eta$ is assumed to be a random variable with distribution  $\nu_0\defeq \Law(\eta)$. This and the forward map $\cG$ define the \emph{likelihood} of observing $y$ given the parameter $u$; we denote by $\nu_{\cG(u)}\defeq\Law(y\condbar u)$ the pushforward of $\nu_0$ under the translation-by-$\cG(u)$ map. The core of any Bayesian procedure is the \emph{prior distribution} specification. The parameter $u$ is modeled as a random variable with distribution $\mu\in\sP(\cU)$, the prior. The prior probability measure $\mu$ represents the modeler's beliefs and existing knowledge before observing any data. The prior $\mu$ is known to regularize the ill-posedness of the inverse problem. Popular choices in infinite dimensions include Gaussian priors, Besov priors, or data-driven pushforwards thereof \cite{dashti2017}.

The Bayesian solution to the inverse problem is not just a single point estimate of the parameter $u\in\cU$, but instead the entire posterior probability distribution $\mu^y\defeq\Law(u\condbar y)\in\sP(\cU)$ obtained by conditioning the prior on the random variable $y$ according to \eqref{eqn:ip_main}. Assume that $\nu_{\cG(u)}$ has density $\exp(-\Phi(u;\slot))$ with respect to $\nu_0$ for all $u$, $\mu$ almost everywhere (a.e.). Then, under some mild conditions that can be verified on a case-by-case basis, the infinite-dimensional Bayes' rule \cite[Thm.~14, pp.~9--10]{dashti2017} shows that the posterior is given by
\begin{align}\label{eqn:posterior}
	\begin{split}
		\mu^y(du)&=\frac{1}{Z_\mu^y}\exp\bigl(-\Phi(u;y)\bigr)\mu(du)\,,\qw \\
		Z_\mu^y&\defeq \int_\cU \exp\bigl(-\Phi(u';y)\bigr)\mu(du')\,.
	\end{split}
\end{align}
The formula \eqref{eqn:posterior} is valid for every $y$ a.e. under the joint distribution of $(u,y)$. The normalizing constant $Z_\mu^y$, also known as the \emph{evidence}, depends on both the data $y$ and the prior $\mu$. The potential function $\Phi\colon\cU\times\cY\to\R$, also called the negative log-likelihood, depends on the forward map $\cG$. A typical example is the Gaussian likelihood. Taking $\cY=\R^J$ and $\nu_0=\normal(0,\Gamma)$ for some positive definite $\Gamma\in\R^{J\times J}$, we have $\Phi(u;y)=\frac{1}{2}\norm{\Gamma^{-1/2}(y-\cG(u))}_{\R^J}^2$ \cite{stuart2010inverse}.

There is by now a mature understanding of the statistical and computational \cite{nickl2022polynomial} performance of Bayesian posterior solutions to linear \cite{agapiou2013posterior,agapiou2014bayesian,agapiou2023heavy,gugushvili2020bayesian,knapik2011bayesian,knapik2016bayes,knapik2018general,ray2013bayesian,trabs2018bayesian} and nonlinear inverse problems \cite{monard2021statistical,nickl2023bayesian,nickl2025bayesian}, especially those arising from PDEs \cite{abraham2020statistical,giordano2020consistency,nickl2020convergence,nickl2020bernstein}. In PDE problems, the map $\cG$ is typically factorized as $\cG=\cO\circ \cF$, where $\cF\colon\cU\to\cV$ is a parameter-to-solution operator and $\cO\colon\cV\to\cY$ is the solution-to-data operator, i.e., the observation operator. Here, the solution refers to the entire PDE solution, which is usually not fully observable in practice. While the interaction between $\cO$ and $\cF$ can substantially influence statistical accuracy~\cite{huang2025operator}, for the purposes of this chapter we are content with working with $\cG$ alone.

\subsubsection{Transport maps}\label{sec:back_ip_prob_transport}
Computational measure transport is a general paradigm that seeks to represent a target probability measure $\nu=\sfT^\star\push\rho$ as the pushforward of some easy-to-sample reference probability measure $\rho$ under a measurable function $\sfT^\star$ \cite{marzouk2016sampling}. Thus, if $u\sim\rho$ is a sample from $\rho$, then $y=\sfT^\star(u)\sim \nu$ is a sample from $\nu$. Such a task can be approximately accomplished by solving the minimum divergence estimation problem
\begin{align}\label{eqn:min_div}
	\min_{\sfT\in\cT} \sfd(\nu, \sfT\push\rho)
\end{align}
for some statistical distance or divergence $\sfd$ and hypothesis class of transport maps $\cT$. One can also add a regularizer to penalize the transport maps. Problem~\eqref{eqn:min_div} falls under the umbrella of variational inference. By adjusting the choice of $\sfd$, $\cT$, and $\rho$, one derives different measure transport algorithms, including those that enable conditional sampling of joint distributions. While the solution to \eqref{eqn:min_div} is generally not ``optimal'' in the sense of optimal transport~\cite{villani2009optimal}, connections to optimal transport can still be made, especially if $\sfd$ is taken to be a Wasserstein distance~\cite{hosseini2025conditional}.

Specializing to the Bayesian setting discussed in~\cref{sec:back_ip_prob_bayes}, to solve inverse problems with measure transport, one typically seeks a prior-to-posterior map $\sfT=\sfT^{(\mu,y)}$ such that $\mu^y=\sfT\push\mu$ \cite[Chp.~5]{bach2024inverse}, where $\mu^y$ is the posterior \eqref{eqn:posterior} corresponding to prior $\mu$ and observed data $y$. The dependence of $\sfT$ on $\mu$ and $y$ is implicit from the previous equality. While not usually done, in some applications it makes sense to explicitly parametrize the dependence of $\sfT$ on $\mu$ or $y$. This way, problem \eqref{eqn:min_div} does not have to be re-solved every time $\mu$ or $y$ is changed, which can be prohibitively expensive. There is work along these lines that develops \textit{amortized solvers} for Bayesian inverse problems \cite{baptista2024conditional,hosseini2025conditional}. This amounts to methods that, given a new piece of data $y$, can rapidly produce a representation of the posterior distribution corresponding to that $y$ without re-training or additional optimization. Using measure transport, this is equivalent to approximating the map $y\mapsto \sfT^{(y)}$, where $\sfT^{(y)}$ is the divergence minimization solution for a particular $y$. The amortized estimation of the data-to-posterior map is arguably a more challenging task than merely emulating the data-to-parameter inverse map, because now the output space is the metric space $\sP(\cU)$ of probability measures supported on $\cU$, rather than just the linear parameter space $\cU$ itself. We return to data-to-posterior estimation in \cref{sec:ee_posterior}.

\subsubsection{Nonparametric statistical inverse problems}\label{sec:back_ip_prob_stat}
Nonparametric statistical inverse problems depart from the nonparametric Bayesian and transport map settings by only modeling the noise $\eta$ as a random variable \cite{cavalier2008nonparametric,bissantz2007convergence,kaipio2005statistical,kaipio2007statistical}. Adopting a frequentist perspective, the unknown parameter $u$ is thought of as some fixed, true function that generates data $y$ according to \eqref{eqn:ip_main}. Several examples can be formulated this way, such as derivative estimation, denoising, deconvolution, or the inverse heat equation~\cite{cavalier2008nonparametric,agapiou2014bayesian}.

The canonical model for a nonparametric statistical inverse problem is
\begin{align}\label{eqn:stat_ip}
	y=\cG(u)+\frac{\sigma}{\sqrt{N}}\xi\,,
\end{align}
where $\xi$ is a (typically Gaussian) white noise process.
The forward operator $\cG\colon\cU\to\cY$ is usually a linear map between function spaces in theoretical studies, although there are some exceptions \cite{nickl2020convergence,bissantz2004consistency}. The parameter $\sigma$ scales the noise, and the number $N$ represents the ``amount of information'' available for recovery. 

The use of the notation $N$ and the choice of white noise are for the following reason. The idealized problem \eqref{eqn:stat_ip} is known to be equivalent---in a precise statistical sense \cite{cavalier2008nonparametric,abraham2020statistical}---to the nonparametric inverse regression problem
\begin{align}\label{eqn:stat_ip_regress}
	y_n=\bigl(\cG(u)\bigr)(x_n) + \sigma \xi_n\qf n\in\{1,\ldots, N\}\,,
\end{align}
where $\xi_n\diid \normal(0,1)$ are independent and identically distributed (i.i.d.) standard Gaussians and the $\{x_n\}$ are some observed design points in the domain of function space $\cY$. Thus, $N$ represents the sample size of available data. The finite-sample problem~\eqref{eqn:stat_ip_regress} is also known as statistical inverse learning \cite{blanchard2018optimal,helin2024least}.

Theoretical analysis is easier for the continuum model \eqref{eqn:stat_ip}. However, it must be carefully interpreted (i.e., in a weak sense) because $\norm{\xi}_{\cY}=\infty$ almost surely if $\cY$ is infinite-dimensional; see the book \cite{gine2021mathematical} for more details about this technical point. Thus, infinite-dimensional statistical inverse problems are characterized by ``large'' noise realizations. The error in reconstructions of the parameter $u$ is measured on average with respect to the noise. The parameter $u$ is assumed to belong to a bounded set of functions, typically the unit ball in some smoothness class such as a Sobolev or Besov space. To solve the inverse problem, one develops estimators that aim to minimize the worst-case average error over this set of functions. This is the minimax framework, which enables the classification of estimators as minimax optimal or not \cite[Sec.~2]{cavalier2008nonparametric}. Linking back to \cref{sec:back_ip_prob_bayes}, there is also work that performs frequentist analysis of Bayesian estimators within the minimax optimality framework \cite{knapik2011bayesian,knapik2016bayes,knapik2018general,nickl2023bayesian}.

\subsubsection{Measure-centric formulations}\label{sec:back_ip_prob_measure}
Distribution-valued data is becoming increasingly common in modern data science \cite{panaretos2020invitation,khurana2024learning,geshkovski2023mathematical,oprea2025distributional,huang2024unsupervised}. Usually numerically realized as unstructured point clouds of arbitrary size, measure-valued data appear in a wide range of scientific applications including probabilistic weather forecasting~\cite{bach2025learning}, biomedical modeling~\cite{haviv2024wasserstein}, and computational geometry~\cite{qi2017pointnet}, to name a few. It is natural to wonder if inverse problems can also be formulated in the space of probability measures in order to take advantage of these substantial developments. 

There is indeed an emerging \emph{measure-centric} line of research on distributional inversion \cite{bredies2013inverse,bredies2025sparse,botvinick2025invariant,li2024differential,yang2023optimal,li2025inverse,butler2012computational}. For example, measure-valued parameters arise in continuum formulations of optimal experimental design \cite{hellmuth2025data,jin2024optimal,jin2024continuous,huan2024optimal}. The calibration of physical models using diverse data sources leads to inverse problems where both the parameters and the measurements are measure-valued \cite{akyildiz2024efficient,vadeboncoeur2025efficient}. Unlike Bayesian inversion from \cref{sec:back_ip_prob_bayes}, it is worth emphasizing that in these frameworks, the underlying state space for the parameters is a space of distributions. Applying Bayesian ideas in the setting of measure-valued parameters would require a notion of posterior distribution over the probability measure space, which would be an interesting yet challenging future direction. The focus in this subsection---and other parts of the chapter more broadly---is on non-Bayesian inverse problems with observed data belonging to the space of probability measures. We now elaborate on this setting.

\paragraph*{Collections of measurements}
The starting point is the assumption that instead of observing a single measurement $y\in\cY$, we observe an unordered collection of data $Y^M\defeq \{y_m\}_{m=1}^M$, i.e., a point cloud. A key point is that $M=M(u)$ is allowed to depend on the unknown parameter $u$ of the inverse problem. This has implications for inversion routines that are repeatedly called to solve the inverse problem for different data realizations. Different measure-valued inverse problem frameworks arise from the particular data-generating model that is assumed for $Y^M$. This subsection now presents three such models.

Suppose that the forward map $\cG$ and the noise distribution $\pi\defeq\Law(\eta)$ are known quantities. For $m\in\{1,\ldots, M(u)\}$, the first data-generating model is
\begin{align}\label{eqn:dist_model1_data}
	y_m=\cG(u)+\eta_m\,,\qw \eta_m\sim\pi\,. 
\end{align}
The goal is to find $u$ from $Y^{M(u)}$. The variation in each $y_m$ is due to a new realization of noise $\eta_m$; the parameter $u$ is the same across realizations.

The second model now varies the parameter realizations via
\begin{align}\label{eqn:dist_model2_data}
	y_m=\cG(u_m)+\eta_m\,,\qw u_m\sim\mu \qa \eta_m\sim\pi\,. 
\end{align}
This setup corresponds to measurements of the fixed system $\cG$ at varied parameter configurations $u_m\sim\mu$ \cite{akyildiz2024efficient}. The goal here is to perform \emph{distributional inversion} by finding $\mu\in\sP(\cU)$ from $Y^M$, where $M=M(\mu)$. The unknown measure $\mu$ could represent a Bayesian prior or some inherent physical uncertainty, for example. Suppose either $\cG$ or $\pi$ is unknown. In that case, the joint recovery of either of these objects in addition to $\mu$ becomes an instance of blind deconvolution \cite{vadeboncoeur2025efficient}, which is a notoriously difficult nonlinear inverse problem.

The last model, which will be particularly relevant for the EIT problem from \cref{ex:eit}, is the generalized regression model
\begin{align}\label{eqn:dist_model3_data}
	y_m=\cG_m(u)+\eta_m\,,\qw \eta_m\sim\pi
\end{align}
and $\{\cG_m\}_{m=1}^{M(u)}$ is a (potentially random) family of ``forward maps.'' The goal of this inverse problem is to find $u$ from the known pairs $\{(\cG_m,y_m)\}_{m=1}^{M(u)}$.

A traditional inversion approach for \eqref{eqn:dist_model1_data} or \eqref{eqn:dist_model3_data} would involve stacking all $M$ equations on $\cY$ into a single equation on the product space $\cY^M$ \cite{sanz2023inverse}. Then, the enlarged inverse problem can be solved with Bayesian or optimization methods. However, this strategy may impose an ordering on the indices $\{1,\ldots,M\}$ that is not fundamental to the inverse problem. It also raises difficulties when applying operator learning to approximate the inverse map $Y^{M(u)}\mapsto u$ because such a map would have to be re-trained from scratch if $M=M(u)$ changes when $u$ changes. Relatedly, a na\"ive attempt to solve the distributional inversion model~\eqref{eqn:dist_model2_data} by solving $M(\mu)$ individual inverse problems does not take into account correlations between points from the statistical model. One path forward is to work with distribution-valued data instead.

\paragraph*{Lifting to the space of probability measures}
We unify the preceding perspectives by lifting equations on the data space $\cY$ to equations on the infinite-dimensional space $\cP(\cY)$ of probability measures over $\cY$. To this end, let $\delta_z$ denote the Dirac measure, or atom, defined by $\delta_z(A)=1$ if $z\in A$ and $\delta_z(A)=0$ otherwise. We identify the point cloud $\{y_m\}_{m=1}^M$ with the empirical probability measure $\nu^M\in\sP(\cY)$ defined by
\begin{align}\label{eqn:empirical_measure}
	\nu^M\defeq \frac{1}{M}\sum_{m=1}^M\delta_{y_m}\,.
\end{align}
By letting $\nu^M$ be our representation of the point cloud data, we can design inversion procedures that are robust to the number $M$ of measurements and are also invariant---or at least equivariant---to permutations of the indices of the atoms in $\nu^M$.
At the population level for \eqref{eqn:dist_model1_data}, we seek $u\in\cU$ such that
\begin{align}\label{eqn:dist_model1_lift}
	\nu^{M(u)}\approx \pi * \left( \cG\push\delta_u\right)\,,
\end{align}
where $\cG\push\colon \sP(\cU)\to\sP(\cY)$ is the pushforward operator under $\cG$ and $*$ denotes convolution of measures \cite{akyildiz2024efficient}. The right-hand side of \eqref{eqn:dist_model1_lift} is the distribution of $y_m$ from \eqref{eqn:dist_model1_data}. Closeness of this distribution to $\nu^{M(u)}$ should be quantified by a statistical distance that is still meaningful when applied to two empirical measures, such as maximum mean discrepancy (MMD) or the family of Wasserstein metrics. This then suggests an optimization principle for finding $u$.

Similarly, the lifted version of the second data model from \eqref{eqn:dist_model2_data} is
\begin{align}\label{eqn:dist_model2_lift}
	\nu^{M(\mu)}\approx \pi * \left( \cG\push\mu\right)\,,
\end{align}
where now $M=M(\mu)$ is allowed to depend on $\mu$. Transport map generative models have been used to find $\mu$ from $\nu^{M(\mu)}$ according to \eqref{eqn:dist_model2_lift} \cite{akyildiz2024efficient,vadeboncoeur2025efficient}.

For clarity, we defer the discussion of lifting the third model, \eqref{eqn:dist_model3_data}, until we describe it in the concrete setting of EIT. 
Regardless, all three of the preceding multi-measurement inverse problems can be written as finding a distribution $\mu\in\sP(\cU)$ from $\sP(\cY)$-valued data
\begin{align}\label{eqn:ip_distributional}
	\nu^{M(\mu)}\approx \nu\defeq \sfG(\mu)\,.
\end{align}
The lifted forward operator $\sfG\colon \sP(\cU)\to\sP(\cY)$ acts on probability measures, but in the case that $u$ is deterministic, we identify $u$ with the Dirac measure $\mu\equiv \delta_u$ in \eqref{eqn:ip_distributional}. In the case of \eqref{eqn:dist_model1_lift} and \eqref{eqn:dist_model2_lift}, the right-hand side of these displays shows that $\sfG=\pi*\cG\push$ is the composition of a convolution with a pushforward. Both convolution and pushforward, and hence $\sfG$, can be interpreted as ``linear'' operations on measures under certain conditions. However, the measure-centric framework is still amenable to more complicated and highly nonlinear operators $\sfG$ if they arise.

The formal $M=\infty$ problem implied by \eqref{eqn:ip_distributional} is to find $\mu$ from population data $\nu=\mathsf{G}(\mu)$. However, rigorously passing to the limit $M\to\infty$ to relate the sequence of finite atom measurements and solutions to population measurements and solutions could be challenging. Such analysis would depend on assumptions placed on the statistical model (e.g., independence, misspecification) and the type of convergence desired. Taking these limits can be technical, as work in a similar setting has shown~\cite{fiedler2023kernel,fiedler2024statistical,fiedler2025recent}. This is a clear direction for future work.

\paragraph*{Application to EIT}
Calder\'on's problem from \cref{ex:eit} is formulated at the continuum level as finding $\gamma\in\cU$ from noisy data $Y\defeq \Lambda_\gamma + \Xi$ belonging to a space of linear operators. Here, $\Xi$ is an abstract operator-valued noise process. The NtD map $\Lambda_\gamma$ represents all possible current-voltage pairs on the boundary. However, practical EIT hardware only permits a finite number $M=M(\gamma)$ of boundary measurements \cite{dunlop2016bayesian,mueller2012linear}. For $m\in\{1,\ldots, M(\gamma)\}$, we observe functions $g_m$ and
\begin{align}\label{eqn:eit_finite_data_model}
	y_m=\Lambda_\gamma g_m +\eta_m\,,\qw g_m\sim\rho \qa\eta_m\sim\pi\,.
\end{align}
The known distribution $\rho\in\sP(\cX)$ over a space $\cX$ of current patterns represents the experimental design of the EIT imaging system. See \cref{subfig:eit_finite} for a visualization. The goal is to recover the conductivity $\gamma$ from the point cloud $\{(g_m,y_m)\}_{m=1}^{M(\gamma)}$. One approach is to perform a two-stage recovery: first, estimate the nuisance parameter $\Lambda_\gamma$ via regression, then use $\Lambda_\gamma$ to perform inversion using the continuum Calder\'on framework. However, it is desirable to avoid the nuisance parameter $\Lambda_\gamma$ and instead directly reconstruct $\gamma$ from the raw data boundary measurement pairs $\{(g_m,y_m)\}_{m=1}^{M(\gamma)}$.

We do so as follows. First, we relate the EIT data model \eqref{eqn:eit_finite_data_model} to the general regression model \eqref{eqn:dist_model3_data} by defining $\cG_m(u)\defeq \Lambda_u g_m$. Then by lifting the observed data to the joint space $\sP(\cX\times \cY)$, we obtain the measure-valued inverse problem
\begin{align}\label{eqn:dist_model3_lift}
	\frac{1}{M(\gamma)}\sum_{m=1}^{M(\gamma)}\delta_{(g_m, y_m)}\approx (\delta_0\otimes \pi) * (\Id, \Lambda_\gamma)\push\rho \eqdef \sfG(\delta_\gamma)\,. 
\end{align}
In the display, $\delta_0\otimes\pi\in\sP(\cX\times\cY)$ is a product measure. Under the measure-centric perspective, the inverse problem solution operator corresponding to \eqref{eqn:dist_model3_lift} maps $\sP(\cX\times\cY)$ to $\cU$, that is, a joint probability measure over functions on the boundary to functions on the interior.
We will return to this measure-to-function viewpoint on EIT in \cref{sec:ee_archparam_nio}.

\paragraph*{Stochastic inverse problems}
A special case of \eqref{eqn:dist_model2_data} and its lifting \eqref{eqn:dist_model2_lift} is when $\pi=\Law(\eta)=\delta_0$. This noiseless setting leads to the so-called \emph{stochastic inverse problem} (SIP)
\begin{align}\label{eqn:sip}
	\nu=\cG\push\mu\,.
\end{align}
The data is $\nu\in\sP(\cY)$ and the unknown is $\mu\in\sP(\cU)$.
There is a growing body of work on SIPs~\cite{breidt2011measure,butler2012computational,butler2014measure} and related methodologies, such as \textit{data-consistent inversion}~\cite{butler2018combining,butler2020data,butler2025optimal,bergstrom2024distributions}. These approaches extend classical inverse problem frameworks by treating both the data and the unknown parameters as probability measures and modeling their relationships through forward operators defined between metric spaces of probability measures. Recent advances in SIPs have leveraged variational optimization in probability measure spaces, explored the effects of different statistical divergences (e.g., $f$-divergences and Wasserstein distances) used for parameter reconstruction, addressed challenges in conditional and marginal recovery, and introduced dynamical solution techniques. These methods provide robust stability estimates, convergence guarantees, regularization effects, and uncertainty quantification~\cite{li2024stochastic,li2024differential,li2025inverse}. A comparison between SIPs and Bayesian inverse problems can be found in \cite{bingham2024inverse}.

Distributional inversion greatly benefits from modern techniques for analysis in probability measure spaces, including structured transport maps~\cite{marzouk2016sampling,baptista2024representation,el2012bayesian} and Wasserstein gradient flows~\cite{li2024stochastic,li2025inverse,chewi2025wassersteinapp,chewi2025wasserstein}. These techniques are also relevant for learning distribution-informed priors in standard Bayesian inversion workflows, especially in scenarios where data from multiple sources must be incorporated~\cite{akyildiz2024efficient,vadeboncoeur2025efficient,white2024building}. Although the forward operator in most SIP formulations is based on a pushforward map that links the parameter distribution to the observed data distribution as in \eqref{eqn:sip}, alternative operators between probability measure spaces have also gained attention. For example, the cryo-EM imaging problem is a notable instance of an SIP where the forward operator itself involves randomness~\cite{cryo-em-sip}.

\subsection{Regularization methods}\label{sec:back_ip_reg}
In contrast to the chapter so far, this subsection adopts a more deterministic viewpoint on inverse problems. Here, the goal is to solve the equation  
\begin{align}\label{eqn:ip_deterministic}  
	y = \cG(u)  
\end{align}  
for the unknown parameter $u \in \cU$. However, instead of the exact observation $y\in\cY$, we are provided with a perturbed measurement  
\begin{align}\label{eqn:deterministic_perturb}  
	y^{\delta} \in \cY\,, \qw \sfd_{\cY}(y, y^{\delta}) \leq \delta\,.
\end{align}  
In \eqref{eqn:deterministic_perturb}, $\sfd_\cY$ is a distance function on the data space and $\delta$ quantifies the magnitude of the measurement error. In this setting, reconstruction accuracy is studied for the worst-case data perturbation satisfying \eqref{eqn:deterministic_perturb}. This departs from the nonparametric statistical inverse problems setting of \cref{sec:back_ip_prob_stat}, which considers average-case error with respect to random noise in the measurements.

Such inverse problems are ill-posed: the forward operator $\cG$ may be non-invertible or ill-conditioned, making the solution highly sensitive to small perturbations in the data. To mitigate this instability, \emph{regularization} incorporates prior information or constraints on the parameter space $\cU$ to stabilize the inversion process~\cite{engl1996regularization,benning2018modern}. Regularization strategies for inverse problems of the form \eqref{eqn:ip_deterministic}--\eqref{eqn:deterministic_perturb} are commonly grouped into three broad paradigms: direct methods (e.g., spectral filtering via the SVD), implicit techniques (e.g., early stopping and the dynamics of gradient descent), and variational formulations (e.g., Tikhonov and total variation regularization).

Of the three classifications, this chapter centers mostly on variational regularization. Variational regularization is an optimization approach \cite{haber2000optimization}. The solution is obtained by minimizing a composite objective:
\begin{align}\label{eqn:reg_opt}
	\min_{u\in\cU} \Bigl\{\sfd_{\cY} (\cG(u), y^{\delta}) + \lambda R(u)\Big\} \,,
\end{align}
where $R\colon \cU\to\R $ is a chosen regularizer and $\lambda \geq 0$ is a regularization parameter. The first term in~\eqref{eqn:reg_opt} quantifies data misfit, while the second term involving $R$ encodes prior knowledge about the parameter $u$.  
Classical choices for the regularization functional $R$ include Tikhonov's quadratic penalty $R(u) = \|u\|_{\cU}^2$~\cite{groetsch1984theory,golub1999tikhonov,bishop1995training} and total variation (TV) regularization $R(u) = \norm{\nabla u}_{L^1}$~\cite{rudin1992nonlinear,osher2005iterative,chan2005aspects}, which promote smoothness and edge sparsity, respectively. While handcrafted regularizers have proven to be effective in specific applications~\cite{bredies2010total, lou2010image, zhang2010bregmanized, zhang2015total}, they often fall short in capturing the complex, high-dimensional structures present in real-world data. Moreover, the choice and tuning of the hyperparameter $\lambda$ remain significant practical challenges~\cite{chirinos2024learning}.

In recent years, data-driven approaches have significantly improved the quality of reconstructions by learning regularizers or priors directly from data~\cite{burger2024learned, duff2024regularising,arridge2019solving}. Modern machine learning techniques for regularization can be categorized based on their training paradigm and model type~\cite{habring2024neural}. In supervised approaches, networks are trained on paired examples $(y, u)$ to learn either the full inverse map or a regularization functional. In contrast, unsupervised methods aim to model the distribution of $u$ (i.e., a prior $\mu$) using only samples of $u$ from historical datasets. There are also untrained approaches, such as Deep Image Prior (DIP)~\cite{ulyanov2018deep}, which impose priors implicitly during instance-specific inversion without requiring any offline training. Most of these methods require or assume some knowledge about the forward operator $\cG$.

Among these advances, score-based diffusion models have emerged as powerful priors for use in solving inverse problems. These models leverage learned gradients of log densities to guide stochastic sampling conditioned on observed data. This enables unsupervised, high-fidelity reconstructions that replace traditional handcrafted regularizers with learned generative dynamics. Current research efforts are focused on improving the efficiency of these models while incorporating physics-based constraints~\cite{feng2023score, chung2023diffusion, daras2024survey, tolooshams2025equireg}.

More recently, operator learning, which focuses on estimating mappings between function spaces, is gaining significant attention in the context of inverse problems. Beyond the more common tasks of learning forward or inverse operators from paired data, the operator learning framework offers new opportunities to derive or incorporate regularization techniques to tackle challenging inverse problems. \Cref{sec:reg} will explore two key categories of methods for learning priors and regularization mechanisms: empirical Bayes and denoising-driven regularization, notably the plug-and-play technique.

Empirical Bayes-type methods aim to approximate prior distributions by learning them from within a parametric family of probability measures. This is achieved by optimizing statistical divergence-based objective functions, which is also known as variational inference. In contrast, plug-and-play methods leverage the insight that denoisers inherently encode prior information. By embedding a learned denoiser into iterative optimization algorithms, typically replacing a proximal operator, plug-and-play methods regularize inverse problems without requiring an explicit prior. As we will examine, these prior learning and regularization approaches can be interpreted within the broader framework of operator learning, where the mathematical properties of the learned operators (e.g., Lipschitz continuity, monotonicity) play a crucial role in ensuring algorithmic stability and convergence when solving inverse problems.

%%%%%%%%%%%%%%%%%%%%%%%%%%%%%%%%%%%%%%%%%%%%%%%%%%%%%%%%%%%%%%%%%%%%%%%%%%%%%%%%%%%%%%%%%%%%%%%%%%%%%%%%%%%%%%%%%%%%%%%%%%%%%%%%%%
\section{Background on operator learning}\label{sec:back_ol}
Supervised operator learning lifts finite-dimensional vector-to-vector regression to infinite-dimensional function-to-function regression. The framework designs approximation architectures at the continuum level. Consequently, these architectures exhibit robustness to the level of discretization used during numerical implementation \cite[Sec.~2.3, pp.~8--9]{kovachki2021neural}. This is a core distinction between operator learning and traditional machine learning. Once trained, operator learning models can be deployed at different discretizations of the input and output spaces without additional effort, while at the same time produce consistent predictions. The primary motivators for operator learning are the acceleration of scientific computation and the discovery of physical laws from real data. Operator learning has been particularly successful in surrogate modeling of parametrized PDEs and their corresponding downstream tasks, such as uncertainty propagation or engineering design optimization, based on these PDE parameter-to-solution map surrogates.

Operator learning developed independently from the field of inverse problems and does not explicitly rely on ideas from that field. However, it is known that forward operator learning \cite{de2023convergence,mollenhauer2022learning,lu2022data,nelsen2024statistical}---and nonparametric regression more generally \cite{de2005learning}---can be interpreted as a linear inverse problem with a non-compact forward map. The non-compact linear forward map arises from point evaluations of the unknown and possibly nonlinear target operator at the randomly sampled input data points. This perspective has proven to be highly effective in the study of linear regression and kernel methods, where sharp error estimates are developed by combining concentration of measure techniques in Hilbert spaces with concrete manipulations of integral operators \cite{caponnetto2007optimal,rosasco2010learning,smale2007learning}. 

Another related line of work exploits the Green's function structure of differential operators to derive data-efficient operator learning algorithms that are inspired by numerical linear algebra \cite{boulle2022data,boulle2022learning,boulle2023learning,boulle2024mathematical}. This framework is generalized to a statistical setting that allows for recovery of kernels in possibly nonlinear operators that do not necessarily admit Green's functions \cite{chada2024data,lu2022data,lu2024nonparametric,zhang2025minimax}.
In contrast, the present chapter focuses on solving a variety of inverse problems that arise from the physical, imaging, and data sciences, rather than the specific non-physical inverse problem associated with the act of supervised learning.

For a more comprehensive overview of operator learning from a mathematical and theoretical perspective, we refer the reader to recent survey articles \cite{boulle2024mathematical,kovachki2024operator,subedi2025operator}. The present section only covers material that is necessary for inverse problem purposes. \Cref{sec:back_ol_learn} presents the basic components of an operator learning workflow, while \cref{sec:back_ol_arch} defines common neural operator architectures that appear frequently in the rest of this chapter.

\subsection{Supervised learning of operators}\label{sec:back_ol_learn}
Suppose that $\cX_1=\cX_1(\Omega_1;\R^{p_1})$ and $\cX_2=\cX_2(\Omega_2;\R^{p_2})$ are input and output Banach spaces of vector-valued functions defined on bounded Euclidean subsets $\Omega_1$ and $\Omega_2$, respectively. Let $\sfF\colon \cX_1\to\cX_2$ denote the target operator that we wish to learn. As a concrete example, consider $\sfF\colon L^\infty(\Omega;\R_{>0})\to H^1(\Omega;\R)$ defined by $\gamma\mapsto \phi$ in \eqref{eqn:elliptic_eit}. This is a nonlinear operator defining the coefficient-to-solution map of the elliptic boundary value problem \eqref{eqn:elliptic_eit}. More generally, we could also consider learning the well-posed forward operator $\cG$ appearing in inverse problems such as \eqref{eqn:ip_main}. These examples are representative of much of the field of operator learning: approximation of well-defined, potentially very smooth, forward processes that map physical parameters to PDE solutions or quantities of interest.

In supervised operator learning, we first assume that a dataset of labeled pairs $\{(u_n, y_n)\}_{n=1}^N\subset\cX_1\times\cX_2$ is available for training. The $u_n\colon x\mapsto u_n(x)$ are the input functions and the $y_n\colon x'\mapsto y_n(x')$ are the possibly noisy output functions
\begin{align}\label{eqn:ol_noise_model}
	y_n=\sfF(u_n)+\xi_n\,.
\end{align}
In statistical learning, the data are assumed to be i.i.d.~samples from some joint distribution over $\cX_1\times\cX_2$. Such an assumption may not hold for highly correlated processes, such as time series or optimal experimental designs. Nonetheless, in this section we further assume that $u_n\diid \mu$ for some $\mu\in\sP(\cX_1)$ and that $\xi_n\equiv 0$ for all $n$. Often in operator learning, the labeled data are not passively given, but instead are actively queried using the map $\sfF$ as a black box. Thus, the training data distribution $\mu$ is frequently specified by the user. While it is of considerable interest to determine optimal choices for $\mu$~\cite{boulle2023elliptic,gao2024adaptive,guerra2025learning,subedi2024benefits,satheesh2025picore}, in practice $\mu$ is usually taken to be a Gaussian measure over $\cX_1$ with a specified covariance structure such as Mat\'ern or squared exponential \cite[Sec.~4.1]{boulle2024mathematical}. Then $\mu$ is sampled i.i.d.~to generate the random functions $\{u_n\}$, and consequently $\{y_n\}$ after querying $u_n\mapsto \sfF(u_n)$. The noise-free condition $\xi_n\equiv 0$ is justified when the data arises from high-fidelity numerical solvers with controlled discretization errors; this is a typical use case in operator learning.

The overall goal in operator learning is to estimate the operator $\sfF$ from the given data. The most popular way to do so is by performing the \emph{empirical risk minimization} (ERM)
\begin{align}\label{eqn:erm_ol}
	\min_{\Psi\in \sH }\frac{1}{N}\sum_{n=1}^N\ell\bigl(y_n, \Psi(u_n)\bigr)\,.
\end{align}
Solving this optimization problem is referred to as \emph{training}. In \eqref{eqn:erm_ol}, $\ell\colon \cX_2\times\cX_2\to\R$ is a loss functional and $\sH$ is the hypothesis space of possible operator approximations. The loss is often taken to be the squared loss $\ell(y,y')=\norm{y-y'}_X^2$ or non-squared norm loss $\ell(y,y')=\norm{y-y'}_X$ in any Banach space $X\supseteq \cX_2$. Usually $X$ equals the Lebesgue space $L^q(\Omega_2;\R^{p_2})$ for $q=1,2$ or the Sobolev space $H^k(\Omega_2;\R^{p_2})$ for $k=1,2$. The relative loss 
\begin{align}\label{eqn:rel_loss_output}
	\ell(y,y')=\frac{\norm{y-y'}_X}{\norm{y}_X+\epsilon}
\end{align}
is also popular; here $\epsilon\geq 0$ is a small numerical stability parameter.
Practically, the function space norms appearing in the loss $\ell$ are discretized with quadrature rules and finite difference schemes.

For parametric hypothesis classes based on deep learning, \eqref{eqn:erm_ol} is typically solved with variants of mini-batch stochastic gradient descent (SGD) and accelerated with graphics processing unit (GPU) computing. However, there are more structured operator learning models, such as those based on kernels or random features \cite{nelsen2024operator,schafer2021sparse,batlle2024kernel}, that do not require the explicit minimization \eqref{eqn:erm_ol} with SGD, e.g., because of closed-form expressions for minimizers. These methods still benefit from GPU acceleration, however.

Once a minimizer or approximate minimizer $\widehat{\Psi}$ of \eqref{eqn:erm_ol} is obtained, it remains to assess its accuracy as an approximation to $\sfF$.
It is theoretically convenient, especially in universal approximation results, and also practically desirable in safety-critical applications, to consider the uniform error over a fixed compact subset $\sfK$ of $\cX_1$:
\begin{align}\label{eqn:loss_unif}
	\sup_{u\in\sfK}\,\ell\bigl(\sfF(u), \widehat{\Psi}(u)\bigr) \qor \sup_{u\in\sfK}\,\norm{\sfF(u)-\widehat{\Psi}(u)}_{X}\,.
\end{align}
However, we cannot expect the worst-case error \eqref{eqn:loss_unif} to be small by only training in the empirical norm \eqref{eqn:erm_ol}, which is an averaged quantity. The most used accuracy metric in practice is the expected risk
\begin{align}\label{eqn:loss_expected_risk}
	\E_{u\sim\mu}\bigl[\ell\bigl(\sfF(u), \widehat{\Psi}(u)\bigr)\bigr]\,.
\end{align}
When $\cX_2\subseteq L^2$ and $\ell$ is chosen as \eqref{eqn:rel_loss_output}, we then obtain the popular \emph{expected relative $L^2$ error}
\begin{align}\label{eqn:loss_L2_rel}
	\E_{u\sim\mu}\left[\frac{\norm{\sfF(u)-\widehat{\Psi}(u)}_{L^2(\Omega_2;\R^{p_2})}}{\norm{\sfF(u)}_{L^2(\Omega_2;\R^{p_2})} + \epsilon}\right]
\end{align}
with respect to the training distribution $\mu$. Other test accuracies are defined by choosing different function spaces $X$ and loss functions $\ell$. The out-of-distribution error with respect to some $\mu'\neq\mu$ is also relevant in applications~\cite{guerra2025learning,de2023convergence}; in this case, one would replace every instance of $\mu$ by $\mu'$ in the preceding displays. See \cite[Sec.~4.2.1]{boulle2024mathematical} for other choices of loss functions or loss augmentations. The expectations in expressions such as \eqref{eqn:loss_expected_risk} and \eqref{eqn:loss_L2_rel} are evaluated numerically by averaging over a finite set of held-out test samples.

More broadly, operator learning as a field extends beyond the supervised setting considered here. The present chapter will touch on some unsupervised operator learning problems in \cref{sec:reg}. We also refer the reader to \cite[Chp.~6.2]{nelsen2024statistical} for other examples.

\subsection{Common architectures}\label{sec:back_ol_arch}
By now, there are dozens of variants of operator learning architectures. For an overview of some of these and their shared structure, see \cite[Sec.~3]{boulle2024mathematical,kovachki2024operator,subedi2025operator} and \cite{lanthaler2023nonlocality,nelsen2024operator}. What these models share is a conceptually continuous design that decouples model parameters from any fixed discretization. Some operator learning methods of note that we do not study in detail in this chapter include continuum transformers for functions and probability measures~\cite{guotransformer,cao2021choose,ovadia2024vito,geshkovski2024measure,bach2025learning,calvello2024continuum}, function space random features \cite{nelsen2021random,nelsen2024operator}, kernel methods and Gaussian processes~\cite{batlle2024kernel}, and nonparametric linear operator parametrizations \cite{boulle2022data,boulle2022learning,boulle2023elliptic,boulle2023learning,boulle2024operator,de2023convergence,mollenhauer2022learning,schafer2021sparse}. The remainder of this subsection describes two important parametric classes: encoder-decoder networks in \cref{sec:back_ol_arch_onet} and neural operators in \cref{sec:back_ol_arch_no}. They correspond to a hypothesis space $\sH$ in \eqref{eqn:erm_ol} given by
\begin{align}\label{eqn:hypothesis_parametric}
	\sH\defeq\set{\Psi(\slot; \theta) \colon \cX_1\to\cX_2}{\theta\in\Theta}
\end{align}
for some parameter set $\Theta$ of high but finite dimension.
Within the two classes are the well-known DeepONet, PCA-Net, and Fourier Neural Operator architectures.

\subsubsection{Encoder-decoder operator networks}\label{sec:back_ol_arch_onet}
A mapping $\Psi\colon \cX_1\to\cX_2$ that takes the form
\begin{align}\label{eqn:encode_decode}
	\Psi\defeq \mathscr{D}\circ \mathscr{A}\circ \mathscr{E}
\end{align}
is said to have \emph{encoder-decoder structure} if the encoder $\sE\colon\cX_1\to \R^{d_1}$, approximator $\sA\colon\R^{d_1}\to\R^{d_2}$, and decoder $\sD\colon \R^{d_2}\to\cX_2$ are all continuous maps \cite[Sec.~2.4 and 4.1]{kovachki2024operator}.
Encoder-decoder structure is a fundamental component of most universal approximation proofs for various operator learning methods \cite{kovachki2021neural,lanthaler2022error,lanthaler2023nonlocality,lanthaler2023operator}. It can also be related to autoencoders on function spaces; see \cite[Fig.~3]{kovachki2024operator} and \cite{bunker2024autoencoders,bhattacharya2021model,huang2025operator}. Encoder-decoder networks are also the nonlinear architectures that enjoy the most comprehensive quantitative theory, including sharp approximate rates \cite{lanthaler2023operator} and sample complexity bounds \cite{liu2024deep,kovachki2024data}.

The \emph{DeepONet} architecture is a canonical encoder-decoder operator network~\cite{lu2021learning}. Building off of earlier work in the shallow case~\cite{chen1995universal}, DeepONet constructs $\Psi$ in \eqref{eqn:encode_decode} from traditional, finite-dimensional, deep feedforward neural networks \cite[Chp.~2]{petersen2024mathematical} as follows. Let $\sE=L\in\cL(\cX_1;\R^{d_1})$ be a continuous linear map. For example, $L$ could be represented by projections onto a finite number of basis functions if $\cX_1$ is a Hilbert space, general linear functionals in $\cX_1^*$ if $\cX_1$ is a Banach space, or point evaluations at $d_1$ fixed (or varying, see~\cite{tretiakov2025setonet,prasthofer2022variable}) sensor points in $\Omega_1$ if $\cX_1$ embeds into the space of continuous functions. The approximator $\sA\colon\R^{d_1}\to\R^{d_2}$ is a deep neural network with $d_2$ output components $a_j\colon\R^{d_1}\to\R$ called the ``branch net.'' For another deep neural network with $d_2$ output components $\varphi_j\colon \Omega_2\to\R^{p_2}$ called the ``trunk net,'' the decoder is the linear map $\sD\colon z\mapsto \sum_{j=1}^{d_2}z_j\varphi_j$. The $\{\varphi_j\}$ represent basis functions for the output space $\cX_2$. Putting together the pieces, we can write the DeepONet as
\begin{align}\label{eqn:deeponet}
	\Psi_{\mathrm{DON}}\bigl(u;(\theta_\sA,\theta_\sD)\bigr)(x)=\sum_{j=1}^{d_2}a_j(Lu;\theta_\sA)\varphi_j(x;\theta_\sD)
\end{align}
for all $u\in\cX_1$ and $x\in \Omega_2$. Since both $\{a_j\}$ and $\{\varphi_j\}$ are parametrized as neural networks, we denote their parameter dependence explicitly in \eqref{eqn:deeponet} by $\theta_\sA$ and $\theta_\sD$, respectively. Note that the encoder map $L$ is fixed \emph{a priori}, not learned.

A related architecture, called \emph{PCA-Net} \cite{bhattacharya2021model}, is based on principal component analysis (PCA)~\cite{reiss2020nonasymptotic}. This requires $\cX_1$ and $\cX_2$ to be separable Hilbert spaces, together with the conditions $\E_{u\sim\mu}\norm{u}^2_{\cX_1}<\infty$ and $\E_{u\sim\mu}\norm{\sfF(u)}^2_{\cX_2}<\infty$. Then $\mu$ and $\sfF\push\mu$ (the pushforward of $\mu$ under $\sfF$) have well-defined trace-class covariance operators $\Cov(\mu)\in\cL(\cX_1;\cX_1)$ and $\Cov(\sfF\push\mu)\in\cL(\cX_2;\cX_2)$ whose eigenfunctions (i.e., the PCA bases) determine the encoder $\sE$ and decoder $\sD$. The approximator $\sA$ is still a neural network. Letting $\{e_j\}_{j\in\N}$ be the PCA eigenbasis of $\Cov(\mu)$ ordered in correspondence with its eigenvalues which we sort to be nonincreasing, PCA-Net takes $\sE=L\colon u\mapsto \{\ip{e_j}{u}_{\cX_1}\}_{j=1}^{d_1}$. Letting $\{\varphi_j\}_{j\in\N}$ be the ordered PCA eigenbasis of $\Cov(\sfF\push\mu)$, PCA-Net is defined by
\begin{align}\label{eqn:pcanet}
	\Psi_{\mathrm{PCA}}(u;\theta)(x)=\sum_{j=1}^{d_2}a_j(Lu;\theta)\varphi_j(x)
\end{align}
for all $u\in\cX_1$ and $x\in \Omega_2$. Unlike DeepONet \eqref{eqn:deeponet}, which uses neural networks, the basis functions $\{\varphi_j\}$ in \eqref{eqn:pcanet} are fixed by PCA and are parameter-free. Since $\mu$ and $\sfF\push\mu$ are only accessible through the finite amount of training data $\{(u_n,\sfF(u_n) )\}_{n=1}^N$, the inaccessible true PCA eigenfunctions $\{e_j\}$ and $\{\varphi_j\}$ are replaced by their empirical PCA counterparts in practice.

DeepONet and PCA-Net have linear decoders defined through basis expansions. This severely restricts their approximation power~\cite{lanthaler2022error,lanthaler2023operator}. Encoder-decoder operator networks with nonlinear decoders include NOMAD~\cite{seidman2022nomad} and Shift-DeepONet~\cite{lanthalernonlinear}, which generalize DeepONet, and PARA-Net~\cite{de2022cost}, which generalizes PCA-Net. NOMAD and PARA-Net both take the form
\begin{align}\label{eqn:deeponet_nonlinear}
	\Psi_{\mathrm{NONLINEAR}}\bigl(u;(\theta_\sA,\theta_\sD)\bigr)(x)=Q\Bigl(\bigl\{a_j(Lu;\theta_\sA)\bigr\}_{j=1}^{d_2}, x;\theta_\sD\Bigr)\,,
\end{align}
where the neural network $Q(\slot;\theta_{\sD})\colon \R^{d_2}\times\Omega_2\to\R^{p_2}$ with parameters $\theta_{\sD}$ represents the nonlinear decoder. Shift-DeepONet adopts the more structured parametrization
\begin{align}\label{eqn:deeponet_shift}
	\Psi_{\mathrm{SHIFT-DON}}(u)(x)=\sum_{j=1}^{d_2}a_j(Lu)\varphi_j\bigl(A_j(Lu)x+b_j(Lu)\bigr)\,.
\end{align}
In \eqref{eqn:deeponet_shift}, we suppress the dependence of the neural networks $\{a_j,A_j,b_j,\varphi_j\}_{j=1}^{d_2}$ on their trainable parameters for brevity.
It should be noted that the potential approximation benefits of nonlinear decoding must be balanced with the increase in computational cost~\cite[Sec.~3.4]{de2022cost}.

\subsubsection{Neural operators}\label{sec:back_ol_arch_no}
The encoder-decoder structure from \eqref{eqn:encode_decode} can be viewed as a form of dimension reduction. In contrast, \emph{neural operator} architectures do not rely on dimension reduction into finite-dimensional latent spaces. They instead lift information into higher-dimensional latent function spaces. At their core, neural operators $\Psi_{\mathrm{NO}}\colon\cX_1\to\cX_2$ generalize the compositional feedforward structure of standard deep neural networks to function spaces \cite{kovachki2021neural}. 

To present the main ideas, we specialize to the setting that $\Omega\defeq \Omega_1=\Omega_2\subset\R^d$. Let $d_\mathrm{c}$ be the so-called \emph{channel dimension}, where we typically take $d_\mathrm{c}\geq \max(p_1,p_2)$. Let $\cH=\cH(\Omega;\R^{d_\mathrm{c}})$ be a latent function space. We define
\begin{equation}\label{eqn:no}
	\Psi_{\mathrm{NO}}(u) \defeq \bigl(\mathcal{Q}\circ \mathscr{L}_{T}\circ \mathscr{L}_{T-1}\circ\cdots \circ \mathscr{L}_{2}\circ\mathscr{L}_{1}\circ \mathcal{S}\bigr)(u)
\end{equation}
for all $u\in\cX_1$. The number of hidden layers is denoted by $T$. The map $\cS\colon \cX_1\to\cH$ is a pointwise-defined local lifting operator given by
\begin{align}\label{eqn:pw_lift}
	\bigl(\cS (u)\bigr)(x)\defeq S\bigl(x,u(x)\bigr)
\end{align}
for some map $S\colon \R^d\times\R^{p_1}\to\R^{d_\mathrm{c}}$ that is typically affine, but could also be a neural network.

Similarly, $\cQ\colon\cH\to\cX_2$ is the pointwise-defined local projection operator
\begin{align}\label{eqn:pw_proj}
	\bigl(\cQ (h)\bigr)(x)\defeq Q\bigl(h(x)\bigr)
\end{align}
for a (typically shallow) neural network $Q\colon \R^{d_\mathrm{c}}\to\R^{p_2}$. The nonlinearity of $Q$ is crucial to ensure that $\Psi_{\mathrm{NO}}$ is a nonlinear method of approximation, which means that its range $\Psi_{\mathrm{NO}}(\cX_1)$ is not a vector space \cite{kovachki2021neural}. In contrast, DeepONet and PCA-Net from \cref{sec:back_ol_arch_onet} are linear approximators because their range is the linear span of the basis functions $\{\varphi_j\}_{j=1}^{d_2}$.

For each $t$, the hidden layers $\sL_t \colon \cH \to \cH$ are the nonlinear operators
\begin{equation}\label{eqn:no_layer_nonlinear}
	\bigl(\mathscr{L}_{t}(h)\bigr)(x) = \sigma_t\bigl(W_t h(x) + (\cK_t h)(x) + b_t(x)\bigr)
\end{equation}
for every $h\in\cH$ and $x\in\Omega$. In \eqref{eqn:no_layer_nonlinear}, $\sigma_t\colon\R\to\R$ is a nonlinear activation that acts pointwise on functions; usually $\sigma_T$ is the identity map while the other activations $\sigma_t$ for $t< T$ are taken to be the same non-polynomial nonlinearity (e.g., ReLU, GELU, etc). Moreover, $W_t \in \R^{d_\mathrm{c}\times d_\mathrm{c}}$ is a weight matrix, $b_t\colon \Omega \to\R^{d_{c}}$ is a bias function (usually parametrized as a constant vector), and $\cK_t\colon\cH\to\cH$ is a linear kernel integral operator given by
\begin{align}\label{eqn:no_kernel}
	(\cK_t h)(x)\defeq\int_\Omega \kappa_t(x,y)h(y)\dd{y}
\end{align}
for $h\in\cH$ and $x\in\Omega$. The particular choice of parametrization of the matrix-valued kernel $\kappa_t\colon\Omega\times\Omega\to\R^{d_\mathrm{c}\times d_\mathrm{c}}$ determines the classification of the neural operator architecture. For example, choosing $\kappa_t$ to be a neural network localized to patches defines the Graph Kernel Neural Operator, while choosing a low-rank kernel $\kappa_t$ with learnable factors leads to the Low-Rank Neural Operator \cite[Sec.~4]{kovachki2021neural}. 
Alternatively, expanding $\kappa_t$ in a wavelet basis or Laplace--Beltrami operator eigenfunctions on a manifold leads to the Wavelet Neural Operator~\cite{tripura2023wavelet,gupta2021multiwavelet} and the Neural Operator on Riemannian Manifolds~\cite{chen2024learning}, respectively.
There are generalizations to nonlinear $\cK_t$ in which its integral kernel $\kappa_t$ depends on the input function $h$ itself \cite{kovachki2021neural}. Regardless of the parametrization used, what is important is that $\cK_t$ is nonlocal \cite{lanthaler2023nonlocality}.

In the context of inverse problems, this chapter later makes use of the \emph{Fourier Neural Operator} (FNO) \cite{li2021fourier,kossaifi2024multigrid}, which corresponds to a neural operator $\Psi_{\mathrm{FNO}}$ of the form \eqref{eqn:no} on the unit hypertorus $\Omega=\T^d\simeq [0,1]^d_{\mathrm{per}}$ with integral operator parametrization
\begin{align}\label{eqn:fno_layer}
	(\cK_th)(x) \defeq \left\{\sum_{k \in \Z^d}\left(\sum_{j=1}^{d_\mathrm{c}} \bigl(P^{(k)}_t\bigr)_{\ell j}\ip[\big]{e^{2\pi \iunit \ip{k}{\slot}_{\R^d}}}{h_j}_{L^2(\mathbb{T}^d;\mathbb{C})} \right) \, e^{2\pi \iunit \ip{k}{x}_{\R^d}}\right\}_{\ell=1}^{d_\mathrm{c}}\,.
\end{align}
This corresponds to enforcing that the matrix-valued function $\kappa_t$ in \eqref{eqn:no_kernel} is a translation-invariant convolution kernel parametrized directly in Fourier space by its Fourier series coefficients $P^{(k)}_t \in \mathbb{C}^{d_{c} \times d_{c}}$ for each $k\in\Z^d$.
A practical discretization of the FNO is usually performed on equally-spaced grids with the Fast Fourier Transform algorithm after truncating the series in \eqref{eqn:fno_layer} to a maximum number of resolved Fourier modes \cite{lanthaler2024discretization,kossaifi2024multigrid}. While this limits the most popular implementations to simple periodic box domains, recent works have extended the FNO to non-periodic domains and complex geometries \cite{de2025extension,li2022fourier,huang2025operator,li2021fourier,li2024physics,lanthaler2023nonlocality,kovachki2021universal,wu2025learning}. Another limitation includes the FNO's handling of spatially local phenomena~\cite{liu2024neural,liu2025enhancing}. Despite these particular drawbacks that deserve future attention, the FNO and its variants have demonstrated vast practical success in scientific applications. Theory for the FNO is also beginning to emerge \cite{kovachki2021universal,lanthaler2024discretization,kovachki2024data, lanthaler2023nonlocality,huang2025operator}.

%%%%%%%%%%%%%%%%%%%%%%%%%%%%%%%%%%%%%%%%%%%%%%%%%%%%%%%%%%%%%%%%%%%%%%%%%%%%%%%%%%%%%%%%%%%%%%%%%%%%%%%%%%%%%%%%%%%%%%%%%%%%%%%%%%%%%%%%%%%%%%%%
\section{Learning inverse problem solvers}\label{sec:ee}
To solve an inverse problem such as~\eqref{eqn:ip_main} is to recover the unknown input $u$ from noisy observations $y = \cG(u) + \eta$, thereby defining a solution operator $y \mapsto u$. Although traditional inverse solvers can approximate this data-to-parameter map with reasonable accuracy, evaluating $y \mapsto u$ even once can be computationally expensive due to the reliance on iterative or sequential algorithms. Direct solvers that are derived from first principles (e.g., the D-bar method for EIT \cite{knudsen2009regularized}) can be more efficient, but they often must be handcrafted for each specific nonlinear inverse problem considered~\cite{mueller2012linear}. This naturally motivates the use of \emph{operator learning} techniques to learn efficient inverse solvers directly from data.

However, a direct, off-the-shelf application of operator learning to inverse problems encounters several immediate challenges. Beyond the need for a sufficiently rich dataset to support the costly offline training step~\eqref{eqn:erm_ol}, the forward map $\cG$ may not be invertible, or if it is, its inverse may be ill-conditioned. As a result, one cannot simply define the target operator to be $\sfF = \cG^{-1}$ without carefully restricting the domain of $\cG$ to a sufficiently small set in order for a continuous inverse to be well-defined. Ensuring that the data distribution assigns probability one to such a set is itself a nontrivial challenge. Moreover, the additive noise model~\eqref{eqn:ol_noise_model} commonly assumed in forward operator learning is inappropriate for the inverse setting. Here, the roles of input and output must be reversed: the inverse map is trained to recover noise-free ground-truth parameters, while the noisy outputs of the forward map now serve as inputs. This reversal leads to a challenging \emph{errors-in-variables} learning problem~\cite{cardot2007smoothing,patel2022error,zhou2019gaussian}.

Despite these challenges, there is growing theoretical and empirical evidence that estimating the inverse map with operator learning is both statistically principled and practically robust~\cite{cao2021choose,de2022deep,li2024physics,nguyen2024tnet,pineda2023deep}. Neural operators can be trained using loss functions that explicitly target desired solution attributes, such as accurate parameter recovery or physical consistency with the forward model. In contrast, plug-in estimators that invert a learned forward map are not goal-oriented and therefore cannot achieve optimal performance; see~\cite[Eqn.~(219), p.~58]{de2022deep} for details. Moreover, even when the underlying inverse problem is severely ill-posed, the average-case nature of the training process tends to suppress the impact of worst-case instability: the extreme inputs and outputs responsible for instability are unlikely under realistic data distributions. End-to-end inverse map learning is especially advantageous in scenarios where data is plentiful and the forward or measurement processes are imperfect or only partially known.

A subset of the research on learning inverse solvers distinguishes itself by blending knowledge of the forward operator into the design of algorithms \cite{li2024physics,nguyen2024tnet,ren2025model,ding2022coupling,cho2025physics,jiang2024resolution,jiao2024solving,zhang2024bilo,jin2022unsupervised}, as advocated by \cite{ghattas2021learning}. This could involve correcting errors in the forward model or inverse operator with learned postprocessing maps \cite{lunz2021learned,hamilton2018deep}.
Other work, such as~\cite{alberti2022continuous,furuya2024can,furuya2023globally}, develops injective architectures that account for the stability of the inverse problem. Additional approaches introduce invertible networks that jointly learn approximate forward and inverse maps while incorporating physical information~\cite{ardizzone2019analyzing,arndt2023invertible,chung2025good,kaltenbach2023semi,long2025invertible,vadeboncoeur2023fully}.
It is also possible to adopt a distribution-matching approach for the inverse map loss component \cite{lingsch2024fuse}, which lifts similar ideas \cite{dunlop2021stability,li2024differential,li2025inverse,yang2023optimal} to the operator learning setting.

This section concentrates on pure data-driven approaches for end-to-end operator learning in linear and nonlinear inverse problems. \Cref{sec:ee_math} formulates the problem of learning inverse maps from noisy paired training data. \Cref{sec:ee_archparam} introduces three neural operator architectures that are tailored to common forms of measurement data that arise in inverse problems. A theoretical analysis develops approximation theory for nonlinear inverse problems and convergence results for linear inverse problems in \cref{sec:ee_theory}. A unifying theme of this theory is that of projection onto finite-dimensional spaces. Finally, \cref{sec:ee_posterior} goes beyond the current emphasis on point estimation by using measure transport to learn the data dependence of posterior distributions.

\subsection{Mathematical formulation}\label{sec:ee_math}
This subsection aims to make the inverse map learning problem more precise. We work under the additive model \eqref{eqn:ip_main}, that is, $y=\cG(u)+\eta$. Suppose that $\cU$ is a Banach space of parameters and $\cY$ is a Banach space of data. Let $\sfd_\cU$ and $\sfd_\cY$ denote the metrics induced by the norms on $\cU$ and $\cY$, respectively.
Both $\cU$ and $\cY$ are assumed to be infinite-dimensional. We view
\begin{align}\label{eqn:forward_map}
	\cG\colon\domain(\cG)\subseteq\cU\to\cY
\end{align}
as a mapping from its domain $(\domain(\cG), \sfd_\cU)$ into $\cY$. The introduction of the set $\domain(\cG)\subseteq\cU$ is necessary in many inverse problems because although $\cG$ may be well-defined on the whole of the Banach space $\cU$, often it is only injective on a smaller set (thus guaranteeing unique solutions of the noiseless inverse problem). This is the case for EIT from \cref{ex:eit}, where even smaller sets are identified such that the restriction of $\cG$ to these sets has a uniformly continuous inverse \cite{harrach2019uniqueness}. Thus, this subsection assumes that $\domain(\cG)$ is chosen such that the inverse map $\Psi^\star$ defined by
\begin{align}\label{eqn:inverse_map}
	\begin{split}
		\Psi^\star\colon\domain(\Psi^\star)\subseteq\cY&\to\cU\\
		\cG(u)&\mapsto u
	\end{split}
\end{align}
is uniformly continuous on its domain
\begin{align}\label{eqn:inverse_map_domain}
	\domain(\Psi^\star)\defeq\cG\bigl(\domain(\cG)\bigr)\,.
\end{align}
We denote the modulus of continuity of $\Psi^\star$ by $\omega^\star$. In \cref{ex:eit} on EIT, $\omega^\star$ is the logarithmic function on the right-hand side of \eqref{eqn:stability_eit_inside_range}. From \eqref{eqn:inverse_map_domain}, we see that $\domain(\Psi^\star)$ equals the range $\image(\cG)$ of the forward map $\cG$ in \eqref{eqn:forward_map}. For nonlinear inverse problems, identification of $\domain(\cG)$ and hence $\domain(\Psi^\star)$ is often a challenging PDE analysis question in its own right.

The goal of end-to-end inverse map learning is to approximate the operator $\Psi^\star$ by some $\widehat{\Psi}$ in the sense that $\widehat{\Psi}(\cG(u)+\eta)\approx u$ for $(u,\eta)$ ranging over an appropriate set. An operator learning framework that seeks to achieve this goal is presented in \cref{sec:ee_math_train}. It is important to appreciate that the true inverse map $\Psi^\star$ in \eqref{eqn:inverse_map} \emph{is only defined} on the range of $\cG$. \Cref{sec:ee_math_noisy} explores the implications of this fact both in theory and in practice.

\subsubsection{Training and evaluating learned inverse maps}\label{sec:ee_math_train}
We work in the setting of \cref{sec:back_ol_learn}, except here the noise plays a crucial role. The simplest training data generation model for inverse map learning is
\begin{align}
	\bigl\{(y_n, u_n)\bigr\}_{n=1}^N\,,\qw u_n\diid\mu \qa y_n=\cG(u_n)+\eta_n\,.
\end{align}
The $\{\eta_n\}$ are i.i.d.~copies of a random variable $\eta$ that models noise or errors. The roles of input and output are reversed from those in \cref{sec:back_ol_learn}. Now the $y_n$ are inputs and the $u_n$ are outputs. Choosing a hypothesis class $\sH$ consisting of operators mapping $\cY$ into $\cU$ and a loss functional $\ell\colon \cU\times\cU\to\R$, we solve
\begin{align}\label{eqn:erm_ip}
	\min_{\Psi\in \sH }\frac{1}{N}\sum_{n=1}^N\ell\bigl(u_n, \Psi(y_n)\bigr)
\end{align}
to obtain an empirical risk minimizer $\widehat{\Psi}\colon \cY\to\cU$. Notice that $\widehat{\Psi}$ is defined on all of $\cY$ even though $\Psi^\star$ is only defined on $\domain(\Psi^\star)\subseteq\cY$ \eqref{eqn:inverse_map_domain}. The hypothesis space $\sH$ could encompass the operator learning architectures from \cref{sec:back_ol_arch}, but also those neural operators tailored to inverse problems that will be introduced in \cref{sec:ee_archparam}. 

Similarly, the same losses $\ell$ from \cref{sec:back_ol_learn} are valid here. However, there is considerably more freedom to choose $\ell$ in the present inverse problem setting. For example, one could choose an $\ell$ with dependence on the forward map $\cG$ to promote desirable problem-adapted structure in the estimator $\widehat{\Psi}$. Moreover, when reconstructing solutions to inverse problems that are discontinuous, such as piecewise constant conductivity fields in EIT or permeability fields in groundwater flow \cite{de2025extension,molinaro2023neural,li2024physics}, it is beneficial to train in a metric on $\cU$ that penalizes sharp transitions less harshly. For instance, the TV norm and its variants are natural candidates, following on its successful introduction in the mathematical imaging community \cite{chan1998total,chan2001active}. Within the class of Lebesgue spaces $L^p$---whose norms are more efficient to compute in a deep learning training loop than TV norms---the $L^1$ distance is best adapted to discontinuous parameter fields \cite{boulle2024mathematical,guerra2025learning,molinaro2023neural,de2025extension}. In analogy with the popular relative $L^2$ loss functions used in forward operator learning, this would involve choosing $X=L^1$ in \eqref{eqn:rel_loss_output}.

A large amount of freedom also exists when choosing how to evaluate the performance of $\widehat{\Psi}$ as a learned inverse solver. In view of the ERM problem \eqref{eqn:erm_ip}, the expected risk
\begin{align}\label{eqn:loss_expected_risk_ip}
	\E_{(u,y)}\bigl[\ell\bigl(u, \widehat{\Psi}(y)\bigr)\bigr]\defeq \E_{(u,\eta)}\bigl[\ell\bigl(u, \widehat{\Psi}\bigl(\cG(u)+\eta\bigr)\bigr)\bigr]
\end{align}
is a natural notion of accuracy. However, we cannot expect \eqref{eqn:loss_expected_risk_ip} to be small unless the noise $\eta$ is small in an appropriate sense. In statistical learning, this is addressed by instead considering the \emph{excess risk}~\cite[Sec.~2.3]{de2023convergence}. Other choices include the noise-free average parameter reconstruction loss
\begin{align}
	\E_{u\sim\mu}\bigl[\ell\bigl(u,\widehat{\Psi}\bigl(\cG(u)\bigr)\bigr)\bigr]
\end{align}
or the average data fitting losses
\begin{align}
	\E_y\bigl[\sfd_\cY\bigl(y, \cG\bigl(\widehat{\Psi}(y)\bigr)\bigr)\bigr] \qor \E_{(u,y)}\bigl[\sfd_\cY\bigl(\cG(u), \cG\bigl(\widehat{\Psi}(y)\bigr)\bigr)\bigr]\,,
\end{align}
and variants thereof.

\subsubsection{Noisy data and the range of the forward map}\label{sec:ee_math_noisy}
We have emphasized that the true inverse map $\Psi^\star$ from \eqref{eqn:inverse_map} is only defined on its domain $\domain(\Psi^\star)=\image(\cG)$. However, estimators of the inverse map, such as the empirical risk minimizer $\widehat{\Psi}$ of \eqref{eqn:erm_ip}, must be defined on a larger set in order to handle perturbed measurement data as input. Indeed, we have that $\cG(u)\in\image(\cG)$, but in general $(\cG(u)+\eta)\notin\image(\cG)$. This challenge is not present in traditional forward operator learning. It is thus important to determine the actual operator that $\widehat{\Psi}$ is really trying to approximate.

Conceptually, a natural approximation $\Psi$---which in some cases will be optimal---of $\Psi^\star$ would take the form $\Psi\defeq \Psi^\star\circ P_{\cG}$, where $P_{\cG}\colon\cY\to\image(\cG)$ is some notion of projection onto $\image(\cG)$. In this way, $\Psi\colon\cY\to\cU$ \emph{extends} $\Psi^\star$, i.e., $\Psi(y)=\Psi^\star(y)$ if $y\in\image(\cG)$. If $y\notin\image(\cG)$, then $P_{\cG}(y)$ is a ``denoised'' version of $y$ in a sense. However, it is practically problematic to construct $P_\cG$ because for many infinite-dimensional nonlinear inverse problems, $\image(\cG)$ is typically a complicated nonconvex set that is difficult to characterize \cite{nickl2023bayesian}. Nonetheless, it is intuitively possible that by minimizing the end-to-end ERM training objective \eqref{eqn:erm_ip}, $\widehat{\Psi}$ is \emph{implicitly} approximating the extension $\Psi$ even though $\Psi$ is not accessible. An interesting research direction is to learn neural operators that \emph{explicitly} project onto $\image(\cG)$ by training on enough data from the range.

\begin{figure}[tb]
	\centering
	\captionsetup{skip=10pt}
	\begin{tikzpicture}[scale=1.0]
		
		% Boxes
		\draw[thick] (-4.5,-3) rectangle (-0.5,3);
		\draw[thick] (1.5,-3) rectangle (5.5,3);
		
		% Labels
		\node at (-2.5,3.3) {Parameter Space};
		\node at (3.5,3.3) {Data Space};
		\node at (-1,2.8) {$\cU$};
		\node at (5,2.8) {$\cY$};
		
		% Curved sets
		\draw[thick, domain=-2:1, smooth, variable=\y] plot ({-4+0.3*(\y+2)^2}, \y);
		\node at (-3.7,-2.4) {$\domain(\cG)$};
		
		\draw[thick, domain=-3:4]
		plot[smooth] coordinates {(1.8,2.1) (2,2) (3.2,1.4) (5,1) (2,-2)};
		\node at (3.4,-2.4) {$\cG\bigl(\domain(\cG)\bigr)=\domain(\Psi^\star)$};
		
		% Points
		\filldraw (-2.1,0.5) circle (1.2pt) node[below] {$u$};
		\filldraw (3.2,1.4) circle (1.2pt) node[above right=-1pt] {$\mathcal{G}(u)$};
		\filldraw (3.0,0.5) circle (1.2pt) node[above] {$y^\delta$};
		\filldraw (-2.5,-0.8) circle (1.2pt) node[below] {$\widetilde{\Psi}^\star(y^\delta)$};
		
		% Arrows
		\draw[arrows = {-Latex[width=5pt, length=10pt]}] (-2.1,0.5) -- (3.2,1.4) node[midway, above, sloped] {$\mathcal{G}$};
		\draw[arrows = {-Latex[width=5pt, length=10pt]}] (3.0,0.5) -- (-2.5,-0.8) node[midway, below, sloped] {$\widetilde{\Psi}^\star$};
		
		% Delta circle
		\draw[dotted,ultra thick] (3.2,1.4) circle (1.3cm);
		\draw[thick] (3.2,1.4) -- (4.2,0.53);
		\node at (3.7,0.75) {$\delta$};
		
		% Spaces
		\node at (-4.6,2.6) {};
		\node at (5.4,-2.6) {};
		
	\end{tikzpicture}
	\caption{Noisy data $y^\delta$ lives outside of the domain $\domain(\Psi^\star)$ of the true inverse map $\Psi^\star\colon \domain(\Psi^\star)\to \cU$. In the illustration, the distance between $y^\delta$ and $\domain(\Psi^\star)$ is less than $\delta$.  A continuous extension $\widetilde{\Psi}^\star\colon\cY\to\cU$ of $\Psi^\star$ is well-defined on noisy data $y^\delta$. As $\delta\to 0$, it holds that $\widetilde{\Psi}^\star(y^\delta)\to \Psi^\star(\cG(u))=u$ with a certain convergence rate. Learned inverse maps $\widehat{\Psi}\colon\cY\to\cU$ conceptually aim to approximate the extension $\widetilde{\Psi}^\star$.
	}
	\label{fig:ip_out}
\end{figure}

The intuition from the preceding thought experiment is instructive. Learned inverse maps $\widehat{\Psi}$ should approximate an \emph{extension} $\widetilde{\Psi}^\star$ of the true inverse $\Psi^\star$. \Cref{fig:ip_out} illustrates this idea. While we cannot expect the extension to exactly project onto the range of the forward map $\cG$, we can hope that the extension inherits the stability of the true inverse map. Few works deliberately address these issues; the papers \cite{pineda2023deep,de2025extension,de2022deep, abhishek2024solving} are exceptions. There are two main analysis strategies in this line of work. The first demands quantitative stability of the trained neural operator \cite{pineda2023deep}. The second strategy puts the onus of stability only on the extended inverse map \cite{abhishek2024solving,de2025extension}.

Consider the setting of \cref{fig:ip_out}. Let $\widetilde{\Psi}^\star\colon\cY\to\cU$ extend $\Psi^\star$. Suppose that $\widehat{\Psi}\colon\cY\to\cU$ has modulus of continuity $\widehat{\omega}$. For a fixed $u\in\domain(\cG)$, we define $y\defeq \cG(u)\in\domain(\Psi^\star)$ and $y^\delta=y+\eta$ with $\sfd_\cY(y^\delta,y)\leq \delta$. To control the approximation error of $\widehat{\Psi}$, the first analysis strategy rests on the following application of the triangle inequality:
\begin{align}\label{eqn:tri1}
	\begin{split}
		\sfd_\cU\bigl(\Psi^\star(y), \widehat{\Psi}(y^\delta)\bigr)&\leq \sfd_\cU\bigl(\Psi^\star(y), \widehat{\Psi}(y)\bigr) + \sfd_\cU\bigl(\widehat{\Psi}(y), \widehat{\Psi}(y^\delta)\bigr)\\
		&\leq \sup_{v \in \domain(\Psi^\star)}\sfd_\cU\big(\widetilde{\Psi}^\star(v),\widehat{\Psi}(v)\bigr) + \widehat{\omega}(\delta)\,.
	\end{split}
\end{align}
The first term in the upper bound \eqref{eqn:tri1} is controlled by how well $\widehat{\Psi}$ approximates $\widetilde{\Psi}^\star$. The second error term is the modulus of continuity of the learned inverse map evaluated at the maximal noise level; this term tends to zero as $\delta\to 0$. While it is a hard problem to quantitatively control Lipschitz or other uniform continuity properties of trained neural networks and neural operators, there do exist specialized architectures that attempt this program \cite{gouk2021regularisation,pineda2023deep,murari2025approximation,anil2019sorting}.

In contrast, the second analysis strategy demands that the extension $\widetilde{\Psi}^\star$ have a modulus of continuity $\widetilde{\omega}^\star$. Typically $\widetilde{\omega}^\star$ bounds $\omega^\star$ from above. An alterative use of the triangle inequality delivers the bound
\begin{align}\label{eqn:tri2}
	\begin{split}
		\sfd_\cU\bigl(\Psi^\star(y), \widehat{\Psi}(y^\delta)\bigr)&\leq \sfd_\cU\bigl(\widetilde{\Psi}^\star(y), \widetilde{\Psi}^\star(y^\delta)\bigr) + \sfd_\cU\bigl(\widetilde{\Psi}^\star(y^\delta), \widehat{\Psi}(y^\delta)\bigr)\\
		&\leq \widetilde{\omega}^\star(\delta) + \sup_{y' \in \mathsf{K}_\delta}\sfd_\cU\big(\widetilde{\Psi}^\star(y'), \widehat{\Psi}(y')\bigr)\,,
	\end{split}
\end{align}
where $\mathsf{K}_\delta$ is a set of perturbed data with distance at most $\delta$ away from $\domain(\Psi^\star)$. The advantage of this approach is that the stability of the learned map $\widehat{\Psi}$ is never needed. However, finer properties of the extension are required \cite{de2025extension}.

One way to decouple the role of regularization from inverse map training is to build regularization into the training data itself \cite{chen2023let}. This method trains approximate inverse maps on data consisting of regularized inverse problem solutions. These data are obtained from classical inverse solvers, e.g., Tikhonov or Bayesian regularization. The insight is that training on pre-regularized data corresponds to learning a regularized inverse map instead of $\Psi^\star$ \cite{chen2023let,nguyen2024tnet}. While this addresses the noisy data problem, the target regularized inverse no longer exactly solves the inverse problem on the range of $\cG$. That is, it does not extend the true inverse. This reflects a tradeoff between accuracy and regularization in the operator learning setting. Similar considerations arise when choosing to define the inverse map outside of its domain as a point estimator of a Bayesian posterior, such as the posterior mean~\cite{marzouk2025discussion}. In practice, trained inverse maps will not extend the true inverse map, regardless of whether the training data is pre-regularized or not.

The need to work with extensions does not arise in forward operator learning. Moreover, traditional inverse solvers handle data outside of the range of $\cG$ by explicitly prescribing regularization. In contrast, the minimization of \eqref{eqn:erm_ip} implicitly encodes regularization~\cite{bishop1995training}. This regularization mechanism of learned inverse solvers is not well understood. Therefore, analyzing and interpreting the properties of regularization in end-to-end learned inverse maps is a vital direction for future research.

\subsection{Architectures for data-to-parameter maps}\label{sec:ee_archparam}
The choice of machine learning architecture for end-to-end learning of inverse maps depends on the parameter space $\cU$ and data space $\cY$. If $\cU$ and $\cY$ are both function spaces, then standard operator learning architectures such as FNO or DeepONet from \cref{sec:back_ol_arch} are appropriate. However, if the data-to-parameter solution operator of the inverse problem maps between more exotic spaces, then additional innovation is required. For simplicity, the current subsection assumes that the parameter of interest is a function belonging to $\cU=\cU(\Omega;\R^{d_{\mathrm{p}}})$. Depending on the inverse problem, the observed data could be represented as a finite vector, a function, a functional, a (linear) operator, or even a probability measure. This subsection explores discretization-invariant architectures for vector-to-function maps in \cref{sec:ee_archparam_fnm}, operator-to-function maps in \cref{sec:ee_archparam_orn}, and distribution-to-function maps in \cref{sec:ee_archparam_nio}.

\subsubsection{Fourier Neural Mappings}\label{sec:ee_archparam_fnm}
Although they might be modeled at a continuum level, finite-dimensional data $y\in\R^J$ for some $J\in\N$ are the actual quantities that are acquired by measurement devices in real applications. Thus, to approximate inverse maps $\Psi^\star\colon\domain(\Psi^\star)\subseteq\R^J\to\cU(\Omega;\R^{d_{\mathrm{p}}})$, we require architectures that map finite vectors to functions. In principle, any encoder--decoder model \eqref{eqn:encode_decode} can be turned into a vector-to-function (VtF) map by removing the encoder. This results in $\Psi=\sD\circ\sA\colon\R^J\to\cU$ for the DeepONet \eqref{eqn:deeponet} and PCA-Net \eqref{eqn:pcanet} upon choosing the approximator $\sA$ to have input dimension $J$. Any other autoencoder method that uses a functional decoder $\sD$ also fits into this framework~\cite{bunker2024autoencoders}.

To produce VtF neural operators, one can append a functional decoder $\sD$ near the beginning of the iteration \eqref{eqn:no} from \cref{sec:back_ol_arch_no}. This is the idea adopted by the Neural Mappings framework~\cite{huang2025operator}, which uses functional decoders $\sD$ of the form 
\begin{align}\label{eqn:vtf}
	\bigl(\sD z\bigr)(x)\defeq\kappa(x)z
\end{align}
for each input vector $z\in \R^{d_\mathrm{c}}$ and point $x\in\Omega$ in the domain. In \eqref{eqn:vtf}, $\kappa\colon\Omega\to\R^{d_\mathrm{c}\times d_\mathrm{c}}$ is a learnable matrix-valued function. Thus, $\sD z\colon \Omega\to\R^{d_\mathrm{c}}$ is a vector-valued function. In particular, Fourier Neural Mappings (FNMs) generalize FNOs---which only map functions to functions---by accommodating VtF or function-to-vector (FtV) maps.  In the former case, these \emph{Fourier Neural Decoders} parametrize the function $\kappa$ in Fourier space by its conjugate symmetric Fourier coefficients $\{P^{(k)}\}_{k\in\Z^d}\subseteq\C^{d_\mathrm{c}\times d_\mathrm{c}}$. Letting $\Omega=\T^d$, this leads to the decoder
\begin{align}\label{eqn:decoder_fnm}
	\bigl(\mathscr{D}z\bigr)(x)=\kappa(x)z= \sum_{k\in\Z^d}\bigl(P^{(k)} z\bigr)\, e^{2\pi\mathrm{i}\ip{k}{x}_{\R^d}}
\end{align}
for $z\in\R^{d_\mathrm{c}}$ and $x\in \T^d$. Then the full VtF FNM $\Psi_{\mathrm{FNM}}\colon \R^{J}\to \cU$ is defined by 
\begin{align}\label{eqn:arch_fnm}
	\Psi_{\mathrm{FNM}}(z)\defeq \bigl(\cQ\circ \sL_T\circ \sL_{T-1}\circ\cdots\circ\sL_2\circ \sD\circ S \bigr)(z)\,,
\end{align}
where $\sD$ is given in \eqref{eqn:decoder_fnm} and $S\colon \R^J\to\R^{d_\mathrm{c}}$ is a shallow neural network; compare to \eqref{eqn:no}, which uses $\sL_1$ instead of $\sD$ and $\cS$ instead of $S$. The $\{\sL_t\}$ are the usual Fourier integral operator layers as in \eqref{eqn:no_layer_nonlinear} and \eqref{eqn:fno_layer}. The map $\cQ$ is as in \eqref{eqn:pw_proj}. FNMs retain the computational benefits and universal approximation properties of FNOs \cite[Sec.~3]{huang2025operator}.

If the measurements in the inverse problem are comprised of both functions and finite vectors, then the extension of FNMs introduced in \cite{bhattacharya2025learning} can be used to handle such multimodal data. The FtV ideas from FNMs \cite{huang2025operator} are also useful when the dimension of the input and output spatial domains do not match, such as in inverse boundary value problems (IBVPs) \cite{de2025extension}. Finally, we note that if $J$ is allowed to vary, it could be more natural to model the inverse operator as a function-to-function map (\cref{sec:back_ol_arch}), operator-to-function map (\Cref{sec:ee_archparam_orn}), or measure-to-function map (\cref{sec:ee_archparam_nio}), depending on the nature of the limiting object as $J\to\infty$. This removes the need to re-train the model every time $J$ changes, which would be the case for VtF architectures.

\subsubsection{Operator Recurrent Neural Networks}\label{sec:ee_archparam_orn}
Operator-to-function maps 
\begin{align}\label{eqn:op_to_func_map}
	\begin{split}
		\Psi^\star\colon\domain(\Psi^\star) \subseteq \cL(\cX;\cW)&\to \cU(\Omega;\R^{d_{\mathrm{p}}})\\
		\Lambda&\mapsto \Psi^\star(\Lambda)
	\end{split}
\end{align}
frequently arise as solution operators of IBVPs. Typically $\cX$ and $\cW$ are spaces of functions on the boundary $\partial\Omega$. The continuous linear operator $\Lambda\in \cL(\cX;\cW)$ is usually a mathematically idealized representation of all possible measurements that one can perform at the boundary. Several concrete IBVPs can be formulated as in \eqref{eqn:op_to_func_map}, including EIT from \cref{ex:eit}, inverse wave scattering based on the Helmholtz equation, optical tomography based on the radiative transport equation, and seismic imaging based on the acoustic wave equation \cite{molinaro2023neural}.

In practice, each single evaluation of the forward map $u\mapsto\Lambda_u$ requires tens, hundreds, or even thousands of PDE solves to probe the boundary and acquire measurements that lead to a discrete representation of $\Lambda=\Lambda_u$. Thus, repeated calls to the forward map can become prohibitively expensive, especially if wrapped around an outer objective. An accurate learned inverse map can substantially reduce the time-to-solution when compared to PDE-constrained, optimization-based, iterative IVBP solvers. At the continuum level, it remains to define operator-to-function architectures.

If linear map $\Lambda$ has a kernel integral operator representation, then viewing its kernel function as input data we return to the familiar function-to-function setting of \cref{sec:back_ol_arch}; see \cite{de2025extension,fan2020solving}, where this approach is taken. A more novel deep learning architecture for representing nonlinear operator functions such as \eqref{eqn:op_to_func_map} is the \emph{Operator Recurrent Neural Network} (ORN-Net) \cite{de2022deep}. This architecture---and existing analysis for it---is currently restricted to a finite-dimensional discretized setting. However, it is straightforward to formally extend the idea of operator recurrence to function spaces by, e.g., following the principles outlined in \cite{berner2025principled}. We do so in the following description of the architecture. 

To this end, we define $\Psi_{\mathrm{ORN}}\colon \cL(\cX;\cW)\to \cU$, which is a neural operator version of ORN-Net, by
\begin{align}\label{eqn:arch_ornn}
	\Psi_{\mathrm{ORN}}(\Lambda)\defeq \bigl(\cQ\circ \sL_T(\slot;\Lambda)\circ \sL_{T-1}(\slot;\Lambda)\circ\cdots\circ\sL_1(\slot;\Lambda)\circ\cS\bigr)(h_0)\,.
\end{align}
Notice that because $\Psi_{\mathrm{ORN}}$ maps a linear operator to a function, the input linear operator $\Lambda\in \cL(\cX;\cW)$ does not appear as input to the first map $\cS$ in the chain of compositions \eqref{eqn:arch_ornn}, but instead \emph{parametrizes} the function-to-function operators $\sL_t(\slot; \Lambda)\colon \cH \to \cH$ acting on the hidden states. We used the notation $\cH=\cH(\Omega;\R^{d_{\mathrm{c}}})$ from \cref{sec:back_ol_arch_no}. The function $h_0$ is an initial hidden state function that serves as a tunable parameter of the architecture and does not depend on the input operator $\Lambda$. As usual, $\cS$ and $\cQ$ are pointwise lifting \eqref{eqn:pw_lift} and projection layers \eqref{eqn:pw_proj}, respectively.

The key innovation of ORN-Net is the parametrization of the hidden layer operators $\sL_t\colon \cH\times\cL(\cX;\cW)\to \cH$. These are given by
\begin{align}
	\begin{split}
		\sL_t(h;\Lambda)&\defeq u_t^{(0)} + \sigma_t\bigl(u_t^{(1)}\bigr)\,,\qw \\
		u_t^{(j)}(x)&\defeq \bigl(\cK_t^{(j)}\Lambda h\bigr)(x) + \bigl(\widetilde{\cK}_t^{(j)} h\bigr)(x) + W_t^{(j)} h(x) + b_t^{(j)}(x)
	\end{split}
\end{align}
for each $j\in\{0,1\}$ and every $h$, $\Lambda$, $t$, and $x$. The $\{\cK_t^{(j)}\}$ and $\{\widetilde{\cK}_t^{(j)}\}$ are nonlocal linear kernel integral operators of the form \eqref{eqn:no_kernel} or \eqref{eqn:fno_layer}. Notice how $\Lambda$ is introduced in the previous display via \emph{operator composition}, i.e., matrix multiplication.
If the activation functions $\{\sigma_t\}$ are set to the identity, then a deep ORN-Net computes operator polynomials and starts to resemble a truncated Neumann series expansion \cite[Sec.~6]{de2022deep}. Thus, the nonlinear version with non-polynomial activation functions generalizes the network's representation capacity to more complicated operator functions.

Extensions to the core architecture include adding memory to the layers instead of settling for the Markovian evolution \eqref{eqn:arch_ornn}; see \cite[Sec.~2.2--2.3]{de2022deep}. The architecture design is further motivated by specific inverse problem structure such as the compositional or iterative nature of the ground truth inverse problem solution map. This leads to inductive biases that promote compactness, sparsity, and low-rankness in the model, which in turn lead to improved statistical generalization bounds. We refer to \cite[Sec.~2--3, 7]{de2022deep} for the details and an application to boundary control of the wave equation.
We note that both the quantitative universal approximation and statistical learning theory for ORN-Net \cite{de2022deep} only handles target functions mapping $\R^{n\times n}$ to $\R^m$ and does not carry over immediately to the infinite-dimensional case. This is an important avenue for future mathematical research. Numerical implementation of ORN-Nets in both finite- and infinite-dimensional settings is also an important future direction to determine the practical benefits of the architecture's multiplicative structure.

In practice, our access to the boundary data map $\Lambda$ may only be through a finite collection of input-output pairs $\{(g_m,\Lambda g_m)\}_{m=1}^{M(\Lambda)}$ obtained by $M=M(\Lambda)\in\N$ operator-function products with some functions $\{g_m\}$. One could then perform a preprocessing step to obtain an integral kernel or matrix representation of $\Lambda$ from the paired boundary data, e.g., via linear regression. Alternatively, one can adopt a measure-centric perspective to avoid the preprocessing step by viewing the paired data $\{(g_m,\Lambda g_m)\}_{m=1}^{M(\Lambda)}$ as a single empirical probability measure (\Cref{sec:back_ip_prob_measure}). This is the subject of the next subsection.

\subsubsection{Neural Inverse Operator}\label{sec:ee_archparam_nio}
Inspired again by IBVPs in which a finite but varying number of function-valued boundary measurements serve as the observed data, we adopt the measure-centric perspective on inverse map learning from \cref{sec:back_ip_prob_measure}. In particular, we describe a general architecture called the \emph{Neural Inverse Operator} (NIO) that is able to approximate measure-to-function maps in a discretization-invariant way \cite{molinaro2023neural}. We then instate the abstract NIO framework in the concrete IBVP of EIT from \cref{ex:eit}.

By now, there are several deep learning methods that can represent probability measure-to-vector operators of the form $\Psi\colon\sP(\R^{d_1})\to\R^{d_2}$. The attention mechanism in modern transformer architectures is known to define a mean-field transport map whose dependence on the underlying input measure leads to such a $\Psi$ \cite{geshkovski2024measure,geshkovski2023mathematical,bach2025learning,huang2024unsupervised}. Of more direct relevance to this subsection are the simpler and more computationally efficient DeepSet~\cite{zaheer2017deep} and PointNet~\cite{qi2017pointnet,qi2017pointnet++} models. These architectures are originally formulated as acting on point clouds in $\R^{d_1}$. The basic idea is to encode the point clouds into a latent space $\R^{d_{\mathrm{c}}}$ by first pushing forward the points with a single neural network, followed by the application of a permutation-invariant aggregation or pooling operation, such as summation \cite{zaheer2017deep} or the entrywise maximum \cite{qi2017pointnet}. Another neural network is used to post-process the latent state into $\R^{d_2}$. For DeepSets, this gives
\begin{align}\label{eqn:deepset}
	\Psi_{\mathrm{DS}}(\mu)\defeq \Psi^{(2)}\biggl(\int_{\R^{d_1}} \Psi^{(1)}(x)\mu(dx)\biggr)\,,
\end{align}
where $\Psi^{(1)}\colon\R^{d_1}\to\R^{d_\mathrm{c}}$ and $\Psi^{(2)}\colon\R^{d_{\mathrm{c}}}\to\R^{d_2}$ are neural networks. In particular, \eqref{eqn:deepset} is valid for genuine probability measures $\mu\in \sP(\R^{d_1})$ instead of just point clouds or, equivalently, empirical measures.

We now generalize DeepSets from finite to infinite dimensions. Let $X$ be an abstract metric space or Banach space and $\mu\in\sP(X)$ be a probability measure over $X$. Consider our usual output function space $\cU=\cU(\Omega;\R^{d_\mathrm{p}})$. In IBVPs, $X$ is typically a space of functions defined on the boundary $\partial\Omega$. Let $\cH$ be a latent space. Choosing operators $\Psi^{(1)}\colon X\to\cH$ and $\Psi^{(2)}\colon\cH\to\cU$ leads to an operator $\Psi\colon\sP(X)\to\cU$ of the form \eqref{eqn:deepset}. An abstract NIO \cite{molinaro2023neural} makes the more constrained choice
\begin{align}\label{eqn:nio}
	\Psi_{\mathrm{NIO}}(\mu)\defeq \Psi_{\mathrm{FNO}}\biggl(\int_{X} \Psi_{\mathrm{DON}}(f)\mu(df)\biggr)\,,
\end{align}
where $\cH\defeq\cU$, $\Psi^{(1)}\defeq \Psi_{\mathrm{DON}}\colon X\to\cU$ is a DeepONet \eqref{eqn:deeponet}, and $\Psi^{(2)}\defeq \Psi_{\mathrm{FNO}}\colon\cU\to\cU$ is an FNO \eqref{eqn:no}. The integral in \eqref{eqn:nio} is a $\cU$-valued Bochner integral that can equivalently be written as the expectation under $\mu$.

Recent work has also proposed measure-to-function architectures by combining standard DeepSets with DeepONets \cite{chiu2024deeposets,prasthofer2022variable,tretiakov2025setonet}. However, the NIO differs from these models because the input measures considered in \cite{chiu2024deeposets,prasthofer2022variable,tretiakov2025setonet} are supported on finite-dimensional spaces. In contrast, the NIO~\eqref{eqn:nio} is the first operator learning architecture that accepts probability measures \emph{supported on functions spaces} as input. Handling two levels of infinite-dimensionality, the first from $X$ and the second from $\sP(X)$, is a notable achievement. In principle, the transformer architecture could also be lifted to infinite-dimensional state spaces and used to solve IBVPs for PDEs, but so far it has been restricted to finite dimensions due to the quadratic complexity of the attention mechanism. Regardless, all such architectures lack strong theoretical guarantees aside from some universal approximation results.

For IBVPs like EIT, we only observe $M=M(\gamma)$ pairs \eqref{eqn:eit_finite_data_model} of Neumann and Dirichlet boundary values for each conductivity $\gamma$. The order of the pairs does not influence what the conductivity is; thus, the inverse map from the pairs to the conductivity is \emph{permutation-invariant}. Modern techniques generalize this permutation-invariance by viewing such mappings as defined on the space of probability measures (\cref{sec:back_ip_prob_measure}). The practical discrete setting, where the data is an empirical measure centered at the pairs of Neumann and Dirichlet data, then becomes a special case of this more general framework.

When applied to IBVPs, the NIO \eqref{eqn:nio} is called \emph{A Measure-theoretic Inverse Neural Operator} (AMINO) \cite{guerra2025learning}. The input to this architecture is a distribution over the joint space $\cX\times\cW$ of paired boundary samples, where $\cX=\cX(\pOmega)$ and $\cW=\cW(\pOmega)$ are boundary function spaces as in \eqref{eqn:op_to_func_map}.
Instead of working with operator-function products as in \cref{sec:ee_archparam_orn}, AMINO processes paired samples from the input measure in a permutation-invariant manner. This results in one-shot inversion directly from the raw data. In contrast, ORN-Net \eqref{eqn:arch_ornn} requires a pre-processing step to obtain the boundary operator $\Lambda$ first.

Instate the EIT setting of \cref{ex:eit}. Then $\Lambda$ is the NtD map $\Lambda_\gamma$. For a distribution $\rho\in\sP(\cX)$ over Neumann boundary conditions, define a forward map
\begin{align}\label{eqn:forward_amino_eit}
	\begin{split}
		\sfF_\rho\colon \domain(\sfF_\rho)\subseteq\cU &\to \sP(\cX\times\cW)\\
		\gamma&\mapsto (\Id,\Lambda_\gamma)\push\rho\,.
	\end{split}
\end{align}
Future work should aim to characterize the distributions $\rho$ for which $\sfF_\rho$ is injective. We nevertheless proceed by formally defining an inverse map $\Psi^\star\colon \sfF_\rho(\gamma)\mapsto\gamma$ for measure-centric EIT. AMINO approximates $\Psi^\star$ via the map $\Psi_{\mathrm{AMINO}}\colon \sP(\cX\times \cW)\to\cU$ defined by $\nu\mapsto \Psi_{\mathrm{AMINO}}(\nu)\defeq\Psi_{\mathrm{NIO}}(\nu)$ as in \eqref{eqn:nio}. Thus, on $\image(\sfF_\rho)$, we have
\begin{align}\label{eqn:amino}
	\Psi_{\mathrm{AMINO}}\bigl(\sfF_\rho(\gamma)\bigr)=\Psi_{\mathrm{FNO}}\biggl(\int_{X} \Psi_{\mathrm{DON}}\bigl((g,\Lambda_\gamma g)\bigr)\rho(dg)\biggr)\,.
\end{align}
Now consider a new realization of the finite collection $\{(g_m,y_m)\}_{m=1}^M$ of noisy input pairs. The $\{g_m\}$ are Neumann boundary values sampled from $\rho$. The $\{y_m\}$ are noisy Dirichlet boundary values as in~\eqref{eqn:eit_finite_data_model}.
By identifying these data with the empirical measure on the left-hand side of~\eqref{eqn:dist_model3_lift}, the action of AMINO on the empirical measure is
\begin{align}\label{eqn:amino_empirical}
	\Psi_{\mathrm{AMINO}}\biggl(\frac{1}{M}\sum_{m=1}^{M}\delta_{(g_m, y_m)}\biggr)= \Psi_{\mathrm{FNO}}\biggl(\frac{1}{M}\sum_{m=1}^{M} \Psi_{\mathrm{DON}}\bigl((g_m,y_m)\bigr)\biggr)\,.
\end{align}
In particular, the trainable parameters of the FNO $\Psi_{\mathrm{FNO}}$ and DeepONet $\Psi_{\mathrm{DNO}}$ in \eqref{eqn:amino_empirical} are \emph{independent} of the number of atoms $M=M(\gamma)$, which may be viewed as a resolution hyperparameter in this measure-centric setting.
During practical training of NIO and AMINO using ERM \eqref{eqn:erm_ip}, additional ideas are used to improve computational cost and generalization accuracy by exploiting the atomic-invariance of the architecture. For example, one can subsample atoms uniformly at random from every input empirical measure in each minibatch of SGD. One can also evaluate the trained AMINO on new point clouds with any number of atoms. These features highlight the power of operator learning for IBVPs.

\subsection{Theoretical results}\label{sec:ee_theory}
Theoretical analysis of end-to-end operator learning for approximating inverse problem solution maps remains relatively scarce. Although there are scattered results for specific nonlinear inverse problems, even the theory for linear inverse problems is far from complete. A key challenge lies in the fact that input data may not necessarily belong to the domain of the inverse map being studied (see \cref{sec:ee_math_noisy}). How to reconcile this issue with the training and optimization procedures used in learning remains an open question.

While new neural operator architectures tailored to the data structures arising in inverse problems have been developed (see~\cref{sec:ee_archparam}), we currently lack a deep understanding of their theoretical properties. Moreover, end-to-end learning remains without a systematic regularization strategy to mitigate the inherent instability of the target inverse map. Current approaches have primarily relied on techniques such as dimension reduction~\cite{pineda2023deep,arridge2024inverse,aspri2020data,chung2025good}, implicit regularization through architecture design or optimizer-induced biases~\cite{de2022deep,dittmer2020regularization}, and training on pre-regularized solutions generated by traditional inverse solvers~\cite{chen2023let}. However, these strategies require further mathematical analysis to rigorously demonstrate their advantages, limitations, and tradeoffs.

This subsection examines theoretical results concerning the approximation of solution operators for inverse problems, highlighting this emerging but critical direction of research. In \cref{sec:ee_theory_eit}, we explore the existence of neural operator approximations in the context of EIT from~\cref{ex:eit}. Although the results are not necessarily constructive, they demonstrate that the specified architectures are well-suited for approximating inverse maps arising in IBVPs such as EIT. Notably, these existence results are general. They hold regardless of data availability or the specific choice of training procedure (e.g., physics-informed losses versus purely data-driven losses). This universality illustrates the flexibility of the theoretical framework. It is compatible with both current and future statistical analyses of different estimators and learning algorithms, as long as the hypothesis space is fixed.

In \Cref{sec:ee_theory_finite}, we shift from EIT to a general deep learning framework for nonlinear inverse problems. This framework builds on the result that, even when the forward operator is infinite-dimensional and nonlinear, and the data are corrupted by noise, it is possible to restrict the forward operator to finite-dimensional subspaces in such a way that the inverse operator becomes Lipschitz continuous. Moreover, there exist neural networks that robustly approximate this projected inverse map in the presence of noise and can be effectively learned using suitably perturbed training data.

Finally, \cref{sec:ee_theory_linear} investigates linear inverse problems in Hilbert spaces. It presents a data-driven regularization-by-projection framework that guarantees convergence of estimated solutions even when the forward operator is unknown or potentially ill-posed. This approach combines input-output training pairs with functional-analytic regularization techniques. By focusing on linear inverse problems, the framework goes beyond existence theorems by demonstrating convergence of estimators in the limit of infinite training data for both noise-free and noisy observations.

\subsubsection{Neural operator approximation theory: A case study on EIT}\label{sec:ee_theory_eit}
We revisit the setting of Calder\'on's problem, i.e., EIT, from \cref{ex:eit}. The goal is to recover the conductivity $\gamma$ from knowledge of the NtD map $\Lambda_\gamma$, subject to the PDE~\eqref{eqn:elliptic_eit}. The extensive mathematical foundation underlying EIT makes it a natural testbed for developing and studying theoretical approximation results in the field of data-driven inverse problems. It is important to note, however, that existing statistical analysis of EIT solvers has largely been limited to the Bayesian framework, which assumes a fixed realization of observed boundary data~\cite{abraham2020statistical}. Under this framework, computationally intensive Bayesian reconstructions must be re-computed for each new realization of observed data.

In contrast, the literature on approximation theory for EIT is significantly more developed~\cite{abhishek2024solving,castro2024calderon,de2025extension}. The central theoretical approach consists of three main steps. First, one shows that the Calder\'on inversion operator maps between subsets of separable Hilbert spaces. This step often involves establishing fine-grained stability estimates for the inverse problem and applying results from PDE theory. Second, one extends the domain of the inverse map to the entire data space, while ideally preserving regularity properties. Third, one characterizes the compactness of the admissible sets of conductivities in suitable function spaces or their inclusion within the support of an appropriate probability measure. These three steps together enable the application of universal approximation theorems for neural operators. In the following, we provide a more detailed discussion of the main results and insights arising from this line of research.

The papers~\cite{abhishek2024solving,castro2024calderon,de2025extension} all formulate the inverse map $\Psi^\star$ of EIT as a mapping between two separable Hilbert spaces. This arises in \cite{abhishek2024solving,castro2024calderon} due to a reliance on encoder-decoder architectures \eqref{eqn:encode_decode} that require use of Hilbert space concepts such as projection and expansion in orthogonal bases. However, the reformulation in Hilbert spaces involves significant effort because $\Psi^\star$ naturally maps a Banach space of continuous linear operators into $L^\infty(\Omega)$, neither of which are Hilbert spaces.

To elaborate further,~\cite{castro2024calderon} develops a generalization of DeepONets from function-to-function mappings to more general operators between separable Hilbert spaces $\cH_1$ and $\cH_2$. This abstraction corresponds to PCA-Net~\eqref{eqn:pcanet}, except with trunk functions $\{\varphi_j\}$ replaced by arbitrary vectors in the output separable Hilbert space $\cH_2$. The paper \cite{castro2024calderon} constrains the architecture even more by enforcing that the encoder $L$, as defined in~\eqref{eqn:pcanet}, projects onto an orthonormal basis (ONB) of the input Hilbert space $\cH_1$ and that the trunk functions $\{\varphi_j\}$ form an ONB of the output Hilbert space $\cH_2$. With this structure, the approach can accommodate Hilbert--Schmidt operators as inputs or outputs while maintaining compatibility with the Hilbertian framework.

We slightly deviate from \cref{ex:eit} by working with Dirichlet-to-Neumann (DtN) maps as input data instead of NtD maps. This choice only affects the exponents of the boundary Sobolev spaces involved and leaves the central methodology unchanged. Define $\cH_2 = L^2(\Omega)$ and $\cH_1 = L^2_\mu(H^{1/2}(\partial\Omega); H^{-1/2}(\partial\Omega))$, where $\mu \in \sP(H^{1/2}(\partial\Omega))$ is a probability measure. That is, the linear DtN maps, which belong to the Banach space $\cL(H^{1/2}(\partial\Omega); H^{-1/2}(\partial\Omega))$ of continuous linear operators, are embedded into the larger Bochner Hilbert space $\cH_1$. We assume that $\cH_1$ is separable. It is a large space and even contains nonlinear discontinuous maps. 

Based on results for the encoder-decoder framework~\cite{lanthaler2022error,castro2023kolmogorov}, an approximation theorem for the inverse map $\Psi^\star$, viewed as a mapping from a subset of $\cH_1$ into $\cH_2$, holds~\cite[Thm.~4.4, p.~786]{castro2024calderon}. This theorem asserts that for every $\ep > 0$, there exists a generalized DeepONet $\Psi \colon \cH_1 \to \cH_2$ such that
\begin{align}\label{eqn:eit_avg_error}
	\int_{\cH_1} \one_{\domain(\Psi^\star)}(\Lambda) 
	\norm{\Psi^\star(\Lambda) - \Psi(\Lambda)}^2_{L^2(\Omega)} \Q(d\Lambda) \leq \ep^2\,.
\end{align}
Here $\one_A$ denotes the indicator function of the set $A$.
The domain $\domain(\Psi^\star)$ consists of DtN maps corresponding to sufficiently regular conductivities that are uniformly bounded above and below. The error is measured on average with respect to a probability measure $\Q$ over the large Bochner space $\cH_1$. The result in~\eqref{eqn:eit_avg_error} relies on a simple, discontinuous zero-extension of the inverse map $\Psi^\star$ outside its domain of definition (recall \cref{sec:ee_math_noisy}). The generalized DeepONet $\Psi$ is then shown to approximate this extended inverse map. The proof only relies on the fact that $\Psi^\star$ is a measurable mapping between two separable Hilbert spaces, which simplifies the required analysis. 

However, there remains an inherent modeling challenge in defining a probability measure supported on the range of the forward map. This range is a highly intricate set that is neither easily parametrized nor straightforwardly related to a distribution over conductivities. Furthermore, the current approach offers no mechanisms to control the approximation when faced with noisy input data. This limitation is somewhat artificial because the map $\Psi^\star$ is uniformly continuous on its domain. Incorporating such finer-grained continuity information could improve the treatment of data outside of the range of the forward map $\gamma \mapsto \Lambda_\gamma$.

To this end, we turn our attention to a quantitative extension approach that does use detailed stability information about the inverse map \cite{de2025extension}. To simplify the presentation, we restrict to a setting in which $d=2$ and $\Omega\subset\R^2$ is the open unit disk. Thus, we identify the boundary $\pOmega$ with the unit torus $\T$. This aligns with \cref{fig:eit_viz} from \cref{ex:eit}.
Next, we consider an infinite-dimensional set $\Gamma \subset L^\infty(\Omega)$ of conductivities identically equal to one outside of a compact subset of $\Omega$, uniformly bounded above and below pointwise a.e., and that have uniformly bounded $L^1(\Omega)$ and TV norms. This set is compact in $L^2(\Omega)$ and contains discontinuous conductivities of practical relevance \cite[Sec.~4.2]{de2025extension}. The compactness of $\Gamma$ enables us to work with uniformly accurate approximations, which interact with stability more favorably than average error estimates such as~\eqref{eqn:eit_avg_error}.

For any $\gamma\in\Gamma$, one can show that the NtD map $\Lambda_\gamma\in \HS(L^2(\T);L^2(\T))$ is a Hilbert--Schmidt operator on $L^2(\T)$ and hence a kernel integral operator. This implies that we can identify $\Lambda_\gamma$ with its integral kernel function $\kappa_\gamma\in L^2(\T\times\T)$. Writing $\domain(\Psi^\star)\defeq \{\kappa_\gamma\condbar\gamma\in\Gamma\}$, we seek an approximation to the inverse map $\Psi^\star\colon \domain(\Psi^\star)\subseteq L^2(\T\times\T)\to L^2(\Omega)$ defined by $\Psi^\star(\kappa_\gamma)\defeq\gamma$. By the Benyamin--Lindenstrauss extension theorem~\cite[Theorem~1.12, p.~18]{benyamini1998geometric}, the fact that $\Psi^\star$ maps between two Hilbert spaces and \eqref{eqn:stability_eit_inside_range} allow us to deduce the existence of a continuous extension with log stability \emph{outside} of $\domain(\Psi^\star)$ \cite[Sec.~3.4, Thm.~3.9]{de2025extension}; recall \cref{fig:ip_in,fig:ip_out}. Specifically, there exists $C>0$, $\al>0$, and $\widetilde{\Psi}^\star\colon L^2(\T\times\T)\to L^2(\Omega)$ such that $\widetilde{\Psi}^\star$ extends $\Psi^\star$ and
\begin{align}\label{eqn:stability_kernel_eit}
	\norm[\big]{\widetilde{\Psi}^\star(\kappa_0)-\widetilde{\Psi}^\star(\kappa_1)}_{L^2(\Omega)}\leq C \, \left(\frac{1}{\log\left(\norm{\kappa_0-\kappa_1}_{L^2(\T\times\T)}^{-1}\right)}\right)^{\al}
\end{align}
for all $\kappa_0\in L^2(\T\times\T)$ sufficiently close to $\kappa_1\in L^2(\T\times\T)$.

The stability estimate \eqref{eqn:stability_kernel_eit} shows that $\widetilde{\Psi}^\star$ is a continuous function-to-function map and hence is amenable to universal approximation theorems for neural operators (\cref{sec:back_ol_arch}). We chose to work with the FNO from \eqref{eqn:no}, \eqref{eqn:no_layer_nonlinear}, and \eqref{eqn:fno_layer}; other architectures would also be valid.
Let $\sfE_\sigma$ be a compact subset of $L^2(\T\times\T)$ with the property that $\sup_{\eta\in\sfE_\sigma}\norm{\eta}_{L^2(\T\times\T)}\leq \sigma$.
Then by following the approach~\eqref{eqn:tri2} from \cref{sec:ee_math_noisy}, we can establish a noise-robust approximation result \cite[Sec.~4.4, Thm.~4.8]{de2025extension} which states that for every $\ep>0$, there exists an FNO $\Psi\colon L^2(\T\times\T) \to L^2(\Omega)$ such that
\begin{align}\label{eqn:fno_eit_bound}
	\sup_{(\gamma,\eta)\in\Gamma\times\sfE_\sigma}\norm[\big]{\gamma-\Psi(\kappa_\gamma+\eta)}_{L^2(\Omega)}\leq \ep + C \, \biggl(\frac{1}{\log(\sigma^{-1})}\biggr)^{\al} .
\end{align}
Although this result demonstrates robustness to noise in the sense that the second term on the right-hand side of \eqref{eqn:fno_eit_bound} goes to zero as $\sigma\to 0$, \emph{it converges slowly}. In other words, the inherent instability of EIT is still present. In particular, it is sufficient to take the noise level $\sigma=O(\exp(-\ep^{-1/\al}))$ exponentially small to match the size of the FNO approximation error $\ep$.

Related work~\cite{abhishek2024solving} exploits quantitative continuity properties of $\Psi^\star$ in the more practical complete electrode model (CEM) of EIT~\cite{dunlop2016bayesian,mueller2012linear}. In the CEM model, measurements are limited to voltages corresponding to currents injected at a fixed finite number $M$ of electrodes located along the boundary. This results in a discretized NtD map $\Lambda_\gamma^M \in \R^{M \times M}$. Instead of working with infinite-dimensional sets of admissible conductivities as in~\cite{castro2024calderon,de2025extension}, another option is to restrict attention to a finite-dimensional manifold $\Gamma \subset L^\infty(\Omega)$ of piecewise analytic conductivities. These conductivities are assumed to be identically equal to one outside of a compact subset of the domain $\Omega$ and uniformly bounded above and below pointwise a.e.~\cite{abhishek2024solving}. Crucially, the finite-dimensionality of this set permits the derivation of \emph{Lipschitz stability} results for the inverse map~\cite{alberti2019calderon,alberti2022infinite,harrach2019uniqueness}, as opposed to logarithmic stability results~\cite{alessandrini1988stable}.

By approximating the EIT inverse map using a DeepONet~\eqref{eqn:deeponet} adapted to linear operator inputs, the following result holds~\cite[Thm.~3.1, p.~7]{abhishek2024solving}. Let $L_M$ be the map $L_M \colon \Lambda_\gamma \mapsto \Lambda_\gamma^M - \Lambda_1^M$. For any $\varepsilon > 0$, there exists a number of electrodes $M$ and a DeepONet $\Psi$ of the form~\eqref{eqn:deeponet} (with $L = L_M$, viewed as a function of $L_M(\Lambda_\gamma)$ instead of $\Lambda_\gamma$) such that
\begin{align}\label{eqn:shift_eit_approx_thm}
	\sup_{\gamma \in \Gamma} \, \norm[\big]{\gamma - \Psi\bigl(\Lambda_\gamma^M - \Lambda_1^M\bigr)}_{L^2(\Omega)} \leq \varepsilon\,.
\end{align}
The proof relies on formulating $\Psi^\star$ as a mapping from a space of Hilbert--Schmidt operators into $L^2(\Omega)$, using Dugundji's extension theorem~\cite{dugundji1951extension} to continuously extend the domain of the map to the entire Hilbert space, and leveraging forward and inverse stability estimates to argue that the image of $\Gamma$ under the parameter-to-data forward map is a compact set.

The DeepONet $\Psi$ in \eqref{eqn:shift_eit_approx_thm} receives as input the raw measurements $\Lambda_\gamma^M \in \R^{M \times M}$ shifted by the background NtD map $\Lambda_1^M$ corresponding to $\gamma \equiv 1$. The use of shifted NtD maps is for technical convenience but is not fundamentally necessary, as demonstrated in~\cite{de2025extension}. Regardless, these matrix inputs are treated as vectors in an $M^2$-dimensional Euclidean space rather than as matrices in a Banach algebra of linear operators. It would be interesting to replace fully connected neural networks with ORN-Nets~\eqref{eqn:arch_ornn} in the branch net of $\Psi$, as these may be better suited for matrix-based data~\cite{de2022deep}; recall \cref{sec:ee_archparam}.

Due to Lipschitz stability arising from the finite dimensionality of $\Gamma$, it should further be possible to prove Lipschitz robustness to noisy NtD maps in~\eqref{eqn:shift_eit_approx_thm} by adapting the ideas in~\cite{de2025extension}. However, this favorable robustness will degrade as the dimensionality of $\Gamma$—and hence the number of measurements $M$—increases~\cite[Sec.~4.1, Thm.~5]{alberti2022infinite}. One would only expect to recover the logarithmic robustness shown in \eqref{eqn:fno_eit_bound}.

To summarize, the three preceding lines of research critically rely on inverse stability estimates for $\Psi^\star$ when viewed as a map between subsets of two separable Hilbert spaces~\cite{abhishek2024solving,castro2024calderon,de2025extension}. Achieving such bounds requires leveraging the PDE structure inherent to EIT. Extending this framework beyond EIT entails developing appropriate stability estimates and extensions for the particular inverse map under consideration. However, these results must be carefully verified on a case-by-case basis; not all PDE models possess the necessary structure to enable such analysis.
Consequently, an important direction for future research is to apply the current approximation theory framework to other PDE-based nonlinear inverse problems, such as inverse wave scattering, optical tomography, or seismic imaging~\cite{molinaro2023neural}.

\subsubsection{Finite-dimensionalizing nonlinear inverse problems}\label{sec:ee_theory_finite}
Although most analysis of deep learning for nonlinear inverse problems proceeds on a case-by-case basis, one notable exception is the work of Pineda~and~Petersen~\cite{pineda2023deep}. By transforming the conceptually infinite-dimensional nonlinear inverse problem into a finite-dimensional one, it is possible to reduce the setting to that of finite-dimensional unknown parameters and Lipschitz continuous inverse maps defined on a finite-dimensional data manifold. This transformation relies on ideas from \cite{alberti2022infinite}. The precise connection between the finite and infinite problems is explained in \cite[Sec.~1.2.1~and~Sec.~2]{pineda2023deep}. A key point is that this procedure is valid for a broad class of nonlinear inverse problems. The main ingredients involve restricting the forward map to an \emph{a priori} finite-dimensional set of admissible parameters and then composing a family of finite-rank ``measurement'' linear operators with the restricted forward map.

By carrying out this program, the original problem reduces to one about function approximation in finite dimensions. The following notable result holds true~\cite[Thm.~1.1]{pineda2023deep}. Write $d_{\mathrm{m}}$ for manifold dimension and $d_{\mathrm{a}}$ for ambient dimension. Let $J\in\N$ and let $\psi^\star\colon\domain(\psi^\star)\subset \R^{d_{\mathrm{a}}}\to\R^J$ be a bounded Lipschitz mapping such that $\domain(\psi^\star)$ is a $d_{\mathrm{m}}$-dimensional relatively compact manifold satisfying certain technical conditions \cite[Thm.~3.9]{pineda2023deep}. There exists $\sigma_0>0$ such that the following holds. If $\sigma^2\leq \sigma_0^2$, then for every $\ep>0$, there exists a ReLU neural network $\psi$ with $O(\ep^{-d_{\mathrm{m}}})$ parameters such that
\begin{align}\label{eqn:robust_nn}
	\E_{\eta\sim \normal(0,\sigma^2 I_{d_{\mathrm{a}}})}\biggl[\sup_{y\in\domain(\psi^\star)}\norm[\big]{\psi^\star(y)-\psi(y+\eta)}_{\R^{J}}^2\biggr]\leq \ep^2 + C^\star d_{\mathrm{m}}\sigma^2\,.
\end{align}
The constant $C^\star$ depends on the Lipschitz constant $\Lip(\psi^\star)$ of $\psi^\star$, the manifold $\domain(\psi^\star)$, and the maximal variance $\sigma_0^2$, but not on the error tolerance $\ep$. 
The result as stated in \eqref{eqn:robust_nn} follows from \cite[Thm.~3.9]{pineda2023deep} and the first choice of triangle inequality~\eqref{eqn:tri1} mentioned in \cref{sec:ee_math_noisy}; compare to \eqref{eqn:fno_eit_bound}, which uses the second triangle inequality~\eqref{eqn:tri2}. As usual, the proof involves extending $\psi^\star$ away from its domain and approximating the extension instead. The extra Lipschitz-type robustness of the network to Gaussian noise is due in part to its construction based on ideas from finite element approximation.

We observe from \eqref{eqn:robust_nn} that the network approximation $\psi$ has a favorable noise-damping property in the following sense. First notice that $d_{\mathrm{m}}\sigma^2=(d_{\mathrm{m}}/d_{\mathrm{a}})\E\norm{\eta}^2_{\R^{d_{\mathrm{a}}}}$. A Lipschitz estimate for a general function on $\R^{d_{\mathrm{a}}}$ would only provide an average-case bound of order $\E\norm{\eta}^2_{\R^{d_{\mathrm{a}}}}$ for the perturbation $\eta$. Here, we instead see that the fine properties of the particular data manifold produce the dampening factor $d_{\mathrm{m}}/d_{\mathrm{a}}\leq 1$. This factor could be significantly smaller than $1$ if the manifold is low-dimensional compared to the ambient space.
Although \eqref{eqn:robust_nn} corresponds only to a quantitative existence theorem, it is shown that a neural network approximation can actually be found through ERM \eqref{eqn:erm_ip} over noisy training pairs generated from the forward map \cite[Thm.~4.3]{pineda2023deep}. This trained network is accurate on average and further enjoys Lipschitz stability properties similar to those of the idealized, unattainable network $\psi$ from \eqref{eqn:robust_nn}.

To connect back to the motivation of inverse problem solvers, it is instructive to interpret the map $\psi^\star$ as a finite-dimensional representation of an underlying infinite-dimensional inverse map $\Psi^\star\colon \domain(\Psi^\star)\subset \cY\to\cU$ between Banach spaces \cite[Cor.~1]{alberti2022infinite}. Indeed, it is shown that the preceding theory applies to several nonlinear inverse problems corresponding to Darcy flow, Euler--Bernoulli beam, and nonlinear Volterra integral models \cite[Sec.~2]{pineda2023deep}.
One limitation of the general approach put forth by \cite{pineda2023deep} is that much of the infinite-dimensional structure, including instability as manifested by a generally non-Lipschitz modulus of continuity of $\Psi^\star$, is hidden in the finite-dimensionalization procedure and associated Lipschitz constants. Indeed, the neural network approximation is intertwined with the number of measurements $d_{\mathrm{m}}$. The estimation becomes harder as the parameter dimension $J$ grows because $\Lip(\psi^\star)$---and hence the number of measurements $d_{\mathrm{m}}$---tends to infinity in this limit \cite{alberti2022infinite}. The network also must change each time $J$ or $d_{\mathrm{m}}$ changes, which violates the desirable discretization invariance principles of operator learning. Adding additional encoder-decoder structure to the approximation could reinstate consistency (see~\cref{sec:back_ol_arch_onet}). Moreover, since the theory is only valid for finite-dimensional parameters, this precludes natural infinite-dimensional parameter models such as bounded transformations of Gaussian random fields \cite{dunlop2016bayesian,stuart2010inverse}. Overcoming such challenges while maintaining the wide applicability of \cite{pineda2023deep} would be a substantial contribution to the field of data-driven inverse problems.

\subsubsection{Finite-dimensionalizing linear inverse problems}\label{sec:ee_theory_linear}
It is fruitful to analyze operator learning for linear inverse maps because linearity enables the uncovering of deep theoretical insights that are widely applicable to a whole family of inverse problems, unlike the nonlinear setting of the previous subsections. Such a program was carried out successfully for linear forward operator learning \cite{de2023convergence,mollenhauer2022learning}. Analogous results for linear inverse problems are still missing, possibly due to the difficulty of handling the unboundedness of the inverse operator in conjunction with the training data and noise models.

One exception where analysis can be pushed through is the method of data-driven regularization-by-projection \cite{arridge2024inverse,aspri2020data}. 
The analysis adopts the framework of functional-analytic regularization \cite{engl2015regularization} and works in the setting that the forward operator is unknown. Specifically, the forward map is only accessible through input-output training data. The theory establishes convergence guarantees in the limit of infinite training data, both for noise-free and noisy observations. We now survey such results.

Let $\cU$ and $\cY$ be real separable Hilbert spaces. Let $\cA\in\cL(\cU;\cY)$ be an injective compact linear operator. For example, $\cA$ could represent the Radon transform from the imaging sciences or a (possibly non-symmetric) linear solution operator of a linear PDE. Suppose we are given a training dataset of noiseless labeled input-output pairs
\begin{align}\label{eqn:training_data_pairs}
	(u_n, y_n)\,,\qw y_n=\cA u_n\qf n\in\{1,\ldots, N\}\,.
\end{align}
Based on the $N$ pairs from \eqref{eqn:training_data_pairs}, our goal is to find an estimator $\cR_N\colon\cY\to\cU$ of the unbounded operator $\cA^{-1}\colon \image(\cA)\subset\cY \to\cU $ such that for any $y\in\image(\cA)$ and any $y^\delta\in\cY$ such that $\norm{y-y^\delta}_{\cY}\leq\delta$, it holds that the reconstruction error
\begin{align}\label{eqn:goal_linear_ip}
	\norm[\big]{\cA^{-1}y-\cR_Ny^\delta}_{\cU}
\end{align}
is small if the noise level $\delta>0$ is sufficiently small.
That is, we have \emph{learned to solve} the inverse problem of finding $u$ from imperfect observations $y^\delta\approx \cA u$. Of course, we must choose $N=N(\delta)$ as a function of $\delta$ in order for the display in \eqref{eqn:goal_linear_ip} to converge to zero as $\delta\to 0$ \cite{engl2015regularization,clason2020regularization}. Thus, the sample size $N$ plays the role of a regularization hyperparameter.

To define our choice of learned inverse map $\cR_N\in\cL(\cY;\cU)$, let the linear span of the training data be written as
\begin{align}
	\cU_N\defeq\Span\{u_n\}_{n=1}^N \qa \cY_N\defeq\Span\{y_n\}_{n=1}^N\,.
\end{align}
Further define $P_{\cU_N}\colon \cU\to \cU_N$ to be the orthogonal projection operator onto $\cU_N$. Similarly, define $P_{\cY_N}\colon \cY\to \cY_N$ to be the orthogonal projection operator onto $\cY_N$.
The \emph{data-driven regularization-by-projection estimator} is defined by
\begin{align}\label{eqn:reg_lin_est}
	\cR_N\defeq \cA^{-1}P_{\cY_N}\,.
\end{align}
The motivation for the formula \eqref{eqn:reg_lin_est} is that the action of the forward map $\cA$ on the span $\cU_N$ of the input training data is completely determined by the data. It is then possible to deduce that $\cR_Ny$ corresponds to the minimum norm solution of the \emph{input-projected} linear operator equation $y=\cA P_{\cU_N} u$. Then $\cR_N$ in \eqref{eqn:reg_lin_est} is seen to be the Moore--Penrose pseudoinverse of $\cA P_{\cU_N}$ \cite[Thm.~4, p.~5]{aspri2020data}.

To evaluate the map $y\mapsto \cR_N y$ only using the training data \eqref{eqn:training_data_pairs} and not $\cA$ itself, we appeal to the Gram--Schmidt orthonormalization procedure. Let $\{\varphi_j\}_{j=1}^N$ be the orthonormal basis of $\cY_N$ obtained from the Gram--Schmidt algorithm applied to $\{y_n\}_{n=1}^N$ (already assumed to be linearly independent). Then
\begin{align}
	P_{\cY_N}=\sum_{j=1}^N \varphi_j\otimes_\cY \varphi_j\,,
\end{align}
where $(a\otimes_\cY b)h\defeq\ip{b}{h}_{\cY}a$ is the outer product operator on $\cY$. It follows that
\begin{align}
	\cR_N=\sum_{j=1}^N \bigl(\cA^{-1}\varphi_j\bigr)\otimes_\cY \varphi_j \eqdef \sum_{j=1}^N v_j\otimes_\cY \varphi_j\,.
\end{align}
The vectors $\{v_j\}_{j=1}^N$ can be computed offline using $v_1\defeq u_1/\norm{y_1}_{\cY}$ and the iteration
\begin{align}
	v_j=\frac{u_j-\sum_{k=1}^{j-1}\ip{\varphi_k}{y_j}_{\cY} \, v_k}{\norm{y_j-P_{\cY_{j-1}}y_j}_{\cY}}\qfa j>1\,.
\end{align}
Thus, we can compute $y\mapsto \cR_Ny=\sum_{j=1}^N\ip{\varphi_j}{y}_{\cY}\, v_j$ in linear time complexity without access to $\cA$ or its inverse. These expensive offline steps are similar to the ones used in PCA-Net \eqref{eqn:pcanet} and other model reduction methods \cite{bhattacharya2021model}.

To obtain convergence guarantees for $\cR_N$, we require some assumptions on the training data \eqref{eqn:training_data_pairs}. Assume that the sequence $\{\norm{u_n}_{\cU}\}$ is bounded away from zero and infinity. For all $N$, suppose that the vectors $\{u_n\}_{n=1}^N$ are linearly independent. Furthermore, assume that
\begin{align}
	\bigl\{(u_n,y_n)\bigr\}_{n=1}^{N+1}=\bigl\{(u_n,y_n)\bigr\}_{n=1}^{N} \cup \bigl\{(u_{N+1}, y_{N+1})\bigr\}
\end{align}
are nested and that $\overline{\cup_{N\in\N}\cU_N}=\cU$. Notice that the training data need not be generated according to a probability distribution; they could be deterministic and highly correlated. For example, one could select the $\{u_n\}$ to be a problem-adapted orthonormal basis or construct them according to some greedy algorithm.
Under further technical assumptions on the compatibility between $\cA$, the linear span $\cU_N$, and a fixed $y\in\image(\cA)$ \cite[Sec.~2.2.1]{arridge2024inverse}, it holds that $\cR_Ny $ converges both weakly and strongly to $\cA^{-1}y$ in the Hilbert space $\cU$ as $N\to\infty$ \cite[Thms.~9, 11, 15]{aspri2020data}. We can also handle the noisy data case $y^\delta\not\in\image(\cA)$. Optimizing a standard bias-variance tradeoff \cite[Thm.~46, p.~33]{clason2020regularization} by choosing $N=N(\delta)\to\infty$ according to the requirement
\begin{align}
	\delta\sqrt{N(\delta)}\sup_{n\in\{1,\ldots, N(\delta)\}}\frac{1}{\abs{\ip{\varphi_n}{y_n}_{\cY}}}\to 0\qas\delta\to 0\,,
\end{align}
it then holds under the aforementioned technical assumptions that
\begin{align}\label{eqn:reg_prog_converge}
	\lim_{\delta\to 0}\, \sup_{\widetilde{y}\in B_\delta(y)} \norm[\big]{\cA^{-1}y- \cR_{N(\delta)}\widetilde{y}}_{\cU}=0 \qfa y\in\image(\cA)
\end{align}
by \cite[Thm.~17, p.~17]{aspri2020data}. In \eqref{eqn:reg_prog_converge}, we write $B_\delta(y)\defeq \set{z\in\cY}{\norm{y-z}_{\cY}\leq \delta}$.

The convergence result \eqref{eqn:reg_prog_converge} is pointwise in the ``input space'' $\image(\cA)$ and uniform over perturbations. Because $\cA^{-1}$ is not continuous, we cannot expect a version of \eqref{eqn:reg_prog_converge} that is uniform over the unit ball of $\image(\cA)$ to be valid \cite[Thm.~4.9, p.~36]{clason2020regularization}. However, comparing to \cref{eqn:eit_avg_error,eqn:fno_eit_bound,eqn:robust_nn} which \emph{are uniform} over the domain of a nonlinear inverse map, it is of interest to determine meaningful compact source sets or probability measures over parameter space that lead to uniform or average case convergence, respectively, of $\cR_{N(\delta)}$ at the level of operators, possibly with respect to new accuracy metrics.

The data-driven regularization-by-projection method can be viewed both as learning the forward operator---via $\cA\approx\cA P_{\cU_N}$---or as learning the inverse operator---via $\cA^{-1}\approx \cA^{-1}P_{\cY_N}$. Both viewpoints can be unified in the special case that the training data arise from the SVD of $\cA$ \cite[Remark~6, p.~6]{aspri2020data}. The current approach is but one of many ways to finite-dimensionalize linear inverse problems by discretizing parameters or data \cite{aspri2021data}. Dual least squares and explicitly regularized variational problems also enjoy convergence guarantees \cite[Secs.~4--5]{aspri2020data}. Beyond deterministic data with span-increasing and linear independence constraints, other choices of projection sets could include the leading eigenspaces of convenient ``prior'' operators that are compatible with $\cA$ \cite{boulle2024operator}. This could help circumvent the technical difficulties associated with the Gram–Schmidt procedure.

The preceding methodologies can all be viewed as simple ways to actively select training data \cite{boulle2023elliptic,gao2024adaptive,guerra2025learning,subedi2024benefits,satheesh2025picore}.
Nevertheless, the development of a theoretical foundation for the i.i.d.~random training data setting is still essential. A natural counterpart to deterministic regularization-by-projection in this stochastic setting is an empirical PCA projection-based inversion method. A sharp comparative analysis could help clarify the benefits of active versus passive data design strategies for linear inverse operator learning.

\subsection{Beyond point estimation: Data-to-posterior maps}\label{sec:ee_posterior}
Compared to data-to-parameter map learning, the literature on learning data-to-posterior maps in the context of infinite-dimensional Bayesian inverse problems (\cref{sec:back_ip_prob_bayes}) is much sparser \cite[Chp.~6]{bach2024inverse}. For a fixed observation $y$, there is abundant work on how to compute the posterior probability measure $\mu^y$ given by \eqref{eqn:posterior}, e.g., by sampling or density estimation methods. However, such traditional inference methods are relatively slow. This cost is magnified if the application requires the repeated calculation of the posterior for many different measurements, which is the case in sequential data assimilation \cite{law2015data} or Bayesian optimal experimental design \cite{huan2024optimal}. By \emph{data-to-posterior operator learning}, we refer to the problem of approximating the map $y\mapsto\mu^y$. This approach effectively solves entire families of Bayesian inverse problems parametrized by the observed data without the need for re-training or repeated optimization solves. Although the upfront training or data generation cost may be high, this offline cost is amortized if the trained map is evaluated multiple times and each evaluation is fast (e.g., near real-time). Similar considerations hold for the prior-to-posterior map $\mu\mapsto \mu^y$, or even the combined Bayesian solution map $(\mu,y)\mapsto \mu^y$, although these are outside the scope of the present chapter.

Some work within the fields of likelihood-free \cite{sisson2018handbook} and simulation-based inference \cite{cranmer2020frontier} can be interpreted within the amortized data-to-posterior framework. For instance, one can build efficient surrogate models for the forward map simulator appearing in the likelihood function to accelerate the Bayesian inference loop. There is other research on data-amortized conditional score-based diffusion models for infinite-dimensional linear inverse problems \cite{baldassari2023conditional,schneider2024unconditional}. Instead, this subsection focuses on \emph{measure transport} methods that directly target the parametric dependence of the posterior on the observed data.

To this end, one approach is to learn data-to-posterior maps with a hypernetwork \cite{zhao2024functional}. The hypernetwork $\Psi$ maps the data $y$ to the parameters $\theta=\Psi(y)$ of another neural network $\mathsf{F}$. We view $\mathsf{F}$ as an approximate transport map from the prior $\mu$ to the posterior $\mu^y$. That is, $\mu^y\approx \mathsf{F}(\slot;\theta)\push\mu$. The map $\mathsf{F}$ only depends on $y$ through its parameters $\theta=\Psi(y)$. The training of $\Psi$ can be performed in an unsupervised way by optimizing the Kullback--Leibler (KL) divergence, but this requires knowledge of the likelihood to evaluate the loss. The overall method is a ``meta'' conditional normalizing flow and has been applied to inverse problems on function spaces \cite{zhao2024functional}. Since the hypernetwork training task is difficult in practice, one should later fine tune the meta-conditional model or use it as a preconditioner for a fixed data realization if higher accuracy is required. \Cref{sec:reg_prior_math} contains related discussion.

An alternative methodology is that of triangular transport maps \cite{el2012bayesian}. This approach falls under the broader purview of computational measure transport \cite{marzouk2016sampling}, as described in \cref{sec:back_ip_prob_transport}. In \eqref{eqn:min_div}, which we recall is the minimization problem $\min_{\sfT\in\cT}\sfd(\nu, \sfT\push\rho)$, the triangular framework takes $\sfd$ to be the Kullback--Leibler (KL) divergence and $\cT$ to be a class of triangular maps with a particular parametrization (e.g., polynomial or wavelet expansions) and structural constraints such as monotonicity \cite{baptista2024representation}. This construction is inherently finite-dimensional, which precludes its scalability to infinite-dimensional inverse problems. One solution involves a relaxation to \emph{block triangular transport maps}, which are well-defined between infinite-dimensional Banach spaces and naturally amortize posterior sampling with respect to the observation data \cite{baptista2024conditional}. Our discussion is restricted to a useful but specific case of a more general setup \cite{hosseini2025conditional}.

Let $\cU$ be the parameter state space and $\cY$ the data space, both assumed to be separable Banach spaces. Suppose there is an unknown joint probability measure $\nu\in\sP(\cY\times\cU)$ over the joint data-parameter space. For instance, the data likelihood model \eqref{eqn:ip_main} and the prior measure $\mu\in\sP(\cU)$ together define such a $\nu$. We wish to solve the minimum divergence estimation problem \eqref{eqn:min_div} such that the class $\cT$ consists of block triangular maps of the form
\begin{align}\label{eqn:block_triangular}
	\sfT(y,u)=
	\begin{pmatrix*}[l]
		y\\
		\Psi(y,u)
	\end{pmatrix*}
\end{align}
for $(y,u)\in\cY\times\cU$ and $\Psi\colon \cY\times\cU\to\cU$. The map $\Psi$ parametrizes the class $\cT$. The nonlinear function \eqref{eqn:block_triangular} is called block triangular due to the variable dependence in each block coordinate that ensures the derivative of $\sfT$ takes values as a block lower-triangular linear operator. We further take $\rho\defeq \nu_{\cY}\otimes \rho_\cU\in\sP(\cY\times\cU)$ to be an independent product reference measure, where $\nu_\cY$ is the $\cY$-marginal of $\nu$ and $\rho_\cU$ is arbitrary. A typical choice is $\rho_\cU=\mu$ for a Bayesian inverse problem with prior $\mu$ \cite[Sec.~4]{hosseini2025conditional}. The form \eqref{eqn:block_triangular} is useful for the following reason. If $\sfT\push\rho=\nu$, then 
\begin{align}\label{eqn:match_conditionals}
	\Psi(y,\slot)\push\rho_\cU=\nu(\slot\condbar y)
\end{align}
for $\nu_\cY$-a.e.~$y\in\cY$ \cite[Thm.~2.4, p.~875]{baptista2024conditional}. That is, we can extract the desired conditional distribution $\nu(\slot\condbar y)$ from the subcomponent $\Psi$ of the full transport map $\sfT$. Moreover, global minimizers $\sfT^\star$ to \eqref{eqn:min_div} over the class \eqref{eqn:block_triangular} exist and achieve $\sfT^\star\rho=\nu$ \cite[Remark~2.11, p.~878]{baptista2024conditional}.

These theoretical underpinnings help motivate the monotone generative adversarial network (M-GAN) architecture \cite[Sec.~3]{baptista2024conditional}, which is a powerful amortized solver for Bayesian inverse problems. M-GAN takes $\sfd$ to be a min-max approximation of the Wasserstein-$1$ distance. For PDE-based inverse problems such as Darcy flow, the transport map $\Psi$ is parametrized using PCA-Net \eqref{eqn:pcanet}. Practically, we only have access to $\nu$ through a finite training dataset of $N$ pairs of i.i.d.~joint samples. Then $\nu$ in \eqref{eqn:min_div} is replaced by the empirical measure $\nu^N$ corresponding to these samples. To avoid overfitting to the training data, the reference $\rho$ is also replaced by its $ M$-sample approximation, with $M$ possibly different from $N$. After training, the M-GAN estimator $\widehat{\Psi}$ satisfies
\begin{align}\label{eqn:match_conditionals_bayesian}
	\widehat{\Psi}(y,\slot)\push\mu\approx \mu^y\,.
\end{align}
Numerical experiments demonstrate the excellent performance of M-GAN on a wide range of Bayesian inverse problems \cite{baptista2024conditional,hosseini2025conditional}.

Although there are sample complexity bounds \cite{schreuder2021statistical} and some approximation theory for measure transport \cite{baptista2025approximation}, at least for fixed $y$, a general error analysis for the performance of learned data-to-posterior operators is still missing. In the block triangular transport setting where training is done at the level of joints \eqref{eqn:min_div}, one important open question is to understand how errors in joint distribution approximation propagate to errors in the desired conditionals or posteriors. A satisfactory answer will be highly dependent on the choice of statistical divergence $\sfd$ used to quantify such errors.

%%%%%%%%%%%%%%%%%%%%%%%%%%%%%%%%%%%%%%%%%%%%%%%%%%%%%%%%%%%%%%%%%%%%%%%%%%%%%%%%%%%%%%%%%%%%%%%%%%%%%%%%%%%%%%%%%%%%%%%%%%%%%%%%%%%%%%%%%%%%%%%%
\section{Learning prior distributions and regularizers}\label{sec:reg}
Learning prior information or regularization terms directly from data has become a powerful alternative to hand-crafted regularization in solving inverse problems~\cite{lunz2018adversarial,afkham2021learning,haltmeier2023regularization,chirinos2024learning,mukherjee2024data,tan2024unsupervised,zhang2025learning,soh2019learning,leong2025optimal}. Traditionally, priors were imposed using explicit penalty terms or probabilistic assumptions; however, modern approaches instead learn these priors from examples of signals or measurements. This allows one to capture complex structures or statistical correlations in data that are difficult to express analytically, especially in high-dimensional or ill-posed settings. The learned prior can be used either in Bayesian inference or to guide iterative algorithms.
This section focuses on two complementary paradigms for learning prior distributions and regularizers: approaches based on empirical Bayes-type prior learning in \cref{sec:reg_prior} and denoising-based regularization (e.g., plug-and-play methods) in \cref{sec:reg_pnp}. 

\subsection{Prior distribution learning}\label{sec:reg_prior}
In the Bayesian formulation, solving an inverse problem involves specifying a prior distribution $\mu$ over the unknown parameter $u$ and a likelihood function $\sfL^y(u) = \exp(-\Phi(u;y))$ that models the relationship between the data $y$ and the parameter $u$ (\cref{sec:back_ip_prob_bayes}). The likelihood represents the forward and measurement processes as well as uncertainties, while the prior encodes our assumptions or background knowledge about $u$ before observing the data.

Bayesian inference combines the prior $\mu$ and likelihood $\sfL^y$ using Bayes' theorem to obtain the posterior distribution $\mu^y$: 
\begin{align}\label{eq:Bayes}
	\mu^y(u) \propto  \sfL^y(u) \, \mu(u)\,,
\end{align}
where, for convenience in this section, we abuse notation by writing $\mu^y(u)$ and $\mu(u)$ for the density of the posterior and the prior measures, respectively.
The posterior represents the updated distribution of $u$ after incorporating the data $y$. However, Bayesian inference is expensive. From a computational efficiency perspective, one goal that does not require probing the full posterior distribution is to extract a point estimate of $u$ that best aligns with the posterior. A common approach is MAP estimation \cite{sanz2023inverse}, which maximizes the posterior probability via  
\begin{align}\label{eq:MAP}
	u_{\mathrm{MAP}} \defeq \argmax_{u\in\cU} \, \mu^y(u)\,.
\end{align}  
Equivalently, using the negative log likelihood potential $\Phi$, MAP estimation minimizes the energy functional given by: 
\begin{align}\label{eq:MAP2}
	\cE(u \condbar y) \defeq -\log  \sfL^y(u) - \log \mu(u) = \Phi(u;y) - \log \mu(u)\,.
\end{align}  
From this optimization perspective, regularization can be interpreted as imposing a prior through the term $R(u)\defeq -\log \mu(u)$. In classical inverse problems, handcrafted regularizers such as quadratic penalties for smoothness or sparsity-promoting norms play the role of $R$. For example, Gaussian priors lead to Tikhonov-type regularizers, while a Laplace prior aligns with sparsity-promoting regularizers. 

Note that in infinite dimensions, the densities appearing in the preceding formulas \eqref{eq:Bayes},~\eqref{eq:MAP},~and~\eqref{eq:MAP2} require careful interpretation \cite{dashti2017}. For certain Bayesian inverse problems on Banach spaces, an infinite-dimensional analog of Bayes' rule holds (recall~\cref{sec:back_ip_prob_bayes} and Eq.~\ref{eqn:posterior}). This permits an analysis similar to that of~\cref{eq:Bayes,eq:MAP,eq:MAP2}, but remains valid in infinite dimensions. Although such extensions cannot rely on probability densities, natural generalizations of MAP estimators \eqref{eq:MAP} exist based on small-ball probabilities; see~\cite{dashti2013map,klebanov2023maximum}.

The fact that the prior acts like a regularization term can also be seen from the variational form of Bayes' theorem~\cite{bach2024inverse,trillos2020bayesian}. Rather than viewing Bayes' rule as a formula for conditioning, one can interpret the posterior $\mu^y$ as the minimizer of a variational problem over the metric space of probability measures. Specifically, the posterior distribution satisfies
\begin{equation}\label{eq:var_Bayes}
	\mu^y   = \argmin_{\nu \in \sP(\cU)} \Bigl\{\mathbb{E}_{u \sim \nu}[\Phi(u;y) ] + \KL{\nu}{\mu} \Bigr\}\,,
\end{equation}
where $\KL{\nu}{\mu}$ is the KL divergence between $\nu$ and the prior $\mu$. 
Equation~\eqref{eq:var_Bayes} requires the absolute continuity of the posterior with respect to the prior and follows from the basic properties of KL divergence.
In this formulation, the KL divergence serves as a regularization term that penalizes deviations from the prior, while the expected negative log-likelihood promotes fidelity with the observed data. The choice of prior $\mu$ directly determines the structure and strength of the regularization term $\KL{\nu}{\mu}$.

In modern settings, data-driven approaches aim to learn the prior $\mu$ directly from training samples. Most of these approaches fall under the \emph{empirical Bayes} classification~\cite{robbins1992empirical,efron2012large} in which the prior $\mu$ is not modeled as random itself, unlike \emph{hierarchical Bayesian} methods~\cite{teh2010hierarchical}.
By modeling complex, high-dimensional distributions of real-world data, prior learning methods offer greater flexibility compared to traditional handcrafted priors. Machine learning techniques, particularly deep generative models~\cite{ruthotto2021introduction}, have been particularly effective in capturing such distributions. Examples include energy-based models, generative adversarial networks, variational autoencoders, and score-based diffusion models~\cite{arridge2019solving,song2021solving,dimakis2022deep,bohra2023bayesian,scarlett2023theoretical,zhao2023generative}. These learned priors provide a data-driven means to enforce regularization in the Bayesian framework and provide state-of-the-art reconstructions in many inverse problems.

\Cref{sec:reg_prior_math} introduces popular regularizers arising from the combination of variational inference with deep learning. To state concrete theoretical results, \cref{sec:reg_prior_theory} first develops an analysis of the effect that learned priors have on downstream posteriors. Specializing to infinite-dimensional linear inverse problems, the subsection then presents an analysis of optimal Tikhonov and spectral regularizers within the nonparametric statistical framework of \cref{sec:back_ip_prob_stat}.

\subsubsection{Mathematical formulation and architecture choices}\label{sec:reg_prior_math}
This subsection focuses on unsupervised neural network-based regularization techniques. In this setting, one assumes access to an i.i.d.~dataset 
\begin{align}\label{eqn:prior_iid}
	\{u_n\}_{n=1}^N \sim \mu_\mathrm{data}^{\otimes N}\,,
\end{align} 
but not paired observation data. The goal is to learn a prior $\mu$ or its parametric representation from this collection of representative ``ground‐truth'' samples distributed according to the inaccessible data distribution $\mu_\mathrm{data}$. The measure $\mu_\mathrm{data}$ is rarely Gaussian. But if it is Gaussian, then the problem becomes equivalent to \emph{covariance estimation} \cite{al2025covariance}, which is one of the simplest examples of unsupervised operator learning. Even in the non-Gaussian case, estimating the covariance operator of $\mu_\mathrm{data}$ from samples \eqref{eqn:prior_iid} is a nontrivial task in high or infinite dimensions. We focus on finding $\mu_\mathrm{data}$ itself. It is also possible, but more challenging, to handle the case when only indirect observations of the model \eqref{eqn:prior_iid} are available \cite{akyildiz2024efficient,vadeboncoeur2023fully}; see also \cref{sec:back_ip_prob_measure} for a measure-centric approach to this distributional inversion problem.

Let $\cZ$ be a finite-dimensional latent space. Most state-of-the-art methods represent the prior as the pushforward under a neural network or neural operator
\begin{align}
	\sfT_\theta\colon \mathcal Z \to \mathcal U
\end{align}
that maps samples $z$ from a latent distribution $\rho$---typically a standard isotropic Gaussian---to the parameter space. That is, the generative model
\begin{align}
	u = \sfT_\theta(z) \sim \mu\qa \mu =\mu_\theta = (\sfT_\theta)\push \rho\,.
\end{align}
Once an appropriate $\sfT_\theta$ is trained, this transport map framework (\cref{sec:back_ip_prob_transport}) converts the inverse problem $y=\cG(u)$ from \eqref{eqn:ip_deterministic} into that of finding $z$ from
\begin{align}\label{eqn:ip_latent}
	y = \cG\bigl(\sfT_\theta(z)\bigr)\,.
\end{align}
Thus, the inference is now over the latent space $\mathcal Z$. Different methods vary in the choice of architecture for $\sfT_\theta$ and the design of the loss function, which together determine the learning preferences for encoding the prior. We now survey several possible options.

\paragraph*{Maximum likelihood and invertible transport maps.}
When $\sfT_\theta$ is invertible (e.g., a normalizing flow), the density of the induced prior $\mu_\theta$ can be obtained via the change-of-variables formula as
\begin{align}\label{eqn:prior_density_map}
	\mu_\theta(u) = \rho\bigl(\sfT^{-1}_\theta (u)\bigr)
	\abs[\big]{\det \bigl(\nabla \sfT^{-1}_\theta(u) \bigr)}\,,
\end{align}
where $\nabla \sfT_{\theta}^{-1}$ denotes the Jacobian matrix of the inverse transformation.
The parameters $\theta \in \Theta$ can then be optimized using maximum likelihood estimation,
\begin{align}\label{eqn:mle_nf}
	\max_{\theta\in\Theta} \frac{1}{N}\sum_{n=1}^N \log \mu_\theta(u_n)\,,
\end{align}
by substituting in the density expression~\eqref{eqn:prior_density_map}.

Normalizing flows (NFs) provide a flexible class of invertible parametrizations $\sfT_\theta$ that make the optimization problem \eqref{eqn:mle_nf} tractable by exploiting the triangular or auto-regressive structure of the Jacobian determinant; see, for example,~\cite{rezende2015variational,dinh2017density,kobyzev2020normalizing}. Related continuous-time generalizations, known as \emph{neural ordinary differential equation} (ODE) flows, define $\sfT_\theta$ as the finite time flow map of a learned ODE and interpret the log-determinant term through the instantaneous trace of the Jacobian~\cite{chen2018neural,grathwohl2018ffjord}.

The preceding constructions of $\sfT_\theta$, however, rely fundamentally on the finite-dimensional change-of-variables formula \eqref{eqn:prior_density_map} and Lebesgue probability densities on $\mathbb{R}^d$. Extending such density-based techniques to infinite-dimensional settings, such as function spaces arising in PDE-based inverse problems, is nontrivial and restrictive because Gaussian reference measures on Banach or Hilbert spaces are typically mutually singular~\cite{stuart2010inverse}. Nevertheless, there are promising recent efforts to generalize normalizing flows to function space~\cite{zhao2024functional}. This functional normalizing flow approach aims to transport Gaussian process priors through invertible transformations while maintaining measure-theoretic rigor, potentially bridging finite-dimensional NF methods with infinite-dimensional Bayesian inverse problem formulations. See also the discussion in \cref{sec:ee_posterior}.

\paragraph*{Autoencoders} 
Architectures based on \emph{autoencoders} consist of an encoder--decoder pair 
\begin{align}
	E_\vartheta\colon \cU \to \mathcal Z\qa
	\sfT_\theta\colon \mathcal Z \to \cU
\end{align}
trained jointly to reconstruct input data. This leads to the optimization problem
\begin{align}\label{eq:ae_loss}
	\min_{(\theta, \vartheta)} \frac{1}{N} \sum_{n=1}^N 
	\norm[\big]{\sfT_\theta\bigl(E_\vartheta(u_n)\bigr) - u_n}_\cU^2\,.
\end{align}
Depending on the setting, $(E_\vartheta, \sfT_\theta)$ may represent a FtV encoder and a VtF decoder pair (recall \cref{sec:ee_archparam_fnm}), or remain
entirely finite-dimensional. In either case, the formulation \eqref{eq:ae_loss} only relies on the Banach
space norm on $\cU$, making it applicable to both discrete and functional representations.
One can also consider regularizing the encoded latent variables in the objective function. This naturally motivates a Gaussian
reference $\rho \in \sP(\mathcal Z)$ that is consistent with the role of latent distributions in
functional autoencoders~\cite{bunker2024autoencoders}. After training, this reference is pushed forward through
the decoder to define the induced distribution $\mu_\theta = (\sfT_\theta)\push \rho$, which serves as a learned prior generative model over the space $\cU$.

\paragraph*{Variational autoencoders}
\emph{Variational autoencoders} (VAEs) extend the deterministic framework of
autoencoders by introducing a probabilistic generative model.
While a standard autoencoder learns deterministic encoder and decoder mappings
$E_\vartheta \colon \cU \to \mathcal Z$ and $\sfT_\theta \colon \mathcal Z \to \cU$
that minimize a reconstruction loss, e.g.,~\eqref{eq:ae_loss}, a VAE treats both mappings as
stochastic maps taking value in spaces of probability measures, i.e., Markov kernels. In particular, the encoder outputs a \emph{distribution}
$q_\vartheta(\slot \mid u)$ over latent variables. The decoder specifies a
conditional distribution $p_\theta(\slot \mid z)$ over reconstructions.

VAEs can be formulated from a \emph{joint-measure and transport} perspective, which is compatible with both finite- and infinite-dimensional settings~\cite[Sec.~13.6]{bach2024inverse}. In this formulation, one defines a joint probability measure
\begin{align}
	p_\theta(dz,du) = p_\theta(du \mid z)\, \rho(dz)
\end{align}
on the product space $\mathcal Z \times \mathcal U$.
Here, $\rho$ is a fixed reference measure on the latent space $\mathcal Z$ (typically $\rho = \normal(0,I)$ in finite dimensions) and $p_\theta(\slot \mid z)$ is a parametrized \emph{decoder operator}
that maps each latent $z$ to a probability measure on~$\mathcal U$. This construction induces a stochastic generative model
\begin{align}\label{eqn:vae_sample}
	z \sim \rho\qa u \sim p_\theta(\slot \mid z)
\end{align}
that specifies the distribution of $u$ through the composition of a latent draw and a conditional sampling step.

In practice, the conditional measures $p_\theta(\slot \mid z)$ are often taken to be Gaussian,
\begin{align}
	p_\theta(\slot \mid z) = \normal\bigl(m_\theta(z), \Sigma_\theta(z)\bigr)\,, 
\end{align}
where the mean $m_\theta \colon \mathcal Z \to \mathcal U$ and covariance operator $\Sigma_\theta \colon \mathcal Z \to \mathcal L(\mathcal U;\mathcal U)$ are represented by neural networks or neural operators. This viewpoint interprets the decoder as a (possibly stochastic) operator that maps latent representations to probability measures over functions or parameters, enabling discretization-invariant modeling of infinite-dimensional data.

To infer latent representations from data, one introduces an \emph{encoder operator} $q_\vartheta(\slot \mid u)$ mapping $u \in \mathcal U$ to distributions on~$\mathcal Z$ and defines a corresponding variational joint measure
\begin{align}
	q_\vartheta(dz,du) = q_\vartheta(dz \mid u)\, \nu(du)\,, 
\end{align}
where $\nu$ denotes the (empirical version of the) data measure $\mu_{\mathrm{data}}$ on $\mathcal U$. The encoder is also often taken to be Gaussian in practice.

The goal of VAE training is to bring the model joint $p_\theta(dz,du)$ and the variational joint $q_\vartheta(dz,du)$ into alignment. This can be achieved by minimizing a divergence such as
\begin{align}
	\KL{q_\vartheta}{p_\theta}
\end{align}
over $(\vartheta,\theta)$, which compares the two joint measures directly. The minimization implicitly balances reconstruction of $u$ through the decoder operator and regularization of the encoder against the latent reference measure $\rho$.
Recent work such as~\cite{bunker2024autoencoders} instates the joint-measure-based VAE framework in infinite-dimensional function spaces. However, such infinite-dimensional VAEs necessitate strong assumptions on the data distribution $\mu_{\mathrm{data}}$ in order for the KL divergence between the joint measures to be finite~\cite[Secs.~2--3]{bunker2024autoencoders}. Samples from the trained VAE may be obtained through \eqref{eqn:vae_sample}.

In finite-dimensional settings, training a VAE can be interpreted as maximizing the
\emph{evidence lower bound} (ELBO), which provides a tractable lower bound on the
marginal log-likelihood $\log p_\theta(u)$.
Here $p_\theta(u)$ denotes the model evidence, which is given by marginalizing over the latent
variable via
\begin{align}
	p_\theta(du) \defeq \E_{z \sim \rho}\bigl[p_\theta(du \mid z)\bigr]\,.
\end{align}
The ELBO takes the form
\begin{align}
	\begin{split}\label{eq:elbo}
		\log p_\theta(u)
		&\geq \mathrm{ELBO}(\theta, \vartheta; u) \\
		&\defeq
		\E_{z \sim q_\vartheta(\slot \mid u)}
		\bigl[\log p_\theta(u \mid z)\bigr]
		- \KL[\big]{q_\vartheta(\slot \mid u)}{\rho}\,.
	\end{split}
\end{align}
The first term encourages the decoder to reconstruct $u$ accurately from samples of the
latent variable $z$, while the second term regularizes the encoder by keeping $q_\vartheta(\slot \mid u)$ close to the reference distribution~$\rho$. After minimizing the ELBO, \eqref{eqn:vae_sample} delivers samples from the trained VAE.

\paragraph*{Generative adversarial networks}
\emph{Generative Adversarial Networks} (GANs) are trained to learn a generator
\begin{align}
	\sfT_\theta\colon \mathcal{Z} \to \cU
\end{align}
that induces a distribution
\begin{align}
	\mu_\theta \defeq (\sfT_\theta)\push\rho
\end{align}
given a known latent distribution $\rho\in\sP(\cZ)$. This is the same setup as before. The point of departure is the training procedure. GANs are trained with an adversarial loss, where the adversary is a discriminator
\begin{align}
	D_\vartheta\colon \cU \to [0,1]    \,.
\end{align}
The original GAN objective is a \emph{minimax game} between the generator and discriminator. This game reads
\begin{align}\label{eqn:gan}
	\min_\theta \max_\vartheta \Bigl\{ \mathbb{E}_{u \sim \mu_{\mathrm{data}}} \bigl[\log D_\vartheta(u)\bigr] + \mathbb{E}_{z \sim \rho} \bigl[\log\bigl(1 - D_\vartheta(\sfT_\theta(z))\bigr)\bigr]\Bigr\}\,,
\end{align}
where $\mu_{\mathrm{data}}$ is the training data distribution \eqref{eqn:prior_iid}. The objective~\eqref{eqn:gan} seeks to minimize the
Jensen--Shannon (JS) divergence between the model distribution
$\mu_\theta = (\sfT_\theta)\push\rho$ and the data distribution
$\mu_{\mathrm{data}}$~\cite{goodfellow2014generative}. The discriminator tries to assign $1$ to real data samples and $0$ to generated ones, while the generator attempts to fool the discriminator by producing samples for which $D_\vartheta(\sfT_\theta(z)) \approx 1$. Variants such as the \emph{Wasserstein GAN}~\cite{arjovsky2017wasserstein} replace the Jensen--Shannon divergence in the GAN loss with an approximation of the dual form of the Wasserstein-$1$ distance. This leads to the minimax game
\begin{align}\label{eqn:wgan}
	\min_\theta \max_{\vartheta\in\set{{\vartheta'}}{\Lip(D_{{\vartheta'}})\leq 1}} \Bigl\{ \mathbb{E}_{u \sim \mu_{\mathrm{data}}} \bigl[D_\vartheta(u)\bigr] - \mathbb{E}_{z \sim \rho} \bigl[D_\vartheta\bigl(\sfT_\theta(z)\bigr)\bigr]\Bigr\} \,.
\end{align}
The $1$-Lipschitz constraint $\Lip(D_\vartheta)\leq 1$ can be enforced using a gradient penalty or weight clipping \cite{gouk2021regularisation}. The Wasserstein GAN tends to exhibit improved stability. The expectations in the preceding displays \eqref{eqn:gan} and \eqref{eqn:wgan} are replaced by empirical averages during training in practice.
After training, the generator $\sfT_\theta$ induces a prior that approximates $\mu_{\mathrm{data}}$.

Recent work has extended GANs to infinite-dimensional setups. In~\cite{rahman2022generative}, Wasserstein GANs were extended to function spaces through the \emph{Generative Adversarial Neural Operator} (GANO), in which the generator and discriminators are formulated as neural mappings~\cite{huang2025operator}: the generator acts as a VtF map producing samples $u \in \cU$, while the discriminator serves as a FtV map evaluating functions to scalar scores. Training uses a Wasserstein-$1$ objective with a discretization-invariant Lipschitz constraint.  A complementary theoretical framework for adversarial learning on infinite-dimensional Banach spaces was provided in~\cite{asatryan2023convenient}, showing convergence of the JS divergence under smoothness and positivity assumptions on the relevant Radon--Nikodym derivatives. A related work on Wasserstein GANs in Banach space is \cite{adler2018banach}.

\paragraph*{Explicit prior or score learning}
Finally, some methods seek to learn an explicit, parametrized regularization functional
$R_\theta\colon\cU\to\R $ that approximates the negative log prior density from \eqref{eq:MAP2}. 
Models of this type are often referred to as \emph{energy-based models} (EBMs), since they represent probability measures through an unnormalized energy function:
\begin{align}\label{eqn:proportional_density}
	\mu_\theta(u) \propto \exp\bigl(-R_\theta(u)\bigr)\,.
\end{align}
In approaches based on maximum likelihood estimation, one fits a density model~\eqref{eqn:proportional_density} to a dataset $\{u_n\}_{n=1}^N \sim \mu_{\mathrm{data}}^{\otimes N}$ by solving
\begin{align}\label{eqn:mle_ebm}
	\max_{\theta\in\Theta} \frac{1}{N}\sum_{n=1}^N \log \mu_\theta(u_n) = \min_{\theta\in\Theta} \frac{1}{N}\sum_{n=1}^N \Bigl[R_\theta(u_n) + \log Z_\theta\Bigr]\,,
\end{align}
where $Z_\theta$ is the normalization constant implied in \eqref{eqn:proportional_density}. The only difference from \eqref{eqn:prior_density_map}~and~\eqref{eqn:mle_nf} is the parametrization of the Lebesgue density of $\mu_\theta$. In practice, computing $Z_\theta$ is often intractable, which necessitates approximate training schemes such as contrastive divergence~\cite{carreira2005contrastive} or score matching~\cite{hyvarinen2005estimation}. Once trained, the resulting $R_\theta$ provides an explicit regularizer suitable for use in downstream inverse problem solvers. Although the Lebesgue density parametrization \eqref{eqn:proportional_density} is only valid in finite dimensions, it is easy to generalize to infinite dimensions by defining the density with respect to a reference probability measure $\rho$ such as a proper Gaussian measure. This leads to
\begin{align}\label{eqn:proportional_density_inf}
	\mu_\theta(du) \propto \exp\bigl(-R_\theta(u)\bigr)\rho(du)
\end{align}
and an analogous optimization problem to \eqref{eqn:mle_ebm}.

Alternatively, score-based models---such as denoising score matching and diffusion models \cite{song2021score,song2021solving}---bypass the need to evaluate or approximate the normalizing constant $Z_\theta$ by instead learning the score function
\begin{align}\label{eqn:score_log_density}
	u\mapsto \bigl(\nabla \log \mu\bigr)(u)
\end{align}
directly with neural networks or neural operators~\cite{baldassari2023conditional,franzese2025generative,schneider2024unconditional}. 
We stress that defining the score is considerably more delicate in function space settings because probability measures on infinite-dimensional Banach spaces lack a Lebesgue density. One must instead identify a generalization of the score that remains well defined without relying on the log-density gradient formula~\eqref{eqn:score_log_density}. In~\cite{lim2025score}, the score is formulated explicitly in function space as the \emph{logarithmic derivative} of a perturbed measure in Cameron--Martin space directions. Sampling is performed via function-valued Langevin dynamics. Alternatively, \cite{pidstrigach2024infinite} introduces a \emph{conditional expectation} representation that plays the role of the score. Both infinite-dimensional score formulas agree with \eqref{eqn:score_log_density} if the Banach space is $\R^d$.

While score-based models do not yield an explicit density for $\mu$, they enable sampling and approximate inference through stochastic dynamics. Consequently, they define a learned prior implicitly through its score or generative process rather than an explicit functional form.

\subsubsection{Theoretical results}\label{sec:reg_prior_theory}
Mathematical guarantees for the \emph{accuracy} of data-driven priors and regularizers in the infinite-dimensional setting, which is the primary focus of this chapter, remain somewhat limited~\cite{arridge2019solving}. For instance, in the variational framework of \eqref{eqn:reg_opt}, much of the existing literature focuses on analyzing the convergence or computational complexity of optimization schemes used to find the regularizer or solve the inverse problem, rather than directly addressing the quality of the inversion itself~\cite{hauptmann2024convergent,leong2025optimal,soh2019learning}. Furthermore, there is a strong preference in existing work for the functional-analytic regularization framework, which combines worst-case analysis with deterministic, bounded noise models~\cite{engl1996regularization}.

In contrast, this subsection surveys theoretical work that has a more probabilistic flavor. In the context of prior learning for Bayesian inverse problems, we study the stability of the posterior distribution with respect to perturbations in the prior. By developing quantitative estimates for the modulus of continuity of the prior-to-posterior map, we can control posterior accuracy---that is, the accuracy of the Bayesian inverse problem solution---by the error incurred by the approximate prior. Due to its generality, this stability approach is widely applicable to various choices of data-driven prior learning architectures or algorithms, including those highlighted in \cref{sec:reg_prior_math}. However, from the perspective of accuracy, such prior-to-posterior bounds can be quite crude because all of the dependence on the likelihood (i.e., the forward map, observed data, and noise level) is hidden in the stability constants.

To achieve possibly sharper results, we then turn our attention to the concrete \emph{average-case} analysis of learned regularization in infinite-dimensional linear inverse problems. The statistical aspects of this problem setting---primarily the training data distribution and noise distribution---are seen to greatly impact inversion accuracy.

\paragraph*{Stability of prior learning}
The quantitative robustness of Bayesian inference has been a subject of extensive study for several decades~\cite{basu1998stability,basu2000uniform,diaconis1986consistency,owhadi2017qualitative}. Recent advancements in the theoretical analysis of Bayesian inverse problems and generative modeling have reignited interest in this topic~\cite{sprungk2020local,garbuno2023bayesian,trevisan2023wide,cvetkovic2025upper}. A key question of relevance to this chapter is the robustness of the posterior distribution with respect to perturbations in the prior distribution. Such prior perturbations can arise from randomized approximations~\cite{dunbar2025hyperparameter,nelsen2021random} or finite element approximations~\cite{girolami2021statistical} of Gaussian process priors, for example. We now overview results of this type.

Recall the setting of Bayesian inverse problems from \cref{sec:back_ip_prob_bayes}. We are concerned with stability estimates for the nonlinear prior-to-posterior map $\mu\mapsto \mu^y$ from \eqref{eqn:posterior} for some fixed $y$. The results we cover are of the type
\begin{align}\label{eqn:prior_perturb_general_bound}
	\sfd(\mu^y,\nu^y)\leq \frac{C(\mu,\nu;y)}{\min(Z_\mu^y, Z_\nu^y)}\sfd_0(\mu,\nu)\,,
\end{align}
where $\sfd$ quantifies closeness of posteriors $\mu^y$ and $\nu^y$, $\sfd_0$ quantifies closeness of priors $\mu$ and $\nu$, $Z_\mu^y$ and $Z_\nu^y$ are the posterior normalizing constants, and $C(\mu,\nu;y)$ is a constant that depends on the priors (typically only through their moments) and the data $y$. See \Cref{fig:stab_prior_learning} for an illustration. Most works take $\sfd=\sfd_0$ as the TV, Hellinger, or Wasserstein-$1$ distances~\cite{sprungk2020local}. There are generalizations to the case $\sfd\neq \sfd_0\equiv \sfd_0^y$, where $\sfd_0^y$ is a \emph{data-dependent} distance function \cite{garbuno2023bayesian}. With respect to the local Lipschitz constants, slightly tighter upper bounds than \eqref{eqn:prior_perturb_general_bound} exist, depending on the choice of $\sfd$. Recently, the sharpness of \eqref{eqn:prior_perturb_general_bound} has been studied from the perspective of lower bounds \cite{cvetkovic2025upper}. These results reveal that the increasing sensitivity as the posterior becomes more concentrated, i.e., as the evidence $Z_\mu^y$ from \eqref{eqn:posterior} approaches zero, is inevitable under certain statistical divergences.

\begin{figure}[tb]
	\centering
 	\captionsetup{skip=10pt}
	\begin{tikzpicture}[scale = 1.2]
		\draw[smooth cycle,tension=1] plot coordinates{(-1,0) (-0.5,1.5) (2,1.6) (2,0)};
		\draw[smooth cycle,tension=1] plot coordinates{(3.5,0) (4.3,2) (7,1.9) (6.6,0)};
		\fill  
		(0,0.1) circle (1.5pt) node[right] {$\nu$} 
		(-0.25,0.6) circle (1.5pt) node[above ]{$\mu$}  
		(4.8,0.98) circle (1.5pt) node[below] {$\mu^y$} 
		(5.3,0.75) circle (1.5pt) node[below] {$\nu^y$};
		\path[->, thick](1, 2) edge [out=45,in=135] node[above] {prior-to-posterior map}(4, 2);
		\draw[blue,ultra thick] (0,0.1)--(-0.25,0.6);
		\draw[red,ultra thick] (4.8,0.98)--(5.3,0.75);
		\node[red] at (5.55,1.05) {\small $\sfd(\mu^y,\nu^y)$};
		\node[blue] at (0.4,0.5) {\small $\sfd_0(\mu,\nu)$};
		\draw[black] (1.26,1.23) node[above]{\small space of priors};
		\draw[black] (5,-0.3) node[above]{\small space of posteriors};
		\draw [black,thick,fill=red, opacity=0.2] (5.31,0.9) circle (0.8);
		\draw [black,thick,fill=blue, opacity=0.2] (0.2,0.6) circle (0.8);
	\end{tikzpicture}
	\caption{The stability of prior learning concerns the relationship between the closeness of the posteriors, measured by $\sfd(\mu^y, \nu^y)$, and the closeness of the priors, measured by $\sfd_0(\mu, \nu)$.}
	\label{fig:stab_prior_learning}
\end{figure}

When $\sfd = \sfW_1$, the Wasserstein-$1$ distance~\cite{villani2009optimal}, it holds that
\begin{align}\label{eqn:wass_prior_stab}
	\sfW_1(\mu^y, \nu^y) \leq \frac{C_p(\mu, \nu; y)}{\min(Z_\mu^y, Z_\nu^y)} \sfW_p(\mu, \nu)
\end{align}
for any $p \in [1, \infty]$; see~\cite[Lem.~5.1, p.~38]{trevisan2023wide}. The constant $C_p(\mu, \nu; y)$ depends on the $p/(p-1)$-th moment of the priors and is known explicitly~\cite{trevisan2023wide}. In particular, the limit $p \to 1$ requires one of the priors to have compact support, which excludes popular unbounded priors such as Gaussian measures. As a result, the scenario $\sfd = \sfd_0 = \sfW_1$ is relatively uncommon. 

A result similar to~\eqref{eqn:wass_prior_stab} can be derived using more sophisticated techniques, particularly in the case $p = 2$~\cite{garbuno2023bayesian}. Both proofs rely on Kantorovich--Rubinstein duality formulas for $\sfW_1$. Thus, it remains an open problem to establish stability estimates for the prior-to-posterior map using the stronger Wasserstein-$2$ distance, $\sfd = \sfW_2$. There has been some recent progress addressing this issue for the data-to-posterior map~\cite{dolera2023lipschitz}. Yet, stability bounds in the Wasserstein-$2$ distance for the prior-to-posterior map would hold particular relevance for uncertainty quantification because $\sfW_2(\mu^y, \nu^y)$ also controls the distance between the posterior covariance operators.

Although $\sfd = \sfd_0 = \sfW_1$ is challenging to satisfy, there are many settings in which we nonetheless wish for $\sfd = \sfd_0$ in \eqref{eqn:prior_perturb_general_bound} to hold. For example, this assumption is especially useful when applying the stability estimate iteratively to a time-indexed sequence of Bayesian inverse problems~\cite{law2015data,sanz2023inverse}. By relaxing $\sfW_1$ to the \emph{Dudley metric}, also referred to as the bounded Lipschitz distance, we show that $\sfd = \sfd_0$ becomes achievable under minimal assumptions.

We now adopt the compact linear functional notation $\pi(f) \defeq \E_{u \sim \pi} [f(u)]$, which represents the integral of a function $f$ with respect to a measure $\pi$.
The bounded Lipschitz distance is an integral probability metric defined as
\begin{align}\label{eqn:dudley_def}
	\sfd_{\mathrm{BL}}(\mu, \nu) \defeq \sup_{\norm{f}_{\mathrm{BL}} \leq 1} \abs{\mu(f) - \nu(f)}\,, \ \text{ where } \ \norm{f}_{\mathrm{BL}} \defeq \max\bigl(\norm{f}_{L^\infty}, \Lip(f)\bigr),
\end{align}
for any probability measures $\mu$ and $\nu$. We recall that $\Lip(f)=\sup_{u\neq v}\abs{f(u)-f(v)}/\norm{u-v}$. The distance $\sfd_{\mathrm{BL}}$ metrizes the topology of weak convergence~\cite{villani2009optimal}. Local Lipschitz stability in Dudley's metric has been established in previous work~\cite{basu1998stability,basu2000uniform}. However, we present a more quantitative version tailored for inverse problems that is inspired by prior work~\cite[Remark~5.2, p.~39]{trevisan2023wide}. 

For brevity, we define $\sfL^y(u) \defeq \exp(-\Phi(u; y))$ from \eqref{eqn:posterior}. We now state the following local Lipschitz continuity result.
\begin{proposition}[prior-to-posterior local Lipschitz stability in Dudley metric]\label{prop:bl_stability}
	Fix $y$ and let $\mu^y$ and $\nu^y$ be posteriors corresponding to priors $\mu$ and $\nu$, respectively, according to \eqref{eqn:posterior}. If the likelihood function $\sfL^y$ is globally Lipschitz continuous and there exists a constant $b^y>0$ such that $0\leq \sfL^y\leq b^y$ everywhere, then
	\begin{align}\label{eqn:bl_stability}
		\sfd_{\mathrm{BL}}(\mu^y,\nu^y)\leq \left(\frac{2b^y + 2\Lip(\sfL^y)}{\max\bigl(Z_\mu^y, Z_\nu^y\bigr)}\right)\, \sfd_{\mathrm{BL}}(\mu,\nu)\,.
	\end{align}
\end{proposition}
We provide the proof of~\Cref{prop:bl_stability} below, which is prototypical of how other basic results of the prior-to-posterior type \eqref{eqn:prior_perturb_general_bound} are established (see, e.g., \cite[Thm. 4.6, p. 43]{bach2024inverse} or \cite[Thm. 8, p. 11]{sprungk2020local}).
\begin{proof}
	Except where required, we ignore all dependence on $y$ for notational convenience. Let $f$ satisfy $\norm{f}_{\mathrm{BL}}\leq 1$. By the definition \eqref{eqn:posterior} of posterior,
	\begin{align*}
		\mu^y(f) - \nu^y(f)&=\frac{\mu(f\sfL)}{\mu(\sfL)} - \frac{\nu(f\sfL)}{\nu(\sfL)}\\
		&=\frac{\mu(f\sfL) - \nu(f\sfL)}{\mu(\sfL)} + \nu^y(f)\left( \frac{\nu(\sfL)-\mu(\sfL)}{\mu(\sfL)}\right)\\
		&\leq \frac{1}{\mu(\sfL)}\Bigl[\abs{\mu(f\sfL) - \nu(f\sfL)} + \abs{\mu(\sfL) - \nu(\sfL)}\Bigr]
	\end{align*}
	because $\abs{f}\leq 1$ everywhere and $\nu^y$ is a probability measure. Next, since $\norm{f}_{\mathrm{BL}}\leq 1$, it holds that $\abs{f\sfL}\leq b$ and
	\begin{align*}
		\abs{f(u)\sfL(u)-f(v)\sfL(v)}&\leq \abs{f(u)\sfL(u)-f(u)\sfL(v)} + \abs{f(u)\sfL(v) - f(v)\sfL(v)}\\
		&\leq \abs{\sfL(u)-\sfL(v)} + b\abs{f(u)-f(v)}\\
		&\leq \bigl(\Lip(\sfL) + b\bigr)\norm{u-v}\,.
	\end{align*}
	Thus, $\norm{f\sfL}_{\mathrm{BL}}\leq \max(b, b + \Lip(\sfL))= b + \Lip(\sfL)$. We also have $\norm{\sfL}_{\mathrm{BL}}\leq \max(b, \Lip(\sfL))$. We deduce using \eqref{eqn:dudley_def} that
	\begin{align*}
		\abs{\mu^y(f) - \nu^y(f)}&\leq \frac{\bigl(b + \Lip(\sfL)\bigr)\sfd_{\mathrm{BL}}(\mu,\nu) + \max\bigl(b, \Lip(\sfL)\bigr)\sfd_{\mathrm{BL}}(\mu,\nu)}{\mu(\sfL)}\\
		&\leq \frac{2b + 2\Lip(\sfL)}{\mu(\sfL)} \sfd_{\mathrm{BL}}(\mu,\nu)
	\end{align*}
	uniformly in $f$. Since $\mu(\sfL)=Z_\mu^y$, the assertion \eqref{eqn:bl_stability} follows by symmetry.
\end{proof}
\begin{remark}
	The hypotheses of \cref{prop:bl_stability} are quite mild. For example, if $\cY$ is finite-dimensional, then a Gaussian likelihood satisfies the assumptions on $\sfL^y$. Moreover, it is easy to generalize the proof of \cref{prop:bl_stability} if the likelihood functional $\sfL^y$ is only H\"older continuous instead of Lipschitz continuous.
\end{remark}

Weak metrics such as $\sfd_{\mathrm{BL}}$ and $\sfW_p$ are particularly well-suited for data-driven prior learning because these distances remain meaningful even when the two input distributions are mutually singular. In contrast, statistical divergences such as TV, KL, and Hellinger attain their maximum value when the measures are mutually singular. This issue often arises in practice because a finite training dataset can lead to an approximate prior that is mutually singular with respect to the reference prior. For instance, using the empirical measure as an estimator results in discrete mismatched support, while employing the closest Gaussian approximation to the reference prior yields a degenerate covariance operator. 
More fundamentally, in infinite-dimensional spaces, probability measures tend to be mutually singular. As an illustrative example, we refer to the Feldman--Hajek equivalence theorem for Gaussian measures~\cite[Thm.~6.13, p.~531]{stuart2010inverse}. 

Results such as \cref{prop:bl_stability} provide \emph{consistency} in the following sense. Suppose that the observed data samples are drawn from a ``true'' prior and that we can consistently estimate this prior using a data-driven method. Consequently, \cref{prop:bl_stability} shows that we can achieve consistency in the resulting posteriors, even if the estimated prior is mutually singular with respect to the true prior.

\paragraph*{Insights from linear regularization analysis}
Theoretical analyses of regularizer learning for linear inverse problems in infinite dimensions have recently gained attention, with notable developments presented in \cite{alberti2021learning,alberti2024learning,burger2024learned,kabri2024convergent,ratti2024learned}. This line of research is more naturally situated within the nonparametric statistical setting described in \cref{sec:back_ip_prob_stat}, as opposed to the functional-analytic viewpoint of \cref{sec:back_ip_reg} or the nonparametric Bayesian perspective in \cref{sec:back_ip_prob_bayes}. However, the research in this direction can still be related to Bayesian inference by interpreting variational point estimators as either the mean or MAP point of an underlying posterior distribution. Such a connection will be made clear in the following discussions. Regardless of viewpoint, the probabilistic models for the parameters and observed data in the inverse problem play a major role in the theoretical insights that are unveiled, as we now explain.

To fix the notation, let $\cU$ and $\cY$ be real separable Hilbert spaces. Consider the inverse problem of finding an unknown parameter $u\in\cU$ from observed data $y$, where $u$ and $y$ are related by
\begin{align}\label{eqn:prior_linear_ip}
	y=\cA u+\eta\,.
\end{align}
Here, $\eta$ denotes additive noise corrupting the observation, as in \eqref{eqn:ip_main}, except now the forward operator $\cA\colon\cU\to\cY$ is assumed to be linear, continuous, injective, and \emph{fully known}. Now departing from the setting of \cref{sec:ee_theory_linear}, we provide a description of the statistical structure of \eqref{eqn:prior_linear_ip}. We model $u$ as a $\cU$-valued random vector with mean and covariance denoted by
\begin{align}\label{eqn:mean_and_cov}
	\overline{u}\defeq \E u\qa \Sigma\defeq \E\bigl[(u-\overline{u})\otimes_{\cU} (u-\overline{u})\bigr]\,,
\end{align}
respectively. Similarly, we model $\eta$ as a zero-mean $\cY$-indexed stochastic process with a continuous covariance operator $\Gamma\colon\cY\to\cY$, independent of $u$. We assume that $\Sigma$ and $\Gamma$ are both strictly positive definite for simplicity.

One can associate to \eqref{eqn:prior_linear_ip} the \emph{generalized Tikhonov reconstruction operator} $\cR_{h,B}^\Gamma\colon\cY\to\cU$ defined by
\begin{align}\label{eqn:tik_operator}
	y\mapsto \cR_{h,B}^\Gamma(y)\defeq \argmin_{u\in\cU} \left\{\sfd_\Gamma^2(y,\cA u) + \norm{B^{-1}(u-h)}_{\cU}^2\right\}\,.
\end{align}
The element $h\in\cU$ is a centering parameter for $u$. The continuous linear operator $B\in\cL(\cU;\cU)$ penalizes the smoothness of $u$. We interpret the second term on the right-hand side of \eqref{eqn:tik_operator} as being equal to $+\infty$ if $(u-h) \notin \image(B)$. The first term can be interpreted formally as $\norm{\Gamma^{-1/2}(y - \cA u)}_{\cY}^2$, but such an expression is infinite almost surely if $\cY$ is infinite-dimensional. To address this, the actual discrepancy function $\sfd_\Gamma$ captures the technical modifications needed to ensure that \eqref{eqn:tik_operator} is well-defined~\cite[Sec.~2]{alberti2021learning}. We emphasize that $\sfd_\Gamma$ depends on the covariance operator $\Gamma$ of the noise, which is assumed known. Thus, $\cR_{h,B}^\Gamma$ itself depends on the inherent noise level in the observations.

Due to the possibility of a nonzero $h$, the map $\cR_{h,B}^\Gamma$ is affine rather than linear. Its closed-form expression is provided in~\cite[Sec.~2, p.~5]{alberti2021learning}. The associated calculations are similar to those used to establish that the posterior distribution for a linear Gaussian Bayesian inverse problem in an infinite-dimensional setting (for both data and parameters) is well-defined~\cite{knapik2011bayesian}. However, it is important to emphasize that no Gaussian assumptions are made in the current context.

Alberti et al.~\cite{alberti2021learning} rigorously characterize the average-case optimal generalized Tikhonov regularizer in terms of the covariance structure assumed on the unknown parameter $u$. To this end, let $(h^\star, B^\star)$ be the optimal pair minimizing the expected squared reconstruction error
\begin{align}\label{eqn:tik_mse}
	\E_{(u,y)\sim\nu} \norm[\big]{u-\cR_{h,B}^\Gamma(y)}^2_{\cU}\,.
\end{align}
The expectation is taken with respect to the joint distribution $\nu\defeq \Law(u,y)$ of the random variable $(u,y)$ under \eqref{eqn:prior_linear_ip}.
Under minor technical assumptions, the unique global minimizer of \eqref{eqn:tik_mse} is $(h^\star, B^\star)=(\overline{u}, \Sigma^{1/2})$ \cite[Thm.~3.1, p.~6]{alberti2021learning}. Notably, even though \eqref{eqn:tik_mse} depends on the full joint distribution $\nu$---and hence the forward map $\cA$---the optimal solution \emph{only depends} on the first two moments of the marginal law of $u$.
The technical contributions leading to this result are substantial, as they require careful analysis to handle the infinite-dimensional noise process $\eta$ and regularized inverses of compact operators.

It is also possible to develop finite-sample estimators for the problem \eqref{eqn:tik_mse} if only a dataset of i.i.d.~pairs $\{(u_n,y_n)\}_{n=1}^N$ from $\nu$ is available. Because the optimal Tikhonov regularizer $(h^\star, B^\star)$ is known to be the mean and square root covariance of the random variable $u$, a natural \emph{unsupervised} strategy is first to compute the empirical mean $\widehat{u}$ and empirical covariance $\widehat{\Sigma}$ from the training data, then define the plug-in estimator $(\widehat{h}, \widehat{B})\defeq (\widehat{u}, \widehat{\Sigma}^{1/2})$.
Assuming that $u$ is strongly-subgaussian, application of tools from covariance estimation \cite[Thm.~4.2, p.~8, and p.~29]{alberti2021learning} deliver the nonasymptotic error bound
\begin{align}\label{eqn:tik_excess_risk}
	\E \Bigl[\E_{(u,y)\sim\nu} \norm[\big]{u-R^\Gamma_{\widehat{h},\widehat{B}}(y)}^2_{\cU}\Bigr]\leq \E_{(u,y)\sim\nu} \norm[\big]{u-R^\Gamma_{h^\star,B^\star}(y)}^2_{\cU} + \frac{C}{\sqrt{N}}
\end{align}
for some $C>0$. The first expectation on the left-hand side of \eqref{eqn:tik_excess_risk} is taken with respect to the randomness in the finite sample set $\{(u_n,y_n)\}_{n=1}^N$. The $C/\sqrt{N}$ term corresponds to a ``slow'' nonparametric rate for the excess error because a fast rate for squared error would instead approach $1/N$. It would be interesting to explore the optimality of the derived rate in \eqref{eqn:tik_excess_risk}.

Alternatively, a fully supervised training strategy based on minimizing an empirical approximation to \eqref{eqn:tik_mse} gives similar convergence guarantees in theory~\cite[Thm.~4.1, p.~7]{alberti2021learning}, but empirically demonstrates worse performance than the unsupervised approach. However, suppose the Tikhonov regularizer is replaced by a more general regularizer, such as a sparsity-promoting one~\cite{alberti2024learning}. In that case, the supervised learning strategy becomes necessary because closed-form expressions for the optimal $(h^\star, B^\star)$ are no longer available. The resulting finite-sample analysis also becomes more difficult, requiring advanced tools from statistical learning theory~\cite{ratti2024learned}.

Although the preceding discussion characterizes the optimal Tikhonov regularizer $(h^\star, B^\star)$ from a linear operator learning perspective, it does not comment on the \emph{accuracy} of the corresponding reconstruction map $\cR_{h^\star,B^\star}^\Gamma$. This requires understanding under what conditions $\norm{u-R^\Gamma_{h^\star, B^\star}(\cA u + \eta)}_{\cU}\to 0$ as $\Gamma\to 0$ in an appropriate sense. Combined with \eqref{eqn:tik_excess_risk}, such a result would also imply convergence for the estimator $\cR_{\widehat{h},\widehat{B}}^\Gamma$. To carry out this consistency analysis, it is instructive to work in a setup in which $\cA^*\cA$ and $\Sigma$, together with $\cA\cA^*$ and $\Gamma$, are simultaneously diagonalizable. This should be compared to a similar approach used to develop sharp convergence rates for forward operator learning~\cite{de2023convergence}.

At the center of the analysis is the SVD of the compact injective forward operator $\cA$, which we expand as
\begin{align}\label{eqn:svd}
	\cA=\sum_{j=1}^\infty \sigma_j \, e_j\otimes_\cU\phi_j\,.
\end{align}
The singular values $\{\sigma_j\}_{j\in\N}\subset \R_{>0}$ of $\cA$ are ordered to be nonincreasing and assumed to each have multiplicity one. The $\{\phi_j\}_{j\in\N}$ form an ONB of $\cU$ and $\{e_j\}_{j\in\N}$ form an ONB of $\cY$. Using these SVD bases, we further assume that
\begin{align}\label{eqn:cov_diag}
	\Sigma=\sum_{j=1}^\infty \varsigma_j\, \phi_j\otimes_\cU\phi_j\qa \Gamma=\delta^2\sum_{j=1}^\infty \gamma_j \, e_j\otimes_\cY e_j\,,
\end{align}
which delivers simultaneous diagonalizability with the normal operators corresponding to $\cA$. The eigenvalues $\{\varsigma_j\}_{j\in\N}$ and $\{\delta^2 \gamma_j\}_{j\in\N}$ from \eqref{eqn:cov_diag} are nonincreasing sequences. We assume that $\gamma_1=1$ \cite[Sec.~2.2]{burger2024learned}. The parameter $\delta>0$ in \eqref{eqn:cov_diag} serves as a noise level. It will be convenient to define $\xi\defeq \eta/\delta$. The covariance operator of $\xi$ is $\sum_{j=1}^\infty \gamma_j\, e_j\otimes_\cY e_j$. We also assume that $\E u =\overline{u}=0$.

Based on the SVD from \eqref{eqn:svd}, \emph{spectral reconstruction operators} regularize the pseudoinverse of $\cA$ by only approximately inverting the singular values. Define the spectral reconstruction operator $\cR_\psi\colon\cY\to\cU$ by
\begin{align}\label{eqn:svd_reg}
	\cR_\psi\defeq \sum_{j=1}^\infty \psi_j \,\phi_j\otimes_\cY e_j\,,
\end{align}
where $\psi\defeq\{\psi_j\}_{j\in\N}\subset\R_{\geq 0}$ is a sequence of filtering coefficients---allowed to depend on $\cA$---such that $\psi_j\approx 1/\sigma_j$ \cite{burger2024learned}. One can then ask what the optimal choice of $\psi$ is with respect to the expected squared error \eqref{eqn:tik_mse} \cite{kabri2024convergent}.

To answer this question, we adopt the common approach of diagonalizing the equation \eqref{eqn:prior_linear_ip} into a sequence space model in the SVD bases of $\cA$ \cite{cavalier2008nonparametric,knapik2011bayesian,de2023convergence}. This delivers the sequence of decoupled scalar linear inverse problems
\begin{align}\label{eqn:seq_space}
	\ip{y}{e_j}=\sigma_j\ip{u}{\phi_j} + \delta\ip{\xi}{e_j}\qfa j\in\N\,.
\end{align}
Due to \eqref{eqn:seq_space} and the definitions in \eqref{eqn:cov_diag}, it holds that
\begin{align}\label{eqn:seq_min}
	\E_{(u,y)\sim\nu} \norm{u-\cR_{\psi}y}^2_{\cU}=\sum_{j=1}^\infty\Bigl[(1-\sigma_j\psi_j)^2\varsigma_j + \delta^2\gamma_j\psi_j^2\Bigr] \,.
\end{align}
Thus, the minimization of \eqref{eqn:seq_min} over the whole sequence $\{\psi_j\}_{j\in\N}$ decouples into scalar minimization problems for each individual $\psi_j$. The minimizers are
\begin{align}\label{eqn:seq_optimal}
	\psi_j(\delta)\defeq \frac{\sigma_j}{\sigma_j^2 + \frac{\delta^2\gamma_j}{\varsigma_j}}\qfa j\in\N\,,
\end{align}
which is seen by explicitly solving the quadratic program. Writing $\psi(\delta)\defeq\{\psi_j(\delta)\}_{j\in\N}$, the optimal spectral regularizer is then $\cR_{\psi(\delta)}$.

Under the simultaneous diagonalizability assumption \eqref{eqn:cov_diag} and the condition $\overline{u}=h^\star=0$, by comparing to \cite[Eqn.~(12), p.~5]{alberti2021learning} when written in the SVD bases \eqref{eqn:svd} of $\cA$, we find that $\cR_{\psi(\delta)}=\cR_{h^\star,B^\star}^\Gamma$ \cite[Remark~2, p.~49]{burger2024learned}. That is, the optimal spectral regularizer---which did not enforce an explicit smoothness penalization---\emph{is identical to} the optimal Tikhonov regularizer, which does enforce a quadratic penalty. It is interesting that Tikhonov structure \emph{emerges} in the optimal spectral reconstruction as a consequence of the additive noise model \eqref{eqn:prior_linear_ip} and the expected squared reconstruction error objective function \eqref{eqn:seq_min}. This finding is consistent with known results for finite-dimensional nonlinear regression~\cite{bishop1995training,webb1994functional}.

A heuristic explanation for this phenomenon in the current linear infinite-dimensional setup is as follows. The optimal generalized Tikhonov reconstruction \eqref{eqn:tik_operator} implicitly computes the MAP point of a Gaussian posterior $\mu^y$ determined by a Gaussian likelihood $\normal(\cA u, \Gamma)$ and Gaussian prior $\normal(0, \Sigma)$, which match the first two moments of $\Law(y\condbar u)$ and $\Law(u)$, respectively. This immediately shows that the optimal Tikhonov reconstruction is the posterior mean $\E_{u\sim \mu^y}[u]$ of $\mu^y$, agreeing with known formulas \cite[Prop.~3.1, p.~2630]{knapik2011bayesian}. To relate back to the optimal spectral regularizer, one would have to argue that since $\cR_\psi$ is linear in $\psi$ and $y$ is linearly related to $u$ in expectation, the minimizer of \eqref{eqn:seq_min} over $\psi$ (i.e., linear maps) should equal the Gaussian posterior mean instead of the Bayes estimator $\E[u\condbar y]$. However, a proper proof based on this intuition in the present infinite-dimensional setting is beyond the scope of the chapter.

We close this subsection by summarizing accuracy guarantees for the optimal spectral regularizer $\cR_{\psi(\delta)}$. The preceding discussion implies that consistency of learned spectral regularization is equivalent to consistency of the Gaussian posterior mean estimator. Consistency of the latter is well known \cite{knapik2011bayesian,knapik2018general}. Under mild distributional assumptions, more direct arguments \cite[Thm.~1, p.~1351]{kabri2024convergent} establish the pointwise average-case convergence
\begin{align}\label{eqn:reg_avg_converge}
	\lim_{\delta\to 0} \E_{\xi\sim\Law(\xi)} \norm[\big]{\cA^{-1}y- \cR_{\psi(\delta)}(y+\delta \xi)}_{\cU}^2=0 \qfa y\in\image(\cA)\,.
\end{align}
Compare this result with \eqref{eqn:reg_prog_converge}, which gives pointwise worst-case convergence. The convergence in \eqref{eqn:reg_avg_converge} is actually uniform over certain bounded sets of noise distributions with noise level at most $\delta$ and covariance $\Gamma$; see \cite[Eqn.~(2.8), p.~48]{burger2024learned}. Going beyond pointwise error, we also have the full mean square convergence
\begin{align}\label{eqn:reg_avg_converge_avg}
	\lim_{\delta\to 0} \E_{(u,\xi)\sim\Law(u,\xi)} \norm[\big]{u- \cR_{\psi(\delta)}(\cA u+\delta \xi)}_{\cU}^2=0\,.
\end{align}

By using well-established Bayesian theory, it should also be possible to obtain \emph{quantitative} convergence rates as $\delta \to 0$ under suitable decay conditions on the singular values of all the involved operators \cite{knapik2011bayesian}.
The misspecified case, in which the noise distribution in the training data possibly differs from the one used to compute the squared reconstruction error, has also been handled \cite{burger2024learned}.
In all cases, the diagonal analysis suggests that choosing the training data to be corrupted by white noise (i.e., $\gamma_j= 1$ for all $j$) leads to better convergence properties for solving the inverse problem. This message is consistent with a related insight for forward linear operator learning in the following sense. In~\cite[Corollary~3.6, pp.~17--18]{de2023convergence}, it was observed that training data $\{y_n\}$ corrupted with white noise worsens convergence rates for learning the forward operator. Since the forward problem becomes harder, the inverse problem intuitively becomes easier as claimed.

Beyond the mean square optimal estimator \eqref{eqn:seq_optimal}, the preceding spectral regularization framework has also been combined with plug-and-play estimators~\cite[Sec.~4]{burger2024learned}. We will overview such plug-and-play methods next.

\subsection{Denoising networks and plug-and-play methods}\label{sec:reg_pnp}
A major breakthrough in bridging model-based and learning-based approaches to inverse problems was the realization that powerful denoisers---including neural networks trained to denoise corrupted signals---implicitly encode high-quality priors~\cite{milanfar2025denoising}. Plug-and-play (PnP) methods capitalize on this insight by integrating pre-trained image denoising networks into iterative optimization algorithms as implicit regularizers. Instead of explicitly specifying a prior through a regularization functional $R$, PnP methods replace the proximal step in classical algorithms such as the alternating direction method of multipliers (ADMM) or proximal gradient descent with a learned denoiser $\cD$. In the following, we briefly review the latest developments in PnP. \Cref{sec:reg_pnp_motiv} introduces the PnP framework and \cref{sec:reg_pnp_theory} surveys existing theoretical guarantees for the framework.

\subsubsection{Motivation for PnP}\label{sec:reg_pnp_motiv}
To motivate PnP, consider a composite objective function
\begin{align}
	\cJ(u) = f(u) + R(u)\,,
\end{align}
where $f$ represents a data-fidelity term and $R$ is a possibly non-smooth regularizer. In classical variational approaches, proximal methods play a central role in handling non-smooth regularization. The proximal operator associated with $R$ is defined by
\begin{equation}\label{eq:prox}
	y\mapsto \mathrm{prox}_{R}(y) = \arg\min_{u\in\cU} \left\{\frac{1}{2}\|y - u\|_{\cU}^2 + R(u) \right\}\,.
\end{equation}
Comparing \eqref{eq:prox} to \eqref{eqn:reg_opt}, we observe that $\mathrm{prox}_{R}(y)$ precisely computes the variational regularization solution of a linear inverse problem with identity forward map, i.e., denoising, under a least-squares misfit and regularizer $R$. The proximal operation also has a natural Bayesian interpretation: it computes the MAP estimate under a Gaussian likelihood centered at $y$ and a prior distribution with density proportional to $\exp(-R)$.

In proximal gradient descent (PGD), the iteration for minimizing $\cJ$ takes the form
\begin{align}
	u_{k+1} = \mathrm{prox}_{R}\bigl(u_k - s_k \nabla f(u_k)\bigr)
\end{align}
for each $k$, where $s_k > 0$ is a step size. PnP methods modify this framework by replacing the proximal operator with a denoising operator $\cD\colon\cU\to\cU$, which implicitly encodes prior knowledge. That is, the update rule becomes
\begin{equation}\label{eq:PnP_basic}
	u_{k+1} = \cD\bigl(u_k - s_k \nabla f(u_k)\bigr)\,,
\end{equation}
where $\cD$ may be a pre-trained convolutional neural network (CNN), neural operator, or any other denoiser. In this formulation, $\cD$ implicitly serves as a \emph{learned proximal operator} that regularizes the solution to lie on a manifold of natural images or signals learned from the training data.

The original formulation of PnP was introduced by Venkatakrishnan et al.~\cite{venkatakrishnan2013plug} in the context of ADMM and using the classical BM3D denoiser. This idea was later extended to incorporate deep learning-based denoisers, including CNNs~\cite{meinhardt2017learning}. More generally, each PnP iteration computes the fixed point of an operator that combines the forward model---appearing in the data-fidelity term $f$---with the denoiser, which acts as a prior. This has close connections to deep equilibrium architectures that are also based on fixed-point iterations \cite{gilton2021deep}.

In the language of operator learning, the denoiser is viewed as a learned nonlinear map $\cD\colon \cU \to \cU$ that is trained to model the structure of realistic signals. For example, in the case of a linear forward map $\mathcal{G}(u) = \mathcal{A}u$ and quadratic data fidelity $f(u) \defeq \frac{1}{2}\|\mathcal{A}u - y\|_\cY^2$, the PnP iteration becomes
\begin{align}
	u_{k+1} = \cD\bigl(u_k - s_k \mathcal{A}^* (\mathcal{A}u_k  - y)\bigr)\,.
\end{align}
This formulation embeds a learned operator into an iterative algorithm that solves an inverse problem governed by the forward operator $\cA$.

Subsequent developments in the PnP framework have focused on theoretical analysis (including fixed-point convergence), stability, and the implicit regularization properties of the denoiser. These investigations aim to understand the conditions under which such iterations are well-posed and whether they yield solutions consistent with the underlying inverse problem in the limit of vanishing noise.

\subsubsection{Theoretical analysis of PnP}\label{sec:reg_pnp_theory}
The theoretical understanding of PnP methods has advanced considerably. There is a growing focus on identifying conditions under which a learned denoiser $\cD$ can be reliably embedded into iterative optimization or sampling algorithms for solving inverse problems. Broadly, these developments can be categorized into four areas: convergence analysis, proximal structure, probabilistic interpretations, and control of regularization strength.

\paragraph*{Convergence via nonexpansiveness and monotonicity}
A foundational line of work views the denoiser $\cD$ as a nonexpansive Lipschitz-continuous operator. Ryu et al.~\cite{ryu2019plug} proves that if $\cD$ is strictly nonexpansive, i.e., $\Lip(\cD)< 1$,
then PnP-ADMM and PnP-FBS (forward-backward splitting) converge to a fixed point without requiring vanishing step sizes. The key idea is that the overall update becomes a contraction under suitable step-size conditions. Belkouchi et al.~\cite{belkouchi2025learning} employ monotone operator theory to analyze PnP iterations. They show that if $\cD$ is the resolvent of a monotone operator, then fixed-point iterations converge under Tseng's splitting scheme. This analysis links learned operators to classical monotone inclusions, which broadens the class of denoisers for which convergence guarantees hold.

\paragraph*{Learned proximal structure}
Another direction seeks to design denoisers that are explicit proximal operators. Fang et al.~\cite{fang2024whats} introduce Learned Proximal Networks (LPNs), in which $\cD$ is trained to be the gradient of a convex function, i.e., $\cD = \nabla \phi$ for some convex $\phi$. This implies that $\cD=\mathrm{prox}_R$ is a proximal map of a convex regularizer $R$. It further allows PnP iterations to be interpreted as classical PGD or ADMM updates whose convergence is inherited from convex optimization theory. Pesquet et al.~\cite{pesquet2021learning} propose learning resolvent operators directly. In this setting, a universal approximation theorem holds: any nonexpansive neural network can approximate the resolvent of a maximally monotone operator. This result connects learned denoisers to variational analysis by showing that even when the underlying regularizer is unknown, a denoiser with resolvent structure ensures convergence and interpretability.

\paragraph*{Probabilistic and sampling-based extensions}
A complementary line of work reframes PnP as a probabilistic inference method. Bouman et al.~\cite{bouman2023generative} introduce Generative PnP (GPnP), which transforms deterministic PnP iterations into a Markov chain by replacing each update with a stochastic proximal generator. 
The resulting dynamics sample from the Bayesian posterior distribution corresponding to the prior specified implicitly by the denoiser, and thus can provide both point estimates and uncertainty
quantification.
Wu et al.~\cite{wu2024principled} demonstrate that any inverse problem can be reformulated as a Gaussian image denoising task. This allows one to embed pretrained diffusion models within lightweight MCMC schemes, thereby integrating score-based generative models into the PnP framework. In this setting, denoisers serve as approximations to the score function. They enable posterior sampling via Langevin dynamics or stochastic proximal steps.

\paragraph*{Regularization control and stability}
Recent theoretical contributions establish that PnP methods can serve as convergent regularization schemes in the classical sense: reconstructions remain stable in the presence of noise and converge to the true solution of the inverse problem as the noise level decreases~\cite{engl1996regularization,benning2018modern}. Ebner and Haltmeier~\cite{ebner2024plug} prove this result for a family of PnP methods with denoisers tailored to the noise level. Hauptmann et al.~\cite{hauptmann2024convergent} provide the first complete regularization theory for PnP with linear denoisers. By introducing a tunable spectral filter into the denoiser, it is shown that the resulting iterations form a Tikhonov-type regularization scheme. As the measurement noise $\delta$ goes to zero, the iterates converge to the exact minimizer of a well-defined variational problem and classical convergence rates are recovered under realistic assumptions. Khelifa et al.~\cite{khelifa2025enhanced} propose Tweedie scaling, which introduces a single tunable parameter to rescale any pretrained denoiser. This parameter allows precise control over the strength of the implicit regularization, guarantees convergence of PnP-PGD, and improves empirical performance on tasks such as inpainting and denoising. The Tweedie scaling approach provides a systematic mechanism for adjusting the influence of the prior.

%%%%%%%%%%%%%%%%%%%%%%%%%%%%%%%%%%%%%%%%%%%%%%%%%%%%%%%%%%%%%%%%%%%%%%%%%%%%%%%%%%%%%%%%%%%%%%%%%%%%%%%%%%%%%%%%%%%%%%%%%%%%%%%%%%%%%%%%%%%%%%%%
\section{Conclusion and outlook}\label{sec:conclusion}
This chapter surveys recent advances in the use of operator learning as a tool to solve inverse problems. The intersection of these two fields is particularly powerful when unknown parameters are infinite-dimensional and the forward map requires solving a partial differential equation. Two complementary paradigms emerge in this subject area. The first is an end-to-end approach that uses operator learning to directly invert for unknown parameters or Bayesian posteriors given new observations. This strategy replaces traditional inverse problem solvers with neural operator models that are trained on noisy data-parameter pairs. The most popular operator learning architectures enjoy theoretical approximation guarantees for a few specific inverse problems. New architectures that are more adapted to the ill-posed character of inverse problems have been proposed, but are less established theoretically. The second framework is more modular. Instead of emulating the entire inverse map, this approach invokes operator learning to selectively estimate priors, regularizers, or denoisers within conventional inversion algorithms. This results in a principled blending of data-driven components with structure encoded by the forward map. Consequently, the nature of the regularization introduced by these hybrid data-driven methods is better understood than that of fully learned inverse maps. By adopting a probabilistic perspective, the chapter unifies both paradigms and invites further mathematical development in this rapidly evolving field.

A distinctive feature of this chapter is its emphasis on measure-centric formulations of inverse problems and their compatibility with operator learning. Although the chapter takes a first step in merging both ideas, for example, by connecting DeepSets to Neural Inverse Operators for inverse boundary value problems, there is still a need for basic theoretical development. Core well-posedness questions regarding existence, uniqueness, and stability largely remain unanswered for inverse problems in the space of probability measures. Novel distribution-based operator learning architectures lack guarantees when applied to measure-valued data. Just as importantly, the deployment of measure-based operator learning in scientific inverse problems involving noisy, real-world experimental data would greatly elevate the impact and relevance of the field.

Sequential Bayesian inverse problems in data assimilation \cite{law2015data,sanz2023inverse} are an application where the probabilistic and measure-centric principles advocated by this chapter are highly relevant. Existing accurate algorithms are hindered in practice by the computational expense of repeatedly applying Bayes' rule to update prior to posterior. Recent work shows the promise of combining operator learning with measure transport to capture complex mean-field distribution dependence in filtering problems~\cite{al2025fast,bach2025learning,zhao2024functional,bach2024learning,bach2024inverse,zhou2024bi}. Despite this encouraging progress toward overcoming the computational barriers associated with sequential Bayesian inference in dynamical systems, the underlying mathematical theory and design of measure-based architectures demand further investigation.

Recent cautionary work draws attention to adversarial instabilities, poor interpretability, and overfitting in learned inverse solvers, particularly in high-dimensional imaging tasks~\cite{colbrook2022difficulty,antun2020instabilities,boche2023limitations,gottschling2025troublesome,evangelista2025or}. These limitations underscore the need for theory that clarifies when and why operator learning can deliver both accurate and stable solutions to ill-posed problems. To be more impactful, such a theory will have to go beyond the current approach of performing a case-by-case analysis for each individual inverse problem. As advocated by \cite{burger2024learning}, methodological innovation in out-of-distribution detection and uncertainty quantification (UQ) will bolster trust of data-driven inversion methods. Regarding the former, preliminary work shows that a careful design of the training data distribution can make learned maps robust to distribution shifts~\cite{guerra2025learning}. For the latter, probabilistic inversion methods already deliver UQ. However, the modeling uncertainty inherent in learned inverse maps or regularizers must also be accounted for in the UQ framework. One possible way to achieve this is to use the uncertainty from operator learning to improve the quality of Bayesian inversion~\cite{cleary2021calibrate,dunbar2025hyperparameter}.

%%%%%%%%%%%%%%%%%%%%%%%%%%%%%%%%%%%%%%%%%%%%%%%%%%%%%%%%%%%
\backmatter

\bmhead*{Funding.}
N.H.N. is supported by the U.S. National Science Foundation (NSF) under award DMS-2402036 and by a Klarman Fellowship through Cornell University's College of Arts \& Sciences. Y.Y. is partially supported by NSF award DMS-2409855.

\bmhead*{Acknowledgments.}
The authors thank Maarten de Hoop, Bamdad Hosseini, Matti Lassas, Youssef Marzouk, Richard Nickl, and Andrew Stuart for numerous discussions that helped to shape this chapter. The authors are also grateful to Nicolas Guerra and two anonymous referees for helpful comments on a previous version of the manuscript.

\bibliography{references}% common bib file

%%\input _main.bbl

\end{document}

%% file: macros.tex
\usepackage{lipsum}
\usepackage{algorithmic}
\usepackage{amsopn}

\usepackage{subcaption} % for multiple subfigures
\usepackage{mathtools} % to DeclarePairedDelimiter for ip, abs, norm

\usepackage{tikz}
\usetikzlibrary{positioning}
\usetikzlibrary{arrows.meta, decorations.pathreplacing}

\usepackage{bbm} % for indicator function
\usepackage{stmaryrd} % for ``mapsfrom''
\let\originalleft\left % Fixes \left and \right spacing issues. See discussion at http://tex.stackexchange.com/questions/2607/spacing-around-left-and-right
\let\originalright\right
\renewcommand{\left}{\mathopen{}\mathclose\bgroup\originalleft}
\renewcommand{\right}{\aftergroup\egroup\originalright}
\usepackage{array}
\newcolumntype{M}[1]{>{\centering\arraybackslash}m{#1}} % w/ horizontal centering
\usepackage{tcolorbox}

\usepackage[nameinlink,capitalize]{cleveref}
\crefname{assumption}{Assumption}{Assumptions}
\crefname{notation}{Notation}{Notation}
\crefname{problem}{Problem}{Problems}
\crefname{hypothesis}{Hypothesis}{Hypotheses}
\crefname{condition}{Condition}{Conditions}
\crefname{fact}{Fact}{Facts}

\newcommand{\N}{\mathbb{N}}
\newcommand{\Z}{\mathbb{Z}}
\newcommand{\R}{\mathbb{R}}
\newcommand{\T}{\mathbb{T}}
\newcommand{\C}{\mathbb{C}}
\newcommand{\Q}{\mathbb{Q}}

\newcommand{\cP}{\mathcal{P}}
\newcommand{\cA}{\mathcal{A}}
\newcommand{\cL}{\mathcal{L}}
\newcommand{\cF}{\mathcal{F}}
\newcommand{\cW}{\mathcal{W}}
\newcommand{\cH}{\mathcal{H}}
\newcommand{\cX}{\mathcal{X}}

\newcommand{\cY}{\mathcal{Y}}

\newcommand{\cS}{\mathcal{S}}
\newcommand{\cG}{\mathcal{G}}
\newcommand{\cT}{\mathcal{T}}
\newcommand{\cR}{\mathcal{R}}
\newcommand{\cE}{\mathcal{E}}
\newcommand{\cD}{\mathcal{D}}

\newcommand{\cJ}{\mathcal{J}}
\newcommand{\cK}{\mathcal{K}}
\newcommand{\cU}{\mathcal{U}}
\newcommand{\cO}{\mathcal{O}}
\newcommand{\cV}{\mathcal{V}}
\newcommand{\cZ}{\mathcal{Z}}
\newcommand{\cQ}{\mathcal{Q}}

\newcommand{\sL}{\mathscr{L}}

\newcommand{\sP}{\mathscr{P}}

\newcommand{\sH}{\mathscr{H}}

\newcommand{\sE}{\mathscr{E}}
\newcommand{\sD}{\mathscr{D}}

\newcommand{\sA}{\mathscr{A}}
\newcommand{\sfL}{\mathsf{L}}
\newcommand{\sfE}{\mathsf{E}}
\newcommand{\sfW}{\mathsf{W}}

\newcommand{\sfd}{\mathsf{d}}

\newcommand{\sfT}{\mathsf{T}}

\newcommand{\sfG}{\mathsf{G}}

\newcommand{\sfF}{\mathsf{F}}
\newcommand{\sfK}{\mathsf{K}}

\newcommand{\pOmega}{\partial\Omega}

\DeclarePairedDelimiterX{\iptemp}[2]{\langle}{\rangle}{#1, #2}
\newcommand{\ip}{\iptemp}
\DeclarePairedDelimiterX{\normtemp}[1]{\lVert}{\rVert}{#1}
\newcommand{\norm}{\normtemp}
\DeclarePairedDelimiterX{\abstemp}[1]{\lvert}{\rvert}{#1}
\newcommand{\abs}{\abstemp}
\DeclarePairedDelimiterX{\trtemp}[1]{(}{)}{#1}

\DeclarePairedDelimiterX{\SEtemp}[2]{(}{)}{#1, #2}

\DeclarePairedDelimiterX{\SGtemp}[1]{(}{)}{#1}

\DeclarePairedDelimiterX{\sdivtemp}[2]{(}{)}{%
  #1\;\delimsize\|\;#2%
}

\DeclarePairedDelimiterX{\KLtemp}[2]{(}{)}{%
  #1\;\delimsize\|\;#2%
}
\newcommand{\KL}{\mathsf{KL}\sdivtemp}

\newcommand{\defeq}{\coloneqq} % definition equal in math mode
\newcommand{\eqdef}{\eqqcolon} % equal definition in math mode
\newcommand{\condbar}{\, \vert \,}

\DeclarePairedDelimiterX{\floor}[1]{\lfloor}{\rfloor}{#1} % floor
\DeclarePairedDelimiterX{\ceil}[1]{\lceil}{\rceil}{#1} % ceiling

\DeclareMathOperator{\Id}{Id} % identity operator
\newcommand{\HS}{\mathrm{HS}} % Hilbert--Schmidt
\DeclareMathOperator{\image}{Ran} % image space
\DeclareMathOperator{\domain}{Dom} % image space
\DeclareMathOperator{\Span}{span} % span

\newcommand{\E}{\operatorname{\mathbb{E}}} % expectation
\newcommand{\Cov}{\operatorname{Cov}} % covariance (operator)
\newcommand{\Law}{\operatorname{Law}} % law
\newcommand{\one}{\mathbbm{1}} % indicator function
\newcommand{\diid}{\stackrel{\mathrm{i.i.d.}}{\sim}} % distributed iid tilde symbol
\newcommand{\normal}{\mathcal{N}} % normal distribution
\newcommand{\push}{_{\#}}

\newcommand{\iunit}{\mathrm{i}}
\newcommand{\diff}{d}

\newcommand{\dd}[1]{\,\diff{#1}}

\DeclareMathOperator*{\argmax}{arg\,max}
\DeclareMathOperator*{\argmin}{arg\,min}
\DeclareMathOperator{\Lip}{Lip}

\newcommand{\al}{\alpha}
\newcommand{\ep}{\varepsilon}

\def\qfa{\quad\text{for all}\quad}

\def\qas{\quad\text{as}\quad}
\def\qa{\quad\text{and}\quad}
\def\qf{\quad\text{for}\quad}
\def\qw{\quad\text{where}\quad}
\def\qin{\quad\text{in}\quad}
\def\qon{\quad\text{on}\quad}
\def\qor{\quad\text{or}\quad}

\newcommand{\set}[2]{{\left\{ #1 \,\middle|\, #2 \right\}}}
\newcommand{\slot}{{\,\cdot\,}}

\usepackage{enumitem}
\setlist[enumerate]{leftmargin=.5in}
\setlist[itemize]{leftmargin=.5in}
\makeatletter
\newcommand{\myitem}[1]{%
\item[#1]\protected@edef\@currentlabel{#1}%
}
\makeatother